\def\Resetstrings{%Clears all strings before processing reference listing.
    \def\present{ }\let\bgroup={\let\egroup=}%primitive TeX
    \def\Astr{}\def\astr{}\def\Atest{}\def\atest{}%
    \def\Bstr{}\def\bstr{}\def\Btest{}\def\btest{}%
    \def\Cstr{}\def\cstr{}\def\Ctest{}\def\ctest{}%
    \def\Dstr{}\def\dstr{}\def\Dtest{}\def\dtest{}%
    \def\Estr{}\def\estr{}\def\Etest{}\def\etest{}%
    \def\Fstr{}\def\fstr{}\def\Ftest{}\def\ftest{}%
    \def\Gstr{}\def\gstr{}\def\Gtest{}\def\gtest{}%
    \def\Hstr{}\def\hstr{}\def\Htest{}\def\htest{}%
    \def\Istr{}\def\istr{}\def\Itest{}\def\itest{}%
    \def\Jstr{}\def\jstr{}\def\Jtest{}\def\jtest{}%
    \def\Kstr{}\def\kstr{}\def\Ktest{}\def\ktest{}%
    \def\Lstr{}\def\lstr{}\def\Ltest{}\def\ltest{}%
    \def\Mstr{}\def\mstr{}\def\Mtest{}\def\mtest{}%
    \def\Nstr{}\def\nstr{}\def\Ntest{}\def\ntest{}%
    \def\Ostr{}\def\ostr{}\def\Otest{}\def\otest{}%
    \def\Pstr{}\def\pstr{}\def\Ptest{}\def\ptest{}%
    \def\Qstr{}\def\qstr{}\def\Qtest{}\def\qtest{}%
    \def\Rstr{}\def\rstr{}\def\Rtest{}\def\rtest{}%
    \def\Sstr{}\def\sstr{}\def\Stest{}\def\stest{}%
    \def\Tstr{}\def\tstr{}\def\Ttest{}\def\ttest{}%
    \def\Ustr{}\def\ustr{}\def\Utest{}\def\utest{}%
    \def\Vstr{}\def\vstr{}\def\Vtest{}\def\vtest{}%
    \def\Wstr{}\def\wstr{}\def\Wtest{}\def\wtest{}%
    \def\Xstr{}\def\xstr{}\def\Xtest{}\def\xtest{}%
    \def\Ystr{}\def\ystr{}\def\Ytest{}\def\ytest{}%
}
\Resetstrings

\def\Refformat{%Determines the kind of reference by the presence or
         \if\Jtest\present
             {\if\Vtest\present\journalarticleformat
                  \else\conferencereportformat\fi}
            \else\if\Btest\present\bookarticleformat
               \else\if\Rtest\present\technicalreportformat
                  \else\if\Itest\present\bookformat
                     \else\otherformat\fi\fi\fi\fi}

\def\Rpunct{%Default punctuation control strings if the punctuation
   \def\Lspace{ }%
   \def\Lperiod{ }%  .
   \def\Lcomma{ }%    ,
   \def\Lquest{ }%     ?
   \def\Lcolon{ }%   :
   \def\Lscolon{ }%   ;
   \def\Lbang{ }%      !
   \def\Lquote{ }%   '
   \def\Lqquote{ }%   "
   \def\Lrquote{ }%    `
   \def\Rspace{}%
   \def\Rperiod{.}%  .
   \def\Rcomma{,}%    ,
   \def\Rquest{?}%     ?
   \def\Rcolon{:}%   :
   \def\Rscolon{;}%   ;
   \def\Rbang{!}%      !
   \def\Rquote{'}%   '
   \def\Rqquote{"}%   "
   \def\Rrquote{`}%    `
   }

\def\Lpunct{%Default punctuation control strings if the punctuation
   \def\Lspace{}%
   \def\Lperiod{\unskip.}%  .
   \def\Lcomma{\unskip,}%    ,
   \def\Lquest{\unskip?}%     ?
   \def\Lcolon{\unskip:}%   :
   \def\Lscolon{\unskip;}%   ;
   \def\Lbang{\unskip!}%      !
   \def\Lquote{\unskip'}%   '
   \def\Lqquote{\unskip"}%   "
   \def\Lrquote{\unskip`}%    `
   \def\Rspace{\spacefactor=1000}%
   \def\Rperiod{\spacefactor=3000}%  .
   \def\Rcomma{\spacefactor=1250}%    ,
   \def\Rquest{\spacefactor=3000}%     ?
   \def\Rcolon{\spacefactor=2000}%   :
   \def\Rscolon{\spacefactor=1250}%   ;
   \def\Rbang{\spacefactor=3000}%      !
   \def\Rquote{\spacefactor=1000}%   '
   \def\Rqquote{\spacefactor=1000}%   "
   \def\Rrquote{\spacefactor=1000}%    `
   }

\def\Refstd{%Standard control strings for formatting bibliography listings.
     \def\Acomma{\unskip, }%between multiple author names
     \def\Aand{\unskip\ and }%between two author names
     \def\Aandd{\unskip\ and }%between last two of multiple author names
     \def\Ecomma{\unskip, }%between multiple editor names
     \def\Eand{\unskip\ and }%between two editor names
     \def\Eandd{\unskip\ and }%between last two of multiple author names
     \def\acomma{\unskip, }%same for authors of reviewed material
     \def\aand{\unskip\ and }%same for authors of reviewed material
     \def\aandd{\unskip\ and }%same for authors of reviewed material
     \def\ecomma{\unskip, }%same for translators
     \def\eand{\unskip\ and }%same for translators
     \def\eandd{\unskip\ and }%same for translators
     \def\Namecomma{\unskip, }%same for citations using authors' names
     \def\Nameand{\unskip\ and }%same for citations using authors' names
     \def\Nameandd{\unskip\ and }%same for citations using authors' names
     \def\Revcomma{\unskip, }%between last and first name of reversed name
     \def\Initper{.\ }%punctuation after initial
     \def\Initgap{\dimen0=\spaceskip\divide\dimen0 by 2\hskip-\dimen0}%
   }

\def\Smallcapsaand{%Smallcaps redefinition of \Aand and \Aandd for \Refstd
     \def\Aand{\unskip\bgroup{\Smallcapsfont\ AND }\egroup}%
     \def\Aandd{\unskip\bgroup{\Smallcapsfont\ AND }\egroup}%
     \def\eand{\unskip\bgroup\Smallcapsfont\ AND \egroup}%
     \def\eandd{\unskip\bgroup\Smallcapsfont\ AND \egroup}%
   }

\def\Smallcapseand{%Smallcaps redefinition of \Eand, \Eeand, etc for Refstd
     \def\Eand{\unskip\bgroup\Smallcapsfont\ AND \egroup}%
     \def\Eandd{\unskip\bgroup\Smallcapsfont\ AND \egroup}%
     \def\aand{\unskip\bgroup\Smallcapsfont\ AND \egroup}%
     \def\aandd{\unskip\bgroup\Smallcapsfont\ AND \egroup}%
   }

   \def\Citefont{}%citations
   \def\ACitefont{}%alternate citations
   \def\Authfont{}%authors
   \def\Titlefont{}%titles
   \def\Tomefont{\sl}%journals or books
   \def\Volfont{}%volume number of journal
   \def\Flagfont{}%citation flag
   \def\Reffont{\rm}%set at beginning of reference listing
   \def\Smallcapsfont{\sevenrm}%small caps
   \def\Flagstyle#1{\hangindent\parindent\indent\hbox to0pt%flag style
       {\hss[{\Flagfont#1}]\kern.5em}\ignorespaces}%        for references

\def\Underlinemark{\vrule height .7pt depth 0pt width 3pc}%for replacing
\def\Citebrackets{\Rpunct%defaults for putting citations in brackets [].
   \def\Lcitemark{\def\Cfont{\Citefont}[\bgroup\Cfont}%mark at left of citation
   \def\Rcitemark{\egroup]}%mark at right of citation
   \def\LAcitemark{\def\Cfont{\ACitefont}\bgroup\ACitefont}%
   \def\RAcitemark{\egroup}%mark at right of alternate citation
   \def\LIcitemark{\egroup}%mark at left of insertion in citation
   \def\RIcitemark{\bgroup\Cfont}%mark at right of insertion in citation
   \def\Citehyphen{\egroup--\bgroup\Cfont}%separater for string of citations
   \def\Citecomma{\egroup,\hskip0pt\bgroup\Cfont}%
   \def\Citebreak{}%mark between parts of citation (e.g. author\Citebreak date)
   }

\def\Citeparen{\Rpunct%defaults for putting citations in parenthesis ().
   \def\Lcitemark{\def\Cfont{\Citefont}(\bgroup\Cfont}%mark at left of citation
   \def\Rcitemark{\egroup)}%mark at right of citation
   \def\LAcitemark{\def\Cfont{\ACitefont}\bgroup\ACitefont}%
   \def\RAcitemark{\egroup}%mark at right of alternate citation
   \def\LIcitemark{\egroup}%mark at left of insertion in citation
   \def\RIcitemark{\bgroup\Cfont}%mark at right of insertion in citation
   \def\Citehyphen{\egroup--\bgroup\Cfont}%separater for string of citations
   \def\Citecomma{\egroup,\hskip0pt\bgroup\Cfont}%
   \def\Citebreak{}%mark between parts of citation (e.g. author\Citebreak date)
   }

\def\Citesuper{\Lpunct%defaults for making superscript citations
   \def\Lcitemark{\def\Cfont{\Citefont}\raise1ex\hbox\bgroup\bgroup\Cfont}%
   \def\Rcitemark{\egroup\egroup}%mark at right of citation
   \def\LAcitemark{\def\Cfont{\ACitefont}\bgroup\ACitefont}%
   \def\RAcitemark{\egroup}%mark at right of alternate citation
   \def\LIcitemark{\egroup\egroup}%mark at left of insertion in citation
   \def\RIcitemark{\raise1ex\hbox\bgroup\bgroup\Cfont}%
   \def\Citehyphen{\egroup--\bgroup\Cfont}%separater for string of citations
   \def\Citecomma{\egroup,\hskip0pt\bgroup%
      \Cfont}%separater for multiple citations
   \def\Citebreak{}%mark between parts of citation (e.g. author\Citebreak date)
   }

\def\Citenamedate{\Rpunct%defaults for making name-date citations
   \def\Lcitemark{%mark at left of citation--also sets internal punctuation
      \def\Citebreak{\egroup\ [\bgroup\Citefont}%separater in citation
      \def\Citecomma{\egroup]; %between multiple citations
         \bgroup\let\uchyph=1\Citefont}(\bgroup\let\uchyph=1\Citefont}%
   \def\Rcitemark{\egroup])}%mark at right of citation
   \def\LAcitemark{%mark at left of alternate citation
      \def\Citebreak{\egroup\ [\bgroup\Citefont}\def\Citecomma{\egroup], %
         \bgroup\ACitefont }\bgroup\let\uchyph=1\ACitefont}%
   \def\RAcitemark{\egroup]}%mark at right of alternate citation
  \def\Citehyphen{\egroup--\bgroup\Citefont}%separater for string of citations
   \def\LIcitemark{\egroup}%mark at left of insertion in citation
   \def\RIcitemark{\bgroup\Citefont}%mark at right of insertion in citation
   }
\def\annotations{}
\def\Flagstyle#1{\par\noindent\hbox to 30pt{[#1]\hfil}%
\hangindent=30pt\hangafter=1}

\def\ifundefined#1{\expandafter\ifx\csname#1\endcsname\relax}
\ifundefined{thesis}\def\bye{\end{document}}\fi

\def\endauthor{:\kern6pt}

\Refstd\Citebrackets%set general formats for reference list and citations
\def\Citefont{\rm}\def\Titlefont{\Reffont}\def\Volfont{\bf}\def\Tomefont{\it}%redefine some fonts

\def\journalarticleformat{\Reffont\let\uchyph=1\parindent=1.25pc\def\Comma{}%

\sfcode`\.=1000\sfcode`\?=1000\sfcode`\!=1000\sfcode`\:=1000\sfcode`\;=1000\sfcode`\,=1000%\frenchspacing
                \par\vfil\penalty-200\vfilneg%\filbreak
      \if\Ftest\present\Flagstyle\Fstr\fi%

\if\Atest\present\bgroup\Authfont\Astr\egroup\def\Comma{\unskip\endauthor}\fi%
        \if\Ttest\present\Comma\bgroup\Titlefont\Tstr\egroup\def\Comma{, }\fi%
         \if\etest\present\hskip.2em(\bgroup\estr\egroup)\def\Comma{\unskip,
}\fi%
          \if\Jtest\present\Comma\bgroup\Tomefont\Jstr\/\egroup\def\Comma{,
}\fi%

\if\Vtest\present\if\Jtest\present\hskip.2em\else\Comma\fi\bgroup\Volfont\Vstr\egroup\def\Comma{, }\fi%
            \if\Dtest\present\hskip.2em(\bgroup\Dstr\egroup)\def\Comma{, }\fi%
             \if\Ptest\present\bgroup, \Pstr\egroup\def\Comma{, }\fi%
              \if\ttest\present\Comma\bgroup\Titlefont\tstr\egroup\def\Comma{,
}\fi%

\if\jtest\present\Comma\bgroup\Tomefont\jstr\/\egroup\def\Comma{, }\fi%

\if\vtest,\present\if\jtest\present\hskip.2em\else\Comma\fi\bgroup\Volfont\vstr\egroup\def\Comma{, }\fi%
                 \if\dtest\present\hskip.2em(\bgroup\dstr\egroup)\def\Comma{,
}\fi%
                  \if\ptest\present\bgroup, \pstr\egroup\def\Comma{, }\fi%
                   \if\Gtest\present{\Comma Gov't ordering no.
}\bgroup\Gstr\egroup\def\Comma{, }\fi%
                    \if\Mtest\present\Comma MR
\#\bgroup\Mstr\egroup\def\Comma{, }\fi%
                     \if\Otest\present{\Comma\bgroup\Ostr\egroup.}\else{.}\fi%
\annotations  %  MODIFICATION
\vskip5ptplus1ptminus1pt}%\smallskip (but 3 is now 5)  MODIFICATION

\def\conferencereportformat{\Reffont\let\uchyph=1\parindent=1.25pc\def\Comma{}%

\sfcode`\.=1000\sfcode`\?=1000\sfcode`\!=1000\sfcode`\:=1000\sfcode`\;=1000\sfcode`\,=1000%\frenchspacing
                \par\vfil\penalty-200\vfilneg%\filbreak
      \if\Ftest\present\Flagstyle\Fstr\fi%

\if\Atest\present\bgroup\Authfont\Astr\egroup\def\Comma{\unskip\endauthor}\fi%
        \if\Ttest\present\Comma\bgroup\Titlefont\Tstr\egroup\def\Comma{, }\fi%
         \if\Jtest\present\Comma\bgroup\Tomefont\Jstr\/\egroup\def\Comma{,
}\fi%
          \if\Ctest\present\Comma\bgroup\Cstr\egroup\def\Comma{, }\fi%
           \if\Dtest\present\hskip.2em(\bgroup\Dstr\egroup)\def\Comma{, }\fi%
            \if\Mtest\present\Comma MR \#\bgroup\Mstr\egroup\def\Comma{, }\fi%
             \if\Otest\present{\Comma\bgroup\Ostr\egroup.}\else{.}\fi%
\annotations  %  MODIFICATION
\vskip5ptplus1ptminus1pt}%\smallskip (but 3 is now 5)  MODIFICATION

\def\bookarticleformat{\Reffont\let\uchyph=1\parindent=1.25pc\def\Comma{}%

\sfcode`\.=1000\sfcode`\?=1000\sfcode`\!=1000\sfcode`\:=1000\sfcode`\;=1000\sfcode`\,=1000%\frenchspacing
                \par\vfil\penalty-200\vfilneg%\filbreak
      \if\Ftest\present\Flagstyle\Fstr\fi%

\if\Atest\present\bgroup\Authfont\Astr\egroup\def\Comma{\unskip\endauthor}\fi%
        \if\Ttest\present\Comma\bgroup\Titlefont\Tstr\egroup\def\Comma{, }\fi%
         \if\etest\present\hskip.2em(\bgroup\estr\egroup)\def\Comma{\unskip,
}\fi%
          \if\Btest\present\Comma in
\bgroup\Tomefont\Bstr\/\egroup\def\Comma{\unskip, }\fi%
           \if\otest\present\Comma\ \bgroup\ostr\egroup\def\Comma{, }\fi%
            \if\Etest\present\Comma\bgroup\Estr\egroup\unskip,
\ifnum\Ecnt>1eds.\else ed.\fi\def\Comma{, }\fi%
             \if\Stest\present\Comma\bgroup\Sstr\egroup\def\Comma{, }\fi%
              \if\Vtest\present\ \bgroup\bf\Vstr\egroup\def\Comma{, }\fi%
               \if\Ntest\present\Comma no. \bgroup\Nstr\egroup\def\Comma{,
}\fi%
                \if\Itest\present\Comma\bgroup\Istr\egroup\def\Comma{, }\fi%
                 \if\Ctest\present\ (\bgroup\Cstr\egroup)\def\Comma{, }\fi%
                  \if\Dtest\present\Comma\bgroup\Dstr\egroup\def\Comma{, }\fi%
                   \if\Ptest\present\Comma pp.\ \Pstr\def\Comma{, }\fi%
\if\ttest\present\Comma\bgroup\Titlefont\Tstr\egroup\def\Comma{, }\fi%
                     \if\btest\present\Comma in
\bgroup\Tomefont\bstr\egroup\def\Comma{, }\fi%
                       \if\atest\present\Comma\bgroup\astr\egroup\unskip,
\if\acnt\present eds.\else ed.\fi\def\Comma{, }\fi%
                        \if\stest\present\Comma\bgroup\sstr\egroup\def\Comma{,
}\fi%
                         \if\vtest\present\Comma vol.
\bgroup\vstr\egroup\def\Comma{, }\fi%
                          \if\ntest\present\Comma no.
\bgroup\nstr\egroup\def\Comma{, }\fi%

\if\itest\present\Comma\bgroup\istr\egroup\def\Comma{, }\fi%

\if\ctest\present\Comma\bgroup\cstr\egroup\def\Comma{, }\fi%

\if\dtest\present\Comma\bgroup\dstr\egroup\def\Comma{, }\fi%
                              \if\ptest\present\Comma\pstr\def\Comma{, }\fi%
                               \if\Gtest\present{\Comma Gov't ordering no.
}\bgroup\Gstr\egroup\def\Comma{, }\fi%
                                \if\Mtest\present\Comma MR
\#\bgroup\Mstr\egroup\def\Comma{, }\fi%

\if\Otest\present{\Comma\bgroup\Ostr\egroup.}\else{.}\fi%
\annotations  %  MODIFICATION
\vskip5ptplus1ptminus1pt}%\smallskip (but 3 is now 5)  MODIFICATION

\def\bookformat{\Reffont\let\uchyph=1\parindent=1.25pc\def\Comma{}%

\sfcode`\.=1000\sfcode`\?=1000\sfcode`\!=1000\sfcode`\:=1000\sfcode`\;=1000\sfcode`\,=1000%\frenchspacing
                \par\vfil\penalty-200\vfilneg%\filbreak
      \if\Ftest\present\Flagstyle\Fstr\fi%

\if\Atest\present\bgroup\Authfont\Astr\egroup\def\Comma{\unskip\endauthor}%

\else\if\Etest\present\bgroup\def\Eand{\Aand}\def\Eandd{\Aandd}\Authfont\Estr\egroup\unskip, \ifnum\Ecnt>1eds.\else ed.\fi\def\Comma{, }%

\else\if\Itest\present\bgroup\Authfont\Istr\egroup\def\Comma{, }\fi\fi\fi%

\if\Ttest\present\Comma\bgroup\Tomefont\Tstr\/\egroup\def\Comma{\unskip, }%

\else\if\Btest\present\Comma\bgroup\Titlefont\Bstr\/\egroup\def\Comma{\unskip,
}\fi\fi%
            \if\otest\present\Comma\ \bgroup\ostr\egroup\def\Comma{, }\fi%

\if\etest\present\hskip.2em(\bgroup\estr\egroup)\def\Comma{\unskip, }\fi%
              \if\Stest\present\Comma\bgroup\Sstr\egroup\def\Comma{, }\fi%
               \if\Vtest\present\ \bgroup\bf\Vstr\egroup\def\Comma{, }\fi%
                \if\Ntest\present\Comma no. \bgroup\Nstr\egroup\def\Comma{,
}\fi%
                 \if\Atest\present\if\Itest\present
                         \Comma\bgroup\Istr\egroup\def\Comma{\unskip, }\fi%
                      \else\if\Etest\present\if\Itest\present
                              \Comma\bgroup\Istr\egroup\def\Comma{\unskip,
}\fi\fi\fi%
                     \if\Ctest\present\ (\bgroup\Cstr\egroup)\def\Comma{, }\fi%
                      \if\Dtest\present\Comma\bgroup\Dstr\egroup\def\Comma{,
}\fi%
                      \if\Ptest\present\bgroup, pp.\ \Pstr\egroup\def\Comma{,
}\fi% MODIFICATION

\if\ttest\present\Comma\bgroup\Titlefont\tstr\egroup\def\Comma{, }%

\else\if\btest\present\Comma\bgroup\Titlefont\bstr\egroup\def\Comma{, }\fi\fi%

\if\stest\present\Comma\bgroup\sstr\egroup\def\Comma{, }\fi%
                           \if\vtest\present\Comma vol.
\bgroup\vstr\egroup\def\Comma{, }\fi%
                            \if\ntest\present\Comma no.
\bgroup\nstr\egroup\def\Comma{, }\fi%

\if\itest\present\Comma\bgroup\istr\egroup\def\Comma{, }\fi%

\if\ctest\present\Comma\bgroup\cstr\egroup\def\Comma{, }\fi%

\if\dtest\present\Comma\bgroup\dstr\egroup\def\Comma{, }\fi%
                                \if\Gtest\present{\Comma Gov't ordering no.
}\bgroup\Gstr\egroup\def\Comma{, }\fi%
                                 \if\Mtest\present\Comma MR
\#\bgroup\Mstr\egroup\def\Comma{, }\fi%

\if\Otest\present{\Comma\bgroup\Ostr\egroup.}\else{.}\fi%
\annotations  %  MODIFICATION
\vskip5ptplus1ptminus1pt}%\smallskip (but 3 is now 5)  MODIFICATION

\def\technicalreportformat{\Reffont\let\uchyph=1\parindent=1.25pc\def\Comma{}%

\sfcode`\.=1000\sfcode`\?=1000\sfcode`\!=1000\sfcode`\:=1000\sfcode`\;=1000\sfcode`\,=1000%\frenchspacing
                \par\vfil\penalty-200\vfilneg%\filbreak
      \if\Ftest\present\Flagstyle\Fstr\fi%

\if\Atest\present\bgroup\Authfont\Astr\egroup\def\Comma{\unskip\endauthor}%

\else\if\Etest\present\bgroup\def\Eand{\Aand}\def\Eandd{\Aandd}\Authfont\Estr\egroup\unskip, \ifnum\Ecnt>1eds.\else ed.\fi\def\Comma{, }%

\else\if\Itest\present\bgroup\Authfont\Istr\egroup\def\Comma{, }\fi\fi\fi%
          \if\Ttest\present\Comma\bgroup\Titlefont\Tstr\egroup\def\Comma{,
}\fi%
           \if\Atest\present\if\Itest\present
                   \Comma\bgroup\Istr\egroup\def\Comma{, }\fi%
                \else\if\Etest\present\if\Itest\present
                        \Comma\bgroup\Istr\egroup\def\Comma{, }\fi\fi\fi%
            \if\Rtest\present\Comma\bgroup\Rstr\egroup\def\Comma{, }\fi%
             \if\Ctest\present\ (\bgroup\Cstr\egroup)\def\Comma{, }\fi%
              \if\Dtest\present\Comma\bgroup\Dstr\egroup\def\Comma{, }\fi%
              \if\Ptest\present\bgroup, \Pstr\egroup\def\Comma{, }\fi%
               \if\ttest\present\Comma\bgroup\Titlefont\tstr\egroup\def\Comma{,
}\fi%
                \if\itest\present\Comma\bgroup\istr\egroup\def\Comma{, }\fi%
                 \if\rtest\present\Comma\bgroup\rstr\egroup\def\Comma{, }\fi%
                  \if\ctest\present\Comma\bgroup\cstr\egroup\def\Comma{, }\fi%
                   \if\dtest\present\Comma\bgroup\dstr\egroup\def\Comma{, }\fi%
                    \if\Gtest\present{\Comma Gov't ordering no.
}\bgroup\Gstr\egroup\def\Comma{, }\fi%
                     \if\Mtest\present\Comma MR
\#\bgroup\Mstr\egroup\def\Comma{, }\fi%
                      \if\Otest\present{\Comma\bgroup\Ostr\egroup.}\else{.}\fi%
\annotations  %  MODIFICATION
\vskip5ptplus1ptminus1pt}%\smallskip (but 3 is now 5)  MODIFICATION

\def\otherformat{\Reffont\let\uchyph=1\parindent=1.25pc\def\Comma{}%

\sfcode`\.=1000\sfcode`\?=1000\sfcode`\!=1000\sfcode`\:=1000\sfcode`\;=1000\sfcode`\,=1000%\frenchspacing
                \par\vfil\penalty-200\vfilneg%\filbreak
      \if\Ftest\present\Flagstyle\Fstr\fi%

\if\Atest\present\bgroup\Authfont\Astr\egroup\def\Comma{\unskip\endauthor}%

\else\if\Etest\present\bgroup\def\Eand{\Aand}\def\Eandd{\Aandd}\Authfont\Estr\egroup\unskip, \ifnum\Ecnt>1eds.\else ed.\fi\def\Comma{, }%

\else\if\Itest\present\bgroup\Authfont\Istr\egroup\def\Comma{, }\fi\fi\fi%
          \if\Ttest\present\Comma\bgroup\Titlefont\Tstr\egroup\def\Comma{,
}\fi%
            \if\Atest\present\if\Itest\present
                    \Comma\bgroup\Istr\egroup\def\Comma{, }\fi%
                 \else\if\Etest\present\if\Itest\present
                         \Comma\bgroup\Istr\egroup\def\Comma{, }\fi\fi\fi%
                 \if\Ctest\present\ (\bgroup\Cstr\egroup)\def\Comma{, }\fi%
                  \if\Dtest\present\Comma\bgroup\Dstr\egroup\def\Comma{, }\fi%
                  \if\Ptest\present\bgroup, \Pstr\egroup\def\Comma{, }\fi%
                   \if\Gtest\present{\Comma Gov't ordering no.
}\bgroup\Gstr\egroup\def\Comma{, }\fi%
                    \if\Mtest\present\Comma MR
\#\bgroup\Mstr\egroup\def\Comma{, }\fi%
                     \if\Otest\present{\Comma\bgroup\Ostr\egroup.}\else{.}\fi%
\annotations  %  MODIFICATION
\vskip5ptplus1ptminus1pt}%\smallskip (but 3 is now 5)  MODIFICATION
\catcode`\@=11
\def\matrixx#1{\null\,\vcenter{\normalbaselines\m@th
    \ialign{\hfil$\displaystyle ##$\hfil&&\quad\hfil$\displaystyle
##$\hfil\crcr
    \mathstrut\crcr\noalign{\kern-\baselineskip}
    #1\crcr\mathstrut\crcr\noalign{\kern-\baselineskip}}}\,}
\catcode`\@=12

\def\annotations{}
\def\Flagstyle#1{\par\noindent\hbox to 30pt{[#1]\hfil}%
\hangindent=30pt\hangafter=1}

\catcode`\@=11
\def\ATstartsection{\@startsection}
\def\resetATfont{\reset@font}
\def\zAT{\z@}
\def\RIfM@{\relax\protect\ifmmode}
\def\fraktur{\mathfrak}

\catcode`\@=12

\def\hexnumberZ#1{\ifnum#1<10 \number#1\else
 \ifnum#1=10 A\else\ifnum#1=11 B\else\ifnum#1=12 C\else
 \ifnum#1=13 D\else\ifnum#1=14 E\else\ifnum#1=15 F\fi\fi\fi\fi\fi\fi\fi}

\newcounter{alphactr}

\newcounter{romanctr}

\newcounter{arabicctr}

\def\geom{g\'eom\'etrie}
\def\EXPOSE{expos\'e}
\def\SGA{S\'em\-in\-aire de G\'eom\'e\-trie Al\-g\'e\-bri\-que}
\def\PLUS{+}
\def\TIMES{*}

\def\begindiagram{\def\normalbaselines{\baselineskip20pt
      \lineskip3pt \lineskiplimit3pt}}

\catcode`\@=11
\def\operatoratfont{\operator@font}
\let\@afterindentfalse\@afterindenttrue
\@afterindenttrue
\catcode`\@=12

\def\op{\operatoratfont}

\catcode`\@=11
\def\newtheoremrm#1{\@ifnextchar[{\@othmrm{#1}}{\@nthmrm{#1}}}
\def\@nthmrm#1#2{\@ifnextchar[{\@xnthmrm{#1}{#2}}{\@ynthmrm{#1}{#2}}}
\def\@xnthmrm#1#2[#3]{\expandafter\@ifdefinable\csname #1\endcsname%
{\@definecounter{#1}\@addtoreset{#1}{#3%
}\expandafter\xdef\csname the#1\endcsname{\expandafter\noexpand%
\csname the#3\endcsname . \noexpand\arabic{#1}}\global\@namedef{#1%
}{\@thmrm{#1}{#2}}\global\@namedef{end#1}{\endtrivlist}}}
\def\@ynthmrm#1#2{\expandafter\@ifdefinable\csname #1\endcsname%
{\@definecounter {#1}\expandafter\xdef\csname the#1\endcsname%
{\noexpand\arabic{#1}}\global\@namedef{#1}{\@thmrm{#1}{#2}%
}\global\@namedef{end#1}{\endtrivlist}}}
\def\@othmrm#1[#2]#3{\expandafter\@ifdefinable\csname #1\endcsname
{\global\@namedef{the#1}{\@nameuse{the#2}}\global\@namedef{#1%
}{\@thmrm{#2}{#3}}\global\@namedef{end#1}{\endtrivlist}}}
\def\@thmrm#1#2{\refstepcounter{#1}\@ifnextchar[{\@ythmrm{#1}{#2}%
}{\@xthmrm{#1}{#2}}}
\def\@xthmrm#1#2{\trivlist \item[\hskip \labelsep{\bf #2\ \csname
the#1\endcsname}]\ignorespaces}
\def\@ythmrm#1#2[#3]{\trivlist \item[\hskip \labelsep{\bf #2\ \csname
the#1\endcsname\ (#3)}]\ignorespaces}
\catcode`\@=12

\newtheoremrm{remark}[theorem]{Remark}

\def\manyxx#1{\if#1\PLUS{+}\else\if#1\TIMES{*}\fi\fi}

\def\xhmode#1{\relax\ifmmode{#1}\else{\snap\hbox{$#1$}}\fi}
\def\snap{\hskip 0pt plus 4cm\penalty1000\hskip 0pt plus -4cm}
\def\PP{\xmode{\Bbb P\kern1pt}}
\def\email{jaffe{\kern0.5pt}@{\kern0.5pt}cpthree.unl.edu}
\def\C{{\Bbb C}\kern1pt}
\def\Q{{\Bbb Q}\kern1pt}
\def\F{{\Bbb F}\kern1pt}
\def\P#1{{\PP^{#1}}}
\def\N{\xmode{\Bbb N}}
\def\Z{\xmode{\Bbb Z}}
\def\O{{\cal O}}
\def\Ann{\mathop{\operatoratfont Ann}\nolimits}
\def\Aut{\mathop{\operatoratfont Aut}\nolimits}
\def\Coker{\mathop{\operatoratfont Coker}\nolimits}
\def\Eq{\mathop{\operatoratfont Eq}\nolimits}
\def\End{\mathop{\operatoratfont End}\nolimits}
\def\GL{\mathop{\operatoratfont GL}\nolimits}
\def\Hom{\mathop{\operatoratfont Hom}\nolimits}
\def\Im{\mathop{\operatoratfont Im}\nolimits}
\def\Ker{\mathop{\operatoratfont Ker}\nolimits}
\def\Mat{\mathop{\operatoratfont Mat}\nolimits}
\def\Pic{\mathop{\operatoratfont Pic}\nolimits}
\def\Spec{\mathop{\operatoratfont Spec}\nolimits}
\def\Tor{\mathop{\operatoratfont Tor}\nolimits}
\def\coker{\mathop{\operatoratfont coker}\nolimits}
\def\fin{\mathop{\operatoratfont fin}\nolimits}
\def\id{\mathop{\operatoratfont id}\nolimits}
\def\AUT{\mathop{\mathbf{Aut}}\nolimits}
\def\HOM{\mathop{\mathbf{Hom}}\nolimits}
\def\END{\mathop{\mathbf{End}}\nolimits}
\def\BIL{\mathop{\mathbf{Bil}}\nolimits}
\def\HOMalg{\mathop{{\mathbf{Hom}}_{\kern1pt\operatoratfont alg}}
     \nolimits}
\def\AUTalg{\mathop{{\mathbf{Aut}}_{\kern1pt\operatoratfont alg}}
     \nolimits}
\def\uHom{\mathop{\ul{\operatoratfont Hom}}\nolimits}
\def\uKer{\mathop{\ul{\operatoratfont Ker}}\nolimits}
\def\uIm{\mathop{\ul{\operatoratfont Im}}\nolimits}
\def\uEnd{\mathop{\ul{\operatoratfont End}}\nolimits}
\def\uPGL{\mathop{\ul{\operatoratfont PGL}}\nolimits}
\def\uAut{\mathop{\ul{\operatoratfont Aut}}\nolimits}
\def\uTor{\mathop{\ul{\operatoratfont Tor}}\nolimits}
\def\uAnn{\mathop{\ul{\operatoratfont Ann}}\nolimits}
\def\mp[[#1||#2||#3]]{\snap\hbox{$#1:\kern3pt #2\ \to\ #3$}}
\def\mapx[[#1||#2]]{\snap\hbox{$#1\ \to\ #2$}}
\def\dmap[[#1||#2||#3]]%
{     \setbox0=\hbox{$\displaystyle #1:\kern3pt #2\ \ \mapE\ \ \ #3$%
}     \ifdim\wd0<\hsize$$#1:\kern3pt #2\ \ \mapE\ \ \ #3$$\else%
     \begin{flushleft} $\displaystyle #1:\kern3pt #2$
     \end{flushleft} \begin{flushright} $\displaystyle\mapE{}\ \ \ #3$
     \end{flushright}\fi}
\def\dmapx[[#1||#2]]%
{     \setbox0=\hbox{$\displaystyle #1\ \ \mapE\ \ \ #2$%
}     \ifdim\wd0<\hsize$$#1\ \ \mapE\ \ \ #2$$\else%
     \begin{flushleft} $\displaystyle #1$ \end{flushleft}
     \begin{flushright} $\displaystyle\mapE{}\ \ \ #2$ \end{flushright}\fi}
\def\inverselim{\displaystyle{\lim_{\overleftarrow{\hphantom{\lim}}}}}
\def\invlim#1{\displaystyle{\lim_{\scriptstyle
    {\overleftarrow{\hphantom{\lim}}}\atop\scriptstyle{#1}}}}
\def\mapE#1{{\smash{\mathop{\longrightarrow}\limits^{#1}}}}
\def\smapE#1{{\smash{\mathop{\kern-10pt\rightarrow\kern-10pt}\limits^{#1}}}}
\def\includeE#1{{\lhook\kern-3.5pt\joinrel\smash{
    \mathop{\longrightarrow}\limits^{#1}}}}
\def\mapW#1{{\smash{\mathop{\longleftarrow}\limits^{#1}}}}
\def\mapS#1{\Big\downarrow\rlap{$\vcenter{\hbox{$\scriptstyle{#1}$}}$}}
\def\mapN#1{\Big\uparrow\rlap{$\vcenter{\hbox{$\scriptstyle{#1}$}}$}}
\def\diagramx#1{{$$\begindiagram\matrixx{#1}$$}}
\def\splitdiagram#1#2{
    \begin{flushleft}$\displaystyle\begindiagram\matrixx{#1}$\end{flushleft}
    \begin{flushright}$\displaystyle\begindiagram\matrixx{#2}$\end{flushright}}
\def\lesmaps#1#2#3#4#5{\diagramx{0&\mapE{}&#1&\mapE{#2}&#3&\mapE{#4}&#5}}
\def\lesx#1#2#3{\diagramx{#1&\rightarrowtail&#2&\mapE{}&#3}}
\def\lesxdot#1#2#3{\diagramx{#1&\rightarrowtail&#2&\mapE{}&#3.}}
\def\res#1#2#3{\diagramx{#1&\mapE{}&#2&\mapE{}&#3&\mapE{}&0}}
\def\resdot#1#2#3{\diagramx{#1&\mapE{}&#2&\mapE{}&#3&\mapE{}&0.}}
\def\resmapsdot#1#2#3#4#5{\diagramx{#1&\mapE{#2}&#3&\mapE{#4}&#5&\mapE{}&0.}}
\def\ses#1#2#3%
{  \setbox0=\hbox{$\matrixx{0&\mapE{}&#1&\mapE{}&#2&\mapE{}&#3&\mapE{}&0}$%
}  \ifdim\wd0<\hsize\diagramx{0&\mapE{}&#1&\mapE{}&#2&\mapE{}&#3&\mapE{}&
  0}\else\setbox0=\hbox{$\matrixx{0&\smapE{}&#1&\smapE{}&#2&\smapE{}&#3&
  \smapE{}&0}$}\ifdim\wd0<\hsize\diagramx{0&\smapE{}&#1&\smapE{}&#2&\smapE{}&
  #3&\smapE{}&0}\else\splitdiagram{0&\mapE{}&#1&\mapE{}&#2}{\mapE{}&#3
  &\mapE{}&0}\fi\fi}
\def\sesdot#1#2#3%
{  \setbox0=\hbox{$\matrixx{0&\mapE{}&#1&\mapE{}&#2&\mapE{}&#3&\mapE{}&0.}$%
}  \ifdim\wd0<\hsize\diagramx{0&\mapE{}&#1&\mapE{}&#2&\mapE{}&#3&\mapE{}&
  0.}\else\setbox0=\hbox{$\matrixx{0&\smapE{}&#1&\smapE{}&#2&\smapE{}&#3&
  \smapE{}&0.}$}\ifdim\wd0<\hsize\diagramx{0&\smapE{}&#1&\smapE{}&#2&\smapE{}&
  #3&\smapE{}&0.}\else\splitdiagram{0&\mapE{}&#1&\mapE{}&#2}{\mapE{}&#3
  &\mapE{}&0.}\fi\fi}
\def\les#1#2#3%
{  \setbox0=\hbox{$\matrixx{0&\mapE{}&#1&\mapE{}&#2&\mapE{}&#3}$%
}  \ifdim\wd0<\hsize\diagramx{0&\mapE{}&#1&\mapE{}&#2&\mapE{}&
  #3}\else\setbox0=\hbox{$\matrixx{0&\smapE{}&#1&\smapE{}&#2&\smapE{}&
  #3}$}\ifdim\wd0<\hsize\diagramx{0&\smapE{}&#1&\smapE{}&#2&\smapE{}&
  #3}\else\splitdiagram{0&\mapE{}&#1&\mapE{}&#2}{\mapE{}&#3}\fi\fi}
\def\lesdot#1#2#3%
{  \setbox0=\hbox{$\matrixx{0&\mapE{}&#1&\mapE{}&#2&\mapE{}&#3.}$%
}  \ifdim\wd0<\hsize\diagramx{0&\mapE{}&#1&\mapE{}&#2&\mapE{}&
  #3.}\else\setbox0=\hbox{$\matrixx{0&\smapE{}&#1&\smapE{}&#2&\smapE{}&
  #3.}$}\ifdim\wd0<\hsize\diagramx{0&\smapE{}&#1&\smapE{}&#2&\smapE{}&
  #3.}\else\splitdiagram{0&\mapE{}&#1&\mapE{}&#2}{\mapE{}&#3.}\fi\fi}
\def\squareSE#1#2#3#4{\diagramx{#1&\mapE{}&#2\cr
		      \mapS{}&&\mapS{}\cr #3&\mapE{}&#4\cr}}
\def\makenull#1{{\setbox0=\hbox{#1}\wd0=0pt\box0}}
\def\nulldot{\makenull{.}}
\def\nullcomma{\makenull{,}}
\def\rowfour#1#2#3#4{#1&\mapE{}&#2&\mapE{}&#3&\mapE{}&#4}
\def\rowfive#1#2#3#4#5{#1&\mapE{}&#2&\mapE{}&#3&\mapE{}&#4&\mapE{}&#5}
\def\Rowfive#1#2#3#4#5{\diagramx{#1&\mapE{}&#2&\mapE{}&#3&\mapE{}&#4&\mapE{}&#5\cr}}
\def\o#1{\if#1\PLUS\oplus\else\if#1\TIMES\otimes\fi\fi}
\def\manyo#1{\o#1 \cdots \o#1}
\def\manycap{\cap \cdots \cap}
\def\manyIN{\IN \cdots \IN}
\def\manytimes{\times \cdots \times}
\def\manycdot{\cdot \ldots \cdot}
\def\tA{{\tilde{A}}}
\def\tF{{\tilde{F}}}
\def\tG{{\tilde{G}}}
\def\ltF{{\tilde{\lowercase{F}}}}
\def\ltK{{\tilde{\lowercase{K}}}}
\def\ltT{{\tilde{\lowercase{T}}}}
\def\ltX{{\tilde{\lowercase{X}}}}
\def\oC{{\overline{C}}}
\def\oS{{\overline{S}}}
\def\uA{{\underline{A}}}
\def\uD{{\underline{D}}}
\def\uJ{{\underline{J}}}
\def\uK{{\underline{K}}}
\def\uM{{\underline{M}}}
\def\uN{{\underline{N}}}
\def\uP{{\underline{P}}}
\def\shA{{\cal{A}}}
\def\shC{{\cal{C}}}
\def\shD{{\cal{D}}}
\def\shF{{\cal{F}}}
\def\shJ{{\cal{J}}}
\def\shM{{\cal{M}}}
\def\shN{{\cal{N}}}
\def\shP{{\cal{P}}}
\def\shS{{\cal{S}}}
\def\ushM{{\underline{\cal M}}}
\def\ushN{{\underline{\cal N}}}
\def\ushP{{\underline{\cal P}}}
\def\lfC{{\xmode{{\fraktur{\lowercase{C}}}}}}
\def\lfM{{\xmode{{\fraktur{\lowercase{M}}}}}}
\def\lfP{{\xmode{{\fraktur{\lowercase{P}}}}}}
\def\lfQ{{\xmode{{\fraktur{\lowercase{Q}}}}}}
\def\tlfM{{\xmode{{\tilde{\fraktur{\lowercase{M}}}}}}}
\def\DVR{discrete valuation ring}
\def\IFF{if and only if}
\def\WRT{with respect to}
\def\FG{finitely generated}
\def\QC{quasi-coherent}
\def\MQC{module-quasi-coherent}
\def\MC{module-coherent}
\def\MV{module-valued}
\def\WMAT{we may assume that}
\def\th#1{\xmode{{#1}^{\operatoratfont th}}}
\def\IN{\subset}
\def\setof#1{\xmode{\{#1\}}}
\def\setofh#1{\xmode{\{\hbox{#1}\}}}
\def\iso{\cong}
\def\qed{{\hfill$\square$}}
\def\formulaqed#1{$$#1\eqno\square$$}
\def\nor#1{#1_{\operatoratfont nor}}
\def\RED#1{#1_{\operatoratfont red}}
\def\vec#1#2#3{#1_{#2}, \ldots, #1_{#3}}
\def\sets#1#2#3{\xmode{\setof{{#1}_{#2}}_{#2 \in #3}}}
\def\xmode#1{\relax\ifmmode{#1}\else{$#1$}\fi}
\def\cat#1{\xhmode{\ll\negthinspace\hbox{#1}\negthinspace\gg}}
\def\opcat#1{\xmode{\ll\negthinspace\hbox{#1}\negthinspace\gg^\circ}}
\def\fun[[ #1 || #2 || #3 ]]{\mp[[ {#1} || \cat{#2} || \cat{#3} ]]}
\def\funx[[ #1 || #2 ]]{\mapx[[ \cat{#1} || \cat{#2} ]]}
\def\pref#1{(\ref{#1})}
\def\block#1{\section{#1}}
\def\abs#1{{\vert{#1}\vert}}
\def\makeaddress{
     \vskip 0.15in
     \par\noindent {\footnotesize Department of Mathematics and Statistics,
                                  University of Nebraska}
     \par\noindent {\footnotesize Lincoln, NE 68588-0323, USA\ \ (\email)}}
\def\meatnebraska{ \par\noindent David B.\ Jaffe \makeaddress}
\newenvironment{proof}{\trivlist \item[\hskip \labelsep{\sc
Proof.\kern1pt}]}{\endtrivlist}%\proof
 \newenvironment{proofnodot}{\trivlist \item[\hskip \labelsep{\sc
Proof}]}{\endtrivlist}%\proofnodot
 \newenvironment{sketch}{\trivlist \item[\hskip \labelsep{\sc
Sketch.\kern1pt}]}{\endtrivlist}%\sketch

\newenvironment{alphalist}{\begin{list}{(\alph{alphactr})}{\usecounter{alphactr}}}{\end{list}}%\alphalist

\newenvironment{arabiclist}{\begin{list}{(\arabic{arabicctr})}{\usecounter{arabicctr}}}{\end{list}}%\arabiclist
 \def\bxc{bounded\/\kern1.5pt$\times$\discretionary{}{}{\kern1.5pt}cobounded}
\def\liL{{\lambda \in \Lambda}}
\def\ul{\underline}
\def\ol{\overline}
\newenvironment{definition}{\trivlist \item[\hskip
       \labelsep{\bf Definition.\kern1pt}]}{\endtrivlist}
\newenvironment{remarks}{\trivlist \item[\hskip
       \labelsep{\bf Remarks.\kern1pt}]}{\endtrivlist}
\newenvironment{problem}{\trivlist \item[\hskip
       \labelsep{\bf Problem.\kern1pt}]}{\endtrivlist}
\documentclass{article}\usepackage{amssymb}
\input xypic.sty
\newtheorem{theorem}{Theorem}[section]
 \newtheorem{lemma}[theorem]{Lemma}%\lemma
 \newtheorem{conjecture}[theorem]{Conjecture}%\conjecture
 \newtheorem{corollary}[theorem]{Corollary}%\corollary
 \newtheorem{prop}[theorem]{Proposition}%\prop
 \newtheorem{problemx}[theorem]{Problem}%numbered problem
 \newtheorem{exampleth}[theorem]{Example}
\newenvironment{example}{\begin{exampleth}\fontshape{n}\selectfont}{\end{exampleth}}

\def\ifundefined#1{\expandafter\ifx\csname#1\endcsname\relax}

\begin{document}
\vskip 0.15in
\def\nullfoot#1{}
\def\thefootnote{\nullfoot{footnote}}
\footnote{{\it 1991 Mathematics Subject Classification.}
Primary: 14C22, 18A25, 14K30, 18A40.}
\footnote{{\it Key words and phrases.\/} Coherent functor, representable
functor, Picard group.}
\footnote{Partially supported by the National Science Foundation.}
\def\thefootnote{\fnsymbol{footnote}}\setcounter{footnote}{0}
\par\noindent{\bf\LARGE Coherent functors, with application
to torsion in the Picard group}
\def\thefootnote{\arabic{footnote}}\setcounter{footnote}{0}
\vspace*{15pt}
\meatnebraska
\vspace*{0.2in}

\medskip
\par\noindent{\bf\large Abstract}
\medskip

Let $A$ be a commutative noetherian ring.  We investigate a class of functors
from \cat{commutative $A$-algebras} to \cat{sets}, which we call
{\it coherent}.  When such a functor $F$ in fact takes its values in
\cat{abelian groups}, we show that there are only finitely many prime numbers
$p$ such that ${}_p F(A)$ is infinite, and that none of these primes are
invertible in $A$.  This (and related statements) yield information about
torsion in $\Pic(A)$.  For example, if $A$ is of finite type over $\Z$, we
prove that the torsion in $\Pic(A)$ is supported at a finite set of primes,
and if ${}_p \Pic(A)$ is infinite, then the prime $p$ is not invertible in $A$.
These results use the (already known) fact that if such an $A$ is normal, then
$\Pic(A)$ is \FG.  We obtain a parallel result for a reduced scheme $X$ of
finite type over $\Z$.  We classify the groups which can occur as the Picard
group of a scheme of finite type over a finite field.

\medskip
\par\noindent{\bf\large Coherent functors (introductory remarks)}
\medskip

Let us say that an {\it $A$-functor\/} is a functor from the category of
commutative $A$-algebras to \cat{sets}.  Some such $A$-functors have
additional structure: they are actually functors from
\cat{commutative $A$-algebras} to \cat{groups}.  We refer to such functors as
{\it group-valued\/} $A$-functors.  We will also consider $A$-functors $F$
such that $F(B)$ is a $B$-module for every $B$; these {\it \MV\/}
$A$-functors are discussed later in the introduction.  For now, all
$A$-functors which we consider will be treated as set-valued functors.

An $A$-functor is {\it coherent\/} if it may be built up as an iterated finite
limit of functors of the form $\uM$, given by $\uM(B) = M \o*_A B$, where $M$
is a \FG\ $A$-module.  We do not know if every coherent functor may be
expressed as a finite limit of such functors $\uM$.  However, the analogous
question regarding \MV\ functors is answered affirmatively below.

The idea of {\it coherent functor\/} was originally devised by Auslander
\Lcitemark 5\Rcitemark \Rspace{}, in a somewhat different setting; his
notion of coherence applied to functors from an abelian category to
\cat{abelian groups}.  Later Artin\Lspace \Lcitemark 2\Rcitemark \Rspace{}
transposed Auslander's notion to a setting closer to that given here.  Artin
also raised a question about coherence of higher direct images as functors.
This question is considered in \S\ref{higher-section}.

If an $A$-functor is representable by a commutative $A$-algebra of finite
type, then it is coherent.  There are many examples of non-representable
$A$-functors which are coherent.  For example, if $M$ is a \FG\ $A$-module,
then $B \mapsto \Aut_{B-{\op mod}}(M \o*_A B)$ defines a coherent $A$-functor.
More examples may be found in \S\ref{examples-section}.

A {\it \MV\ $A$-functor\/} is an (abelian group)-valued $A$-functor $F$,
together with the following additional structure: for each commutative
$A$-algebra $B$, $F(B)$ has the structure of a $B$-module, such that for any
homomorphism \mapx[[ B_1 || B_2 ]] of commutative $A$-algebras, the induced
map \mapx[[ F(B_1) || F(B_2) ]] is a homomorphism of $B_1$-modules.
The \MV\ $A$-functors form an abelian category.

A \MV\ $A$-functor $F$ is {\it \MC\/} if there exists a homomorphism
\mp[[ f || M || N ]] of \FG\ $A$-modules such that $F$ is isomorphic to the
\MV\ $A$-functor given by $B \mapsto \ker(f \o*_A B)$.  The \MC\ $A$-functors
form a full subcategory of the category of \MV\ $A$-functors.

Most examples of \MV\ $A$-functors are induced naturally by functors from
\cat{$A$-modules} to \cat{$A$-modules}.  Certainly for many purposes it makes
more sense to study the latter sort of functor.  On the other hand, (as pointed
out by Artin\Lspace \Lcitemark 2\Rcitemark \Rspace{}) one can set up a
correspondence
between \MC\ $A$-functors and functors from \cat{$A$-modules} to
\cat{$A$-modules} which satisfy an analogous coherence axiom.  This creates a
bridge to the ideas of Auslander\Lspace \Lcitemark 5\Rcitemark \Rspace{} and
Grothendieck (\Lcitemark 18\Rcitemark \ \S7).  We have not exploited this point
of view.

A key result is that if \mp[[ \sigma || F || G ]] is a
morphism of \MC\ $A$-functors, then $\ker(\sigma)$ and $\coker(\sigma)$ are
\MC.  (These are to be computed in \cat{\MV\ $A$-functors}.)  One deduces
easily from this that any finite limit or finite colimit of \MC\ $A$-functors
is \MC.  In particular, the iterated finite limit construction which we used
in the definition of {\it coherent\/} $A$-functor is not necessary here,
although in the body of the paper we find it convenient to begin with a
definition of \MC\ which uses iterated finite limits.

Let $X$ be a noetherian scheme.  An {\it $X$-functor\/} is a functor from
\snap\hbox{\opcat{$X$-schemes}} to \cat{sets}.  We define the notion of
coherent $X$-functor by analogy with the definition for $A$-functors.  When
$X = \Spec(A)$, the theory of coherent $X$-functors is identical to the theory
of coherent $A$-functors.  Similarly, we define \MC\ $X$-functors.  We show
that the property of being a \MC\ $X$-functor is local on $X$, assuming that
$X$ is separated.  We conjecture that the property of being a coherent
$X$-functor is local on $X$.

%\section*{Finiteness theorems (introductory remarks)}
\medskip
\par\noindent{\bf\large Finiteness theorems (introductory remarks)}
\medskip

If $F$ is an (abelian group)-valued coherent $X$-functor, we prove that there
are only finitely many primes $p$ such that ${}_p F(X)$ is infinite, and that
none of these primes are invertible in $\Gamma(X,\O_X)$.  A stronger form of
this statement holds if $X$ is essentially of finite type over $\Z$ or over
$\Z_p$ for some prime $p$.  For example, if $X$ is of finite type over $\Q$ or
over $\Q_p$, then the torsion subgroup of $F(X)$ is finite.

Assuming that $X$ is reduced and that the canonical map \mapx[[ \nor{X} || X ]]
is finite, consider the quotient sheaf $F = \O_{\nor{X}}^*/\O_X^*$ on $X$.  We
extend $F$ to an (abelian group)-valued $X$-functor, also denoted here by $F$.
We do not know if $F$ is coherent, but we are able (more or less) to find an
(abelian group)-valued coherent $X$-functor $G$, and a morphism
\mp[[ \psi || F || G ]] such that $\psi(X)$ is injective.
We say ``more or less'' because the actual proof works via a
sequence of partial normalizations.  The end result however is the same:
there are only finitely many primes $p$ such that ${}_p F(X)$ is infinite,
and none of these primes are invertible in $\Gamma(X,\O_X)$.  There is also a
stronger form for certain $X$ as discussed in the previous paragraph.

It follows that if the group $Q = \Gamma(\O_{\nor{X}}^*)/\Gamma(\O_X^*)$ is
\FG, and if $K = \Ker[\Pic(X)\ \mapE{}\ \Pic(\nor{X})]$, then there are only
finitely many primes $p$ such that ${}_p K$ is infinite, and none of these
primes are invertible in $\Gamma(X,\O_X)$.  This holds for instance if $X$ is
of finite type over $\Z$, thereby yielding one of the results stated in the
summary.

When $X$ is of finite type over a field $k$, it has been shown
\Lcitemark 21\Rcitemark \Rspace{} that ${}_n \Pic(X)$ is finite for every $n$
which is invertible in $k$.  For
a finite field $k$ of characteristic $p$, we prove a strengthened form of this
statement: modulo $p$-power torsion, $\Pic(X)$ is \FG.  We completely describe
the structure of $\Pic(X)$ as an abstract abelian group.

The following related result is relevant.  Claborn\Lspace \Lcitemark
12\Rcitemark \Rspace{} has shown that
every abelian group occurs as the Picard group of some Dedekind domain over
$\Q$.  (See also\Lspace \Lcitemark 15\Rcitemark \Rspace{}\ \S14.)  In
particular, for suitable $X$,
$\Pic(X)$ itself has infinite $n$-torsion, for every $n$.

It would be interesting to know to what extent the results of this paper
on $\Pic(X)$ can be obtained via \'etale cohomology.

\vspace{0.1in}
\par\noindent{\footnotesize{\it Acknowledgements.}  I thank Deligne
and Ogus for much help on \S\ref{higher-section}.}

\vspace*{0.1in}

\par\noindent{\bf Conventions}
\
\begin{itemize}
\item $A$ denotes an arbitrary commutative noetherian ring (unless specified
otherwise); $B$ usually denotes an arbitrary commutative $A$-algebra;
\item $X$ denotes an arbitrary noetherian scheme;
\item By a {\it Zariski sheaf}, we mean a sheaf for the Zariski topology.
\item If $S$, $T$ are sets, and $M$ is an abelian group, by a
{\it left exact sequence}
\diagramx{S&\rightarrowtail&T&\mapE{f}&M\cr%
}we mean that the map from $S$ to $T$ is injective, and that $S$ is the
kernel of the map from $T$ to $M$, meaning that
$S = \setof{t \in T: f(t) = 0}$.  Similar language applies when $S$, $T$,
and $M$ are functors.
\end{itemize}
\block{Coherent functors}

In this section we develop the basic theory of coherent functors.
We have already defined the notion of {\it $A$-functor\/} in the introduction.
These form a category \cat{$A$-functors} whose morphisms
are natural transformations.

If $M$ is an $A$-module, then there is an $A$-functor $\uM$ given by
$\uM(B) = M \o*_A B$.  If a given $A$-functor $F$ is isomorphic to $\uM$
for some \FG\ $A$-module $M$, we shall say that $F$ is {\it strictly coherent}.

\begin{definition}
Let $\shC$ be a category.  Let $S$ be a collection of objects in $\shC$.
Let $S_0 = S$, and for each $n \geq 0$, let $S_{n+1}$ be the collection of
all objects of $\shC$, which may be obtained as limits (in $\shC$) of diagrams
involving finitely many objects in $S_n$ and finitely many morphisms.  Let
$S_\infty = \cup_{n=0}^\infty S_n$.  Then we say that the objects in $S_\infty$
are {\it iterated finite limits\/} of objects in $S$.
\end{definition}

\begin{definition}
An $A$-functor is {\it coherent\/} if it may be obtained as an iterated finite
limit of strictly coherent $A$-functors, where the limits are taken
in \cat{$A$-functors}.
\end{definition}

We may define $\cat{coherent $A$-functors}$: it is a full subcategory of
\cat{$A$-functors}, which may be thought of as the
{\it finite completion\/} of the subcategory
\cat{strictly coherent $A$-functors} of \cat{$A$-functors}.

Let $\shC_A^0$ denote the collection of strictly coherent
$A$-functors.  For each $n \geq 0$, let $\shC_A^{n+1}$ denote the
collection of $A$-functors which may be obtained as finite limits (inside
\cat{$A$-functors}) of objects in $\shC_A^n$.  We have:
$$\shC_A^0 \IN \shC_A^1 \IN \shC_A^2 \IN \cdots.$%
$Let $\shC_A = \cup_{n=0}^\infty \shC_A^n$.  Then the objects in
$\shC_A$ are exactly the coherent $A$-functors.

\begin{definition}
Let $F$ be a coherent $A$-functor.  Then the {\it level\/} of $F$ is the
smallest integer $n$ such that $F \in \shC_A^n$.
\end{definition}

In some proofs, we will need to induct on the level of a given coherent
$A$-functor.  This process will be facilitated by the following lemma, whose
proof is left to the reader:

\begin{lemma}\label{lesx-exists}
Let $F$ be a coherent $A$-functor of level $n \geq 1$.  Then there exists
a left exact sequence:
\lesx{F}{G}{\uM%
}in which $G$ is a coherent $A$-functor of level $n-1$ and $M$ is a
\FG\ $A$-module.  Moreover, $F$ may be embedded as a
subfunctor of a strictly coherent $A$-functor.
\end{lemma}

If $k$ is a field, then every coherent $k$-functor is representable, and from
this one deduces easily that every coherent $k$-functor has level $\leq 1$.
Later \pref{main-theorem} we shall prove that a large class of coherent
$k$-functors have level $\leq 1$.  However, we do not know the answer to the
following basic question:

\begin{problemx}
Does every coherent $A$-functor have level $\leq 1$?
\end{problemx}

Note that for a given $A$-functor $F$, this is the case
\IFF\ there exist \FG\ $A$-modules $M$ and $N$,
together with a morphism \mp[[ \phi || \uM || \uN ]] of $A$-functors such that
$F$ is the ``kernel'' of $\phi$, meaning that
$$F(B) = \setof{x \in M \o*_A B: \phi(x) = 0}.$%
$The maps $\phi(B)$ need not be homomorphisms of $B$-modules.

In some situations, for a given coherent $A$-functor $F$, it will be necessary
to consider the set $\shS = \setof{\vec M1n}$ of all $A$-modules which enter
into its construction.  This set is not uniquely determined by $F$.  Also it
does not carry information about multiplicity: $\shS$ might consist of a single
module $M$, but many copies of $M$ might enter into the construction of $F$.
If $F$ has level $0$, we can choose $\shS$ to have one element.  If $F$ has
level $1$, then we may view $F$ as a limit of strictly coherent $A$-functors,
and thus we may choose $\shS$ to consist of the corresponding modules.  If $F$
has level $2$, then $F$ is a limit of level $1$ coherent $A$-functors
$\vec F1k$, and we may choose $\shS$ to be the union of the sets corresponding
(as just considered) to $\vec F1k$.  In any case, we shall say that $F$ is
{\it built up from\/} $\vec M1n$.  Conversely, given an arbitrary class $\shS$
of \FG\ $A$-modules, we may speak of coherent $A$-functors which are
{\it built up from\/} $\shS$, meaning that such $A$-functors are built up from
finite subsets of $\shS$.
\block{Module-coherent functors}

In this section we develop the basic theory of module-coherent functors.
The definition given initially will not be the same as that given in
the introduction.  It is only after considerable work that we will find
\pref{main-theorem} that the two definitions agree.

We have already defined the notion of {\it \MV\/} $A$-functor
in the introduction.  If $F$ and $G$ are \MV\ $A$-functors, then
a {\it morphism\/} \mp[[ \sigma || F || G ]] is a natural transformation of
functors, in the following sense.  It is a system of homomorphisms
\mp[[ \sigma(B) || F(B) || G(B) ]] of $B$-modules, for each commutative
$A$-algebra $B$, such that for any $A$-algebra homomorphism
\mp[[ f || B_1 || B_2 ]], the diagram:
\diagramx{F(B_1)&\mapE{\sigma(B_1)}&G(B_1)\cr
          \mapS{F(f)}&&\mapS{G(f)}\cr
          F(B_2)&\mapE{\sigma(B_2)}&G(B_2)\cr%
}commutes.

With this definition of morphism, the \MV\ $A$-functors form a
category, which is abelian.  Kernels and cokernels are computed in the
obvious way; if \mp[[ \sigma || F || G ]] is a morphism, we have:
$$[\ker(\sigma)](B)\ =\ \ker(\sigma(B))
                   \ =\ \setof{x \in F(B): \sigma(x) = 0},$%
$$$[\coker(\sigma)](B)\ =\ \coker(\sigma(B))\ =\ G(B)/\sigma(F(B)).$%
$One sees that $\sigma$ is a {\it monomorphism\/} \IFF\ $\sigma(B)$ is
injective for every $B$, and $\sigma$ is an {\it epimorphism\/}
\IFF\ $\sigma(B)$ is surjective for every $B$.

Evidently, any \MV\ $A$-functor may be viewed also as an $A$-functor.  In some
situations we shall want to consider \mp[[ \phi || F || G ]] in which
$F$ and $G$ are \MV\ $A$-functors but $\phi$ is a morphism of $A$-functors,
not necessarily preserving the module structure.  For clarity, we may say that
$\phi$ is {\it linear}, if we wish to assume that it is a morphism of
\MV\ $A$-functors.  In this section, all morphisms are linear.

If $M$ is an $A$-module, then $\uM$ is a \MV\ $A$-functor.  If a given
\MV\ $A$-functor $F$ is isomorphic to $\uM$ for some \FG\ $A$-module $M$, we
shall say that $F$ is {\it strictly \MC}.

\begin{definition}
A \MV\ $A$-functor is {\it \MC\/} if it may be obtained as an iterated finite
limit of strictly \MC\ $A$-functors.
These limits are all taken in \cat{\MV\ $A$-functors}.
\end{definition}

We will show \pref{main-theorem}, that in fact this definition is equivalent
to the (much simpler) definition of \MC\ given in the introduction.

\begin{problemx}
If a \MV\ $A$-functor $F$ is coherent (when thought of simply as an
$A$-functor), does it follow that $F$ is \MC?
\end{problemx}

We may define $\cat{\MC\ $A$-functors}$: it is a full subcategory of
\cat{\MV\ $A$-functors}, which may be thought of as the
{\it finite completion\/} of the subcategory \cat{strictly \MC\ $A$-functors}
of \cat{\MV\ $A$-functors}.

Let $\shM\shC_A^0$ denote the collection of strictly
\MC\ $A$-functors.  For each $n \geq 0$, let $\shM\shC_A^{n+1}$ denote the
collection of \MV\ $A$-functors which may be obtained as finite limits (inside
\cat{\MV\ $A$-functors}) of objects in $\shM\shC_A^n$.  We have:
$$\shM\shC_A^0 \IN \shM\shC_A^1 \IN \shM\shC_A^2 \IN \cdots.$%
$Let $\shM\shC_A = \cup_{n=0}^\infty \shM\shC_A^n$.  Then the objects in
$\shM\shC_A$ are
exactly the \MC\ $A$-functors.

\begin{definition}
Let $F$ be a \MC\ $A$-functor.  Then the {\it level\/} of $F$ is the
smallest integer $n$ such that $F \in \shM\shC_A^n$.
\end{definition}

To show that our definition of \MC\ is equivalent to the definition given
in the introduction, we will show \pref{main-theorem} that the level of
$F$ is always $\leq 1$.  In the meantime, however, we will employ induction
on the level of a given \MC\ $A$-functor.  For this we use the
following analog of \pref{lesx-exists}, whose proof is left to the reader:

\begin{lemma}\label{les-exists}
Let $F$ be a \MC\ $A$-functor of level $n \geq 1$.  Then there exists
a left exact sequence:
\les{F}{G}{\uM%
}in which $G$ is a \MC\ $A$-functor of level $n-1$ and $M$ is a \FG\
$A$-module.
Moreover, $F$ may be embedded as a
sub-module-valued functor of a strictly \MC\ $A$-functor.
\end{lemma}

There is a functor $i_A$ from \cat{\FG\ $A$-modules} to
\snap\hbox{\cat{\MC\ $A$-functors}}, given by $M \mapsto \uM$.
It is easily seen that $i_A$ is fully
faithful, so we may view module-coherent $A$-functors as a sort of
generalization of \FG\ $A$-modules.  The functor $i_A$ is cocontinuous:
it preserves colimits.  However, $i_A$ does not carry
monomorphisms to monomorphisms and is not continuous: it does not preserve
limits.  For example, if $J \IN A$ is an ideal, then $J$ is the kernel of
the canonical map \mapx[[ A || A/J ]] of modules, but $\uJ$ is not the kernel
of the induced map \mp[[ \phi || \uA || \ul{A/J} ]].  The kernel of $\phi$ is
instead given by $B \mapsto JB$.

Let $H$ be a \MC\ $A$-functor.  It would be very convenient if
there existed an epimorphism \mapx[[ \ul{A^n} || H ]] for some $n$.
Unfortunately, this is not always the case.  For example, it is not the case
if $H(B) = \Ann_B(x)$, where $A = \C[x]$.  (The module-coherence of this
functor follows from example \pref{annihilator} of \S\ref{examples-section}.)
As a compromise, we are lead to the following notion:

\begin{definition}
A \MV\  $A$-functor $F$ is {\it linearly representable\/} if there
exists a left exact sequence:
\les{F}{\ul{A^n}}{\ul{A^k}%
}in \cat{\MV\ $A$-functors}, for some $n, k \geq 0$.
\end{definition}

If $F$ is linearly representable, it is \MC, and it is
representable by an $A$-algebra of the form
$$A[\vec x1n]/(\vec f1k),$%
$where $\vec f1k$ are linear and homogeneous in $\vec x1n$.
The following proposition is a basic tool, because it exhibits any
\MC\ $A$-functor as a quotient of ``something simple''.

\begin{prop}\label{dog-eats-dog}
Let $F$ be a \MC\ $A$-functor.  Then there exists a linearly
representable $A$-functor $R$ and an epimorphism \mapx[[ R || F ]].
\end{prop}

There are some preliminaries.

\begin{lemma}\label{gerbil}
Let $F$ be a linearly representable $A$-functor, represented by
$C = A[\vec x1n]/(\vec f1k)$, where $\vec f1k$ are linear and homogeneous.
We assume that the module structure on $F$ is the canonical one, induced
from the embedding in $\ul{A^n}$ defined by $\vec x1n$.  Let $C_1$ denote the
degree $1$ part of $C$.  Let $N$ be an $A$-module.  Then morphisms from
$F$ to $\uN$ are in bijective correspondence with elements of $N \o*_A C_1$.
\end{lemma}

\begin{proof}
We can think of $F$ and $\uN$ as functors from the category of commutative
$A$-algebras
to \cat{sets}.  If we take this point of view, then some morphisms
(i.e.\ natural transformations) from $F$ to $\uN$ define morphisms in
\cat{module-valued $A$-functors}, and some do not.  Those which do will be
called {\it linear}, for purposes of this proof.  The linear morphisms are
those which preserve the module structure.

If we take this point of view, then morphisms from $F$ to $\uN$ are in
bijective correspondence with elements of $\uN(C) = N \o*_A C$, and it is clear
that the elements of $N \o*_A C_1$ define linear morphisms.  To complete the
proof, we must show that if an element $\eta \in N \o*_A C$ corresponds to a
linear morphism, then $\eta \in N \o*_A C_1$.

The grading of $C$ induces a grading of $N \o*_A C$.
Let $\eta_1 \in N \o*_A C_1$ denote the degree $1$ part of $\eta$.  Let
$\eta_0 = \eta - \eta_1$.  Then the degree $1$ part of $\eta_0$ is $0$ and
$\eta_0$ defines a linear morphism \mp[[ \psi || F || \uN ]].  We must show
that $\psi = 0$.

Let $B$ be a commutative $A$-algebra.
Let $D = B[t]$.  Then $\psi(dx) = d\psi(x)$ for all $d \in D$ and all
$x \in F(D)$.  In particular, $\psi(tx) = t\psi(x)$ for all $x \in F(B)$.  We
have $\psi(x) \in N \o*_A B$, $\psi(tx) \in N \o*_A D = (N \o*_A B)[t]$.  Since
$\psi(tx) = t\psi(x)$, $\psi(tx)$ is a homogeneous linear polynomial in $t$.
The element $tx \in F(D)$ defines a ring homomorphism \mp[[ \rho || C || D ]],
which maps each generator $x_i$ of $C$ to a homogeneous linear polynomial
in $t$.  The map $\rho$ induces a map \mapx[[ N \o*_A C || N \o*_A D ]],
which sends $\eta_0$ to $\psi(tx)$.  But $\eta_0$ has no linear part, so it
follows that $\psi(tx)$ has no linear part.  Hence $\psi(tx) = 0$.  Hence
$t\psi(x) = 0$, so $\psi(x) = 0$.  Hence $\psi = 0$.  \qed
\end{proof}

\begin{prop}\label{make-arrow}
Let $P$ be an $A$-module.  Suppose given a diagram:
\diagramx{\rowfour{0}{F}{\ul{A^n}}{\ul{A^k}}\cr
          &&\mapS{\phi}\cr
          &&\uP\cr%
}in \cat{\MV\ $A$-functors}, with the row exact.  Then there exists a
morphism \mp[[ h || \ul{A^n} || \uP ]] which makes the diagram commute.
\end{prop}

\begin{proof}
Let $C = A[\vec x1n]$,
$$\oC = A[\vec x1n]/(\vec f1k),$%
$where $\vec f1k$ are homogeneous linear elements determined by the given map
from $\ul{A^n}$ to $\ul{A^k}$.  Then $F$ represents $\oC$.  According to
\pref{gerbil}, $\phi$ corresponds to an element of $P \o*_A \oC_1$.  The
canonical map \mapx[[ P \o*_A C_1 || P \o*_A \oC_1 ]] is surjective, so $h$
exists.  \qed
\end{proof}

\begin{corollary}\label{limit-of-linearly-rep}
Let $\shD$ be a finite diagram in \cat{module-valued $A$-functors}, in
which the objects are linearly representable.  Then the limit of $\shD$
is linearly representable.
\end{corollary}

\begin{proof}
It is clear that a product of finitely many linearly representable
$A$-functors is linearly representable.  Since any linearly representable
$A$-functor embeds in $\ul{A^r}$ for some $r$, we may reduce to showing that
if $F$ is linearly representable and \mp[[ \varphi || F || \ul{A^r} ]] is
a morphism, then $\ker(\varphi)$ is linearly representable.  Let
\lesmaps{F}{}{\ul{A^n}}{h}{\ul{A^k}%
}be as in the definition of {\it linearly representable}.  By
\pref{make-arrow},
$\varphi$ extends to a morphism \mp[[ \psi || \ul{A^n} || \ul{A^r} ]].
Hence $\ker(\varphi) = \ker(\psi) \cap \ker(h)$,
so $\ker(\varphi)$ is linearly representable.  \qed
\end{proof}

\begin{proofnodot}
(of \ref{dog-eats-dog}.)
For purposes of the proof, let us say that a \MV\ $A$-functor $G$
{\it dominates\/} a \MV\ $A$-functor $F$ if there exists an epimorphism
\mapx[[ G || F ]], and that a \MV\ $A$-functor $F$ is
{\it linearly-affine-dominated\/}
if it is dominated by a linearly representable $A$-functor.

Let $n$ be the level of $F$.  The case $n = 0$ is clear -- in that case $F$
is dominated by $\ul{A^r}$ for some $r$.

Suppose that $n \geq 1$.  By \pref{les-exists}, we may find a left exact
sequence:
\les{F}{G}{\uM%
}in which $G$ is a \MC\ $A$-functor of level $n-1$ and $M$ is a
\FG\ $A$-module.  By induction on $n$, \WMAT\ there exists a linearly
representable $A$-functor $H$ and an epimorphism \mapx[[ H || G ]].  Let $P$ be
the fiber product of $H$ with $F$ over $G$.  Then $P$ dominates $F$, so it
suffices to show that $P$ is linearly-affine-dominated.

We have a left exact sequence:
\lesdot{P}{H}{\uM%
}Choose an epimorphism \mp[[ h || A^r || M ]].  Let $L$ be the fiber product
of $H$ with $\ul{A^r}$ over $\uM$.  Let $Q$ be the fiber product of $L$ with
$P$ over $H$.  We have a diagram with cartesian squares, in which some maps
are labelled:
\diagramx{Q&\rightarrowtail&L&\mapE{\sigma}&\ul{A^r}\cr
          \mapS{}&&\mapS{}&&\mapS{\pi}\cr
          P&\rightarrowtail&H&\mapE{f}&\uM\nulldot\cr%
}The bottom row (but not the top) is exact.  The vertical arrows are all
epimorphisms.  Since $Q$ dominates $P$, it suffices
to show that $Q$ is linearly-affine-dominated.

Let $K = \ker(\pi)$.  Then $Q = \sigma^{-1}(K)$, so we have a cartesian diagram
\squareSE{Q}{K}{L}{\ul{A^r}\nulldot%
}We will show that $K$ and $L$ are linearly-affine-dominated.  It will follow
(by taking suitable fiber products) that $Q$ is dominated by a fiber product
of linearly-representable $A$-functors.  By \pref{limit-of-linearly-rep}, it
will follow that $Q$ is linearly-affine-dominated.

There is a canonical epimorphism \mapx[[ \ul{\ker(h)} || K ]].
Choose an epimorphism \mapx[[ A^s || \ker(h) ]].  Then $K$ is dominated by
$\ul{A^s}$.  To complete the proof, we will show that $L$ is
linearly-affine-dominated.

By \pref{gerbil}, it follows that $f$ factors through $\pi$; let
\mp[[ g || H || \ul{A^r} ]] be such that $f = \pi \circ g$.  Then
\begin{eqnarray*}
L(B) & = & \setof{(x,y) \in H(B) \times \ul{A^r}(B): f(x) = \pi(y)}\\
     & = & \setof{(x,y) \in H(B) \times \ul{A^r}(B): \pi(g(x) - y) = 0}.
\end{eqnarray*}
Hence the morphism of \MV\ $A$-functors \mapx[[ H \times K || L ]] given by
$$(x,y) \mapsto (x,g(x)+y)$%
$is an isomorphism.  Since $K$ is linearly-affine-dominated, and $H$ is
linearly representable, it follows that $L$ is linearly-affine-dominated.  \qed
\end{proofnodot}

\begin{theorem}\label{coherence-of-cokernel}
Let \mp[[ \varphi || F || G ]] be a morphism of \MC\ $A$-funct\-ors.  Then
$\Coker(\varphi)$ is \MC.
\end{theorem}

\begin{proof}
Let $n$ be the level of $G$.
First suppose that $n = 0$, so \WMAT\ $G = \uM$ for some \FG\ $A$-module $M$.
Choose an epimorphism \hbox{\mapx[[ A^m || M ]]} of $A$-modules and thus an
epimorphism \mapx[[ \ul{A^m} || \uM ]].
Let $F'$ be the fiber product of $F$ with $\ul{A^m}$ over $\uM$.  Let
$C = \Coker(\varphi)$.  Then we have a right exact sequence:
\resmapsdot{F'}{f}{\ul{A^m}}{}{C%
}By \pref{dog-eats-dog}, \WMAT\ $F'$ is linearly representable.  Choose a
left exact sequence:
\lesdot{F'}{\ul{A^n}}{\ul{A^k}%
}By \pref{make-arrow}, there exists a morphism \mapx[[ \ul{A^n} || \ul{A^m} ]]
which makes the following diagram commute:
% xypic used here
$$\diagram F' \rto^f \dto & \ul{A^m} \rto & C \rto & 0 \\
           \ul{A^n} \dto \urto \\
           \ul{A^k}\nulldot \enddiagram$%
$Let $D$ be the co-fiber product of $A^m$ and $A^k$ over $A^n$, computed in
the category of $A$-modules.  Then in fact
$\uD$ is the co-fiber product of $\ul{A^m}$ and $\ul{A^k}$ over $\ul{A^n}$.
Let \mp[[ g || \ul{A^m} || \uD ]] be the canonical map.  Then
$\ker(g) = \Im(f)$, so $C$ is isomorphic to $\Im(g)$, which is module-coherent
by example \pref{image} from \S\ref{examples-section}.  This completes the
case $n = 0$.

Now suppose that $n \geq 1$.  By \pref{les-exists}, we may choose a left exact
sequence:
\les{G}{H}{\uM%
}in which $H$ is a \MC\ $A$-functor of level $n-1$ and $M$ is a \FG\
$A$-module.
Abusing notation slightly, we have a left exact sequence:
\lesdot{G/F}{H/F}{\uM%
}By induction on $n$, \WMAT\ $H/F$ is \MC.  But then $\Coker(\varphi) = G/F$ is
exhibited as the kernel of a morphism of \MC\ $A$-functors, so it too is \MC.
\qed
\end{proof}

\begin{corollary}\label{result-A}
If \mp[[ \phi || F || G ]] is a morphism of \MC\ $A$-functors,
then $\Im(\phi)$ is \MC.
\end{corollary}

\begin{proof}
Let $K = \Ker(\phi)$.  Then $K$ is \MC.  Hence \snap $\Coker[K \mapE{} F]$ is
\MC\ by \pref{coherence-of-cokernel}, but this equals $\Im(\phi)$.  \qed
\end{proof}

The next result says that the definition of \MC\ given in this section
coincides with the simpler definition given in the introduction.

\begin{corollary}\label{main-theorem}
Let $F$ be a \MC\ $A$-functor.  Then $F$ has level $\leq 1$.
\end{corollary}

\begin{proof}
By \pref{les-exists}, we may embed $F$ as a subfunctor of $\uM$, for some
\FG\ $A$-module $M$.  Let $Q = \uM/F$.  By \pref{coherence-of-cokernel},
$Q$ is \MC.   By \pref{les-exists}, we may embed $Q$ as a subfunctor of $\uN$,
for some \FG\ $A$-module $N$.  Hence $F$ is the kernel of a morphism from
$\uM$ to $\uN$.  \qed
\end{proof}
\block{Quasi-coherent functors}

In this section we sketch a theory (parallel to the last two sections) of
\QC\ and \MQC\ $A$-functors.  The results of this section will be used
in \S\ref{global-section}.  In particular, it is the case that
\QC\ $A$-functors (which are not coherent) are useful in the study of coherent
$A$-functors.  However, the reader interested only in the Picard
group results may ignore this section and everything from \pref{locally-mc}
to the end of \S\ref{global-section}.  The reason for this is explained in
the paragraph preceding \pref{locally-mc}.

If a given $A$-functor is isomorphic to $\uM$ for some $A$-module $M$, we
shall say that $F$ is {\it strictly \QC}.  Similarly,
if a given \MV\ $A$-functor $F$ is isomorphic to $\uM$
for some $A$-module $M$, we shall say that $F$ is {\it strictly \MQC}.

An $A$-functor is {\it \QC\/} if it may be obtained as an iterated finite limit
of strictly coherent $A$-functors.  These limits are all taken in
\cat{$A$-functors}.  We shall not have much more to say about \QC\ $A$-functors
per se in this paper.

A \MV\ $A$-functor is {\it \MQC\/} if it may be obtained as an iterated finite
limit of strictly \MQC\ $A$-functors.  These limits are all taken in
\cat{\MV\ $A$-functors}.

The rest of this section is about \MQC\/ $A$-functors.  All morphisms will be
linear.

We may define the {\it level\/} of a \MQC\ $A$-functor, as we have done for
\MC\ $A$-functors.  As we shall see \pref{main-theorem-mqc}, any \MQC\
$A$-functor has level $\leq 1$.
The analog of \pref{les-exists} for \MQC\ $A$-functors is:

\begin{lemma}\label{les-exists-qc}
Let $F$ be a \MQC\ $A$-functor of level $n \geq 1$.  Then there exists
a left exact sequence:
\les{F}{G}{\uM%
}in which $G$ is a \MQC\ $A$-functor of level $n-1$ and $M$ is an $A$-module.
Moreover, $F$ may be embedded as a
sub-module-valued functor of a strictly \MQC\ $A$-functor.
\end{lemma}

We make the following convention: if $N$ is a set, then $A^N$ denotes a
{\it direct sum\/} of copies of $A$, one for each element of $N$.  If
$S$ is a subset of $N$, then we may view $A^S$ as a submodule of $A^N$, and
thence we may view $\ul{A^S}$ as a subfunctor of $\ul{A^N}$.

\begin{prop}\label{make-arrow-generalized}
Let $P$ be an $A$-module.  Let $N$ and $K$ be sets.  Suppose given a diagram:
\diagramx{0&\mapE{}&F&\mapE{}&\ul{A^N}&\mapE{g}&\ul{A^K}\cr
          &&\mapS{\phi}\cr
          &&\uP\cr%
}in \cat{module-valued $A$-functors}, with the row exact.  Then there exists a
morphism \mp[[ h || \ul{A^N} || \uP ]] which makes the diagram commute.
\end{prop}

\begin{proof}
Let $\fin(N)$ denote the collection of finite subsets of $N$.
Then
$$\setof{\ul{A^S}}_{S \in \fin(N)}$%
$forms a directed system of subfunctors of $\ul{A^N}$,
whose union is $\ul{A^N}$.  Let $F_S = F \cap \ul{A^S}$ for each $S$.
Then $\setof{F_S}_{S \in \fin(N)}$ forms a directed system of subfunctors
of $F$, whose union is $F$.

It is clear that $g|_{\ul{A^S}}$ factors through the subfunctor $\ul{A^{S^*}}$
of $\ul{A^K}$, for some finite subset $S^*$ of $K$.  It follows that
$F_S$ is linearly representable by a ring $\ol{C_S} = A[\sets xsS]/I_S$, where
$I_S$ is generated by linear homogeneous elements.  Let $C_S$ be the
polynomial ring $A[\sets xsS]$.

By \pref{gerbil}, $\phi|_{F_S}$ corresponds to an element of
$P \o*_A (\ol{C_S})_1$.  By lifting this element to an element of
$P \o*_A (C_S)_1$, we see that $\phi|_{F_S}$ can be extended to a
morphism \mp[[ h_S || \ul{A^S} || \uP ]].  As $S$ varies, we have to choose
these extensions $h_S$ so that they are compatible with each other.  To do
this is equivalent to showing that the canonical map:
\dmapx[[ \invlim{S \in \fin(N)} P \o*_A (C_S)_1 ||
         \invlim{S \in \fin(N)} P \o*_A (\ol{C_S})_1 ]]%
is surjective.  In general, it is not true that an inverse limit of surjective
module maps is surjective, but (\Lcitemark 4\Rcitemark \ 10.2) it is the case
if
the transition maps in the system of kernels are surjective.  To show this, it
suffices to show that if $S, S' \in \fin(N)$, with $S \IN S'$, then the
canonical map \mapx[[ P \o*_A (I_{S'})_1 || P \o*_A (I_S)_1 ]] is surjective.
This follows from the fact that the canonical map \mapx[[ I_{S'} || I_S ]] is
surjective.  \qed
\end{proof}

\begin{definition}
A module-valued $A$-functor $F$ is {\it linearly quasi-representable\/} if
there exist sets $N$ and $K$ and a left exact sequence:
\les{F}{\ul{A^N}}{\ul{A^K}%
}in \cat{\MV\ $A$-functors}.
\end{definition}

Using \pref{make-arrow-generalized} one may prove the following analog
of \pref{limit-of-linearly-rep}:

\begin{corollary}\label{limit-of-linearly-rep-gen}
Let $\shD$ be a finite diagram in \cat{\MV\ $A$-functors}, in
which the objects are linearly quasi-representable.  Then the limit of $\shD$
is linearly quasi-representable.
\end{corollary}

Evidently, any linear quasi-representable $A$-functor is \MQC.

\begin{prop}\label{rat-eats-rat}
Let $F$ be a \MQC\ $A$-functor.  Then there exists a linearly
quasi-representable $A$-functor $R$ and an epimorphism \mapx[[ R || F ]].
\end{prop}

\begin{sketch}
Take the proof of \pref{dog-eats-dog}, and modify it in the following ways.
In the various places where $A^r$ is written, one has to allow $r$ to be an
arbitrary set.  Do the same with $A^s$.  Change each reference to {\it \MC\/}
to {\it \MQC}.  Drop the assumption that $M$ is \FG.  The construction of $g$
requires the use of \pref{make-arrow-generalized}.  Use
\pref{limit-of-linearly-rep-gen} instead of \pref{limit-of-linearly-rep}.  Use
\pref{les-exists-qc} instead of \pref{les-exists}.  \qed
\end{sketch}

\begin{prop}\label{quasi-coherence-of-cokernel}
Let \mp[[ \varphi || F || G ]] be a morphism of \MQC\ $A$-functors.  Then
$\Coker(\varphi)$ is \MQC.
\end{prop}

\begin{sketch}
Take the proof of \pref{coherence-of-cokernel}, and modify it in the following
ways.  Change each reference to {\it \MC\/} to {\it \MQC}.
Drop the assumption that $M$ is \FG.  In the notations $A^m$, $A^n$, and
$A^k$, one must allow $m$, $n$, and $k$ to be arbitrary sets.  Use
\pref{make-arrow-generalized} instead of \pref{make-arrow}.  Use
\pref{rat-eats-rat} instead of \pref{dog-eats-dog}.  Use
\pref{les-exists-qc} instead of \pref{les-exists}.  Modify reference to example
\pref{image} from \S\ref{examples-section} appropriately.  \qed
\end{sketch}

\begin{corollary}\label{result-A-MQC}
Let \mp[[ \phi || F || G ]] be a morphism of \MQC\ $A$-functors.  Then
$\Im(\phi)$ is \MQC.
\end{corollary}

\begin{corollary}\label{main-theorem-mqc}
Let $F$ be a \MQC\ $A$-functor.  Then $F$ has level $\leq 1$.
\end{corollary}

\begin{lemma}\label{noetherian}
Let $F$ be a \MC\ $A$-functor.  Let $\sets F\lambda\Lambda$ be a system of
subfunctors of $F$, whose union is $F$.  Then $F = F_{\lambda_0}$ for some
$\lambda_0 \in \Lambda$.
\end{lemma}

\begin{proof}
By \pref{dog-eats-dog}, we may find a linearly-representable $A$-functor $H$,
and an epimorphism \mp[[ \pi || H || F ]].  Let
$H_\lambda = \pi^{-1}(F_\lambda)$,
for each $\liL$.  (That is, $H_\lambda$ is the fiber product of $F_\lambda$
with $H$ over $F$.)  Then $H$ is the union of the $H_\lambda$.
Identify $H$ with the functor representing some $A$-algebra $C$:
$H(B) = \Hom_{A-{\op alg}}(C,B)$.  Choose $\lambda_0 \in \Lambda$ such that
$1_C \in H_{\lambda_0}(C)$.  It follows that there exists a natural
transformation of functors \mp[[ s || H || H_{\lambda_0} ]] (possibly not
preserving module structures) such that $i \circ s = 1_H$, where
\mp[[ i || H_{\lambda_0} || H ]] is the inclusion.  Hence $i$ is an
isomorphism.  Hence $H_{\lambda_0} = H$.  Since $\pi$ is an epimorphism, it
follows that $F_{\lambda_0} = F$.  \qed
\end{proof}

\begin{remark}\label{not-coh-example}
If $F$ is \MC\ and $G$ is \MQC, and \mp[[ \varphi || F || G ]] is a
morphism, then $\ker(\varphi)$ need not be \MC.  For an example,
let $A = \Z$, $F = \ul{\Z}$, $G = \ul{\Q}$, and let
$\varphi$ be the map given by $n \mapsto n$.  Indeed,
$$\ker(\varphi)(B)\ =\ \setofh{$b \in B: nb = 0$ for some $n \in \N$}.$%
$Then $\ker(\varphi)$ is the direct limit of its subfunctors
$F_n = \setof{b \in B: nb = 0}$, for $n \geq 1$,
and since $\ker(\varphi) \not= F_n$ for all $n$, it follows from
\pref{noetherian} that $\ker(\varphi)$ is not \MC.
\end{remark}
\block{Examples}\label{examples-section}

First we give some examples, all of which are easily seen to satisfy the
simple definition of \MC\ which we gave in the introduction.

\begin{arabiclist}
\item $B \mapsto M \o*_A B$, where $M$ is a finitely generated $A$-module;
\item\label{kernel}
      $B \mapsto \Ker(f \o*_A B)$, where \mp[[ f || M || N ]] is a homomorphism
      of \FG\ $A$-modules
      \par[We shall denote this $A$-functor by $\uKer(f)$.]
\item $B \mapsto IB$, where $I$ is an ideal of $A$
      \par[Consider the kernel of the map \mapx[[ \uA || \underline{A/I} ]].]
\item\label{image}
      $B \mapsto \Im(f \o*_A B)$, where \mp[[ f || M || N ]] is a homomorphism
      of finitely generated $A$-modules.
      \par[We shall denote this functor by $\uIm(f)$.  Let
      \mp[[ g || N || N/\Im(f) ]] be the canonical map.  Then
      $\uIm(f) = \uKer(g)$, so $\uIm(f)$ is \MC.]
\item\label{hom-example}
      $B \mapsto \Hom_{B-{\op mod}}(M \o*_A B, N \o*_A B)$, where
      $M$ and $N$ are finitely generated $A$-modules
      \par[We shall denote this functor by $\uHom(M, N)$.  To see why it is
      \MC, choose a presentation $A^k \mapE{} A^n \mapE{} M \mapE{} 0$ for
      $M$. Consider the induced map \mp[[ f || \Hom(A^n, N) || \Hom(A^k, N) ]].
      Then $\uHom(M, N)$ is isomorphic to $\uKer(f)$.]
\item\label{end-example}
      $B \mapsto \End_{B-{\op mod}}(M \o*_A B)$, where $M$ is a
      \FG\ $A$-module \par[We shall denote this functor by $\uEnd(M)$.]
\item\label{annihilator}
      $B \mapsto \Ann_B(I)$, where $I \IN A$ is a fixed ideal
      \par[We shall denote this functor by $\uAnn(I)$.  The point is that
      $\Ann_B(I) = \Hom_{B-{\op mod}}(B/I,B)$;
      use example \pref{hom-example}.]
\end{arabiclist}

\par\noindent Now we give some simple examples of coherent $A$-functors.

\begin{arabiclist}
\setcounter{arabicctr}{7}
\item\label{representable-example}
      $B \mapsto \Hom_{A-{\op alg}}(C,B)$, where $C$ is a
      commutative $A$-algebra of finite type
      \par[Choose a presentation $C = A[\vec x1n]/(\vec f1k)$.  Then
      $\vec f1k$ define a morphism of $A$-functors
      \mapx[[ \ul{A^n} || \ul{A^k} ]], whose kernel is the given functor.]
\item $B \mapsto \setof{x \in B: x^2 \in IB}$
      \par[Consider the kernel of the map \mapx[[ \uA || \underline{A/I} ]]
      given by $x \mapsto x^2$.];
\item $B \mapsto \setof{a \in A^n: f(a) = 0}$, where $f \in M[\vec x1n]$
      \par[Consider kernels of maps \mapx[[ \underline{A^n} || \uM ]]; this
      generalizes the preceding example.]
\end{arabiclist}

The coherence of the remaining examples of this section follows without great
difficulty from the tools developed so far.  However, the remaining examples
seem to be deeper, in the sense that their coherence cannot be deduced directly
from the definitions.

Presumably, all of the usual linear algebra operations ($\Hom$, $\o*$,
$\Lambda^n$, $\ldots$) have analogs for \MV\ $A$-functors.  A thorough study
(which we do not give) would include definitions of these operations and an
analysis of which preserve module-coherence.  We restrict our attention to
some special cases.

Let $F$ and $G$ be \MV\ $A$-functors.  Then there is a \MV\ $A$-functor
$F \o* G$, given by $B \mapsto F(B) \o*_B G(B)$.  If $F$ and $G$ are
\MC, one can ask if $F \o* G$ is \MC.  It turns out \pref{tensor-mc} that
this is not the case.  However, there is the following special case:

\begin{prop}\label{tensor-basic}
Let $F$ be a \MC\ $A$-functor, and let $M$ be a \FG\ $A$-module.  Then
$F \o* \uM$ is \MC.
\end{prop}

\begin{proof}
Choose a right exact sequence:
\res{A^k}{A^n}{M%
}of $A$-modules.  We obtain a right exact sequence:
\res{F \o* \ul{A^k}}{F \o* \ul{A^n}}{F \o* \uM%
}of \MV\ $A$-functors, and thence a right exact sequence:
\resdot{F^k}{F^n}{F \o* \uM%
}Since $F^k$ and $F^n$ are \MC, it follows by \pref{coherence-of-cokernel}
that $F \o* \uM$ is \MC.  \qed
\end{proof}

To show that tensor products do not (in general) preserve module-coherence, we
need the following lemma, which will also be used in a counterexample presented
in \S\ref{higher-section}.

\begin{lemma}\label{ding-dong}
Let $A$ be a noetherian local ring of dimension $d$, having maximal ideal
$\lfM$.  Let $F$ be a \MC\ $A$-functor.  Then there exists a constant $c$, such
that for each $n \in \N$, and every ideal $I$ of $A$ with $\lfM^n \IN I$, we
have $\mu[F(A/I)] \leq c n^d$, where $\mu$ gives the minimal number of
generators of an $A$-module.
\end{lemma}

\begin{remark}
Perhaps the bound $c n^d$ can be replaced by $c n^{d-1}$.
\end{remark}

\begin{proofnodot}
(of \ref{ding-dong}.)
We may assume that $F = \uKer(f)$, for some homomorphism \mp[[ f || M || N ]]
of \FG\ $A$-modules.  Let $\lambda$ denote {\it length}.  Choose a surjection
\mapx[[ A^s || M ]].  Since we have
$$\lfM^{n+1}M\ \IN\ \lfM f^{-1}(IN)\ \IN\ f^{-1}(IN) \IN\ M,$%
$it follows that:
\begin{eqnarray*}
\mu[\Ker(f \o*_A A/I)] & = &    \mu[f^{-1}(IN)/IM]\ \leq\ \mu[f^{-1}(IN)]\\
                       &   &    =\ \mu[f^{-1}(IN)/\lfM f^{-1}(IN)]\\
                       &   &    =\ \lambda[f^{-1}(IN)/\lfM f^{-1}(IN)]\\
                       &   &    \leq\ \lambda(M/\lfM^{n+1}M)
                                \ \leq\ \lambda(A^s/\lfM^{n+1} A^s)\\
                       &   &    =\ s\lambda(A/\lfM^{n+1}).
\end{eqnarray*}
The lemma follows from the theory of the Hilbert-Samuel polynomial.  \qed
\end{proofnodot}

Now we show that the tensor product of two \MC\ $A$-functors need not be \MC:

\begin{prop}\label{tensor-mc}
Let $A = \C[[s,t,u]]$, and let $I$ be the ideal $(s)$ of $A$.  Then the
\MV\ $A$-functor $\uAnn(I) \o* \uAnn(I)$ is not \MC.
\end{prop}

\begin{proof}
Let $F$ denote the given functor.  Fix $n \in \N$, and let $B = A/(s,t,u)^n$.
Then a minimal generating set for $\Ann_B(s)$ is
$$\setof{s^i t^j u^k}_{0 \leq i,j,k \leq n-1,\ i+j+k = n-1},$%
$which has cardinality $n(n+1)/2$.
Then $\mu[F(B)] = \mu[\Ann_B(s)]^2 = [n(n+1)/2]^2$.  By \pref{ding-dong},
$F$ is not \MC.  \qed
\end{proof}

Let $M$ and $N$ be \FG\ $A$-modules.  For $n \geq 0$, one can ask if
the functor $\uTor_{\kern1pt n}(M,N)$ given by
$B \mapsto \Tor_n^B(M \o*_A B, N \o*_A B)$ is \MC.  This seems unlikely for
$n \geq 2$ (but we do not have a counterexample).  For $n = 0$,
\pref{tensor-basic} applies.  For $n = 1$, we have:

\begin{prop}
Let $M$ and $N$ be \FG\ $A$-modules.  Then $\uTor_1(M,N)$ is \MC.
\end{prop}

\begin{proof}
Choose an epimorphism \mapx[[ A^n || M ]] and thence a short exact sequence:
\ses{K}{\ul{A^n}}{\uM%
}in which $K$ is \MC.  One obtains a left exact sequence:
\lesdot{\uTor_1(M,N)}{K \o* \uN}{\ul{A^n} \o* \uN%
}The corollary follows now from \pref{tensor-basic}.  \qed
\end{proof}

If $F$ and $G$ are \MV\ $A$-functors, we let $\HOM(F,G)$ denote the
\MV\ $A$-functor given by
$$B \mapsto \Hom_{\kern2pt{\op \MV}\ B-{\op functors}} (F, G),$%
$where $F$ and $G$ may be viewed as \MV\ $B$-functors because any
$B$-algebra is an $A$-algebra.  In a natural way, $\HOM(F,G)$ is itself a
\MV\ $A$-functor.  We let $\END(F)$ denote $\HOM(F,F)$.

\begin{example}\label{uHomHOM}
Let $M$ and $N$ be \FG\ $A$-modules.  Then there is a
canonical isomorphism of \MV\ $A$-functors
\dmapx[[ \uHom(M,N) || \HOM(\uM,\uN). ]]
\end{example}

\begin{prop}\label{HOM-coherent}
Let $F$ and $G$ be \MC\ $A$-functors.  Then \snap $\HOM(F,G)$ is \MC.
\end{prop}

\begin{proof}
First observe that the bifunctor
\dmap[[ \HOM || \opcat{module-valued $A$-functors} \times
      \cat{module-valued $A$-functors} || \cat{module-valued $A$-functors} ]]%
is left exact in both variables.

By \pref{main-theorem}, there is a left exact sequence:
\les{G}{\ul{M_1}}{\ul{M_2}%
}for some \FG\ $A$-modules $M_1$ and $M_2$.  This yields a left exact sequence:
\lesdot{\HOM(F,G)}{\HOM(F,\ul{M_1})}{\HOM(F,\ul{M_2})%
}Thus it suffices to show that $\HOM(F,\ul{M_i})$ is \MC\ for each $i$.
Let $M = M_i$.

It follows from \pref{dog-eats-dog} that there exists a right exact sequence
\res{L_2}{L_1}{F%
}in which $L_1$ and $L_2$ are linearly representable.  We obtain a left exact
sequence:
\lesdot{\HOM(F,\uM)}{\HOM(L_1,\uM)}{\HOM(L_2,\uM)%
}Therefore we may reduce to proving that $\HOM(L,\uM)$ is \MC\ for any linearly
representable $A$-functor $L$.  The $A$-functor $L$ is representable by a ring
$C = A[\vec x1n]/(\vec f1k)$, where $\vec f1k$ are homogeneous linear
polynomials in $\vec x1n$.  From \pref{gerbil}, it follows that
$\HOM(L,\uM) \iso \ul{M \o*_A C_1}$, so $\HOM(L,\uM)$ is \MC.  \qed
\end{proof}

\begin{corollary}
For any \MC\ $A$-functor $F$, $\END(F)$ is \MC.
\end{corollary}

If $F$ is a \MV\ $A$-functor, we let $\AUT(F)$ denote the $A$-functor given by
$$B \mapsto \Aut_{\kern2pt{\op \MV} \ B-{\op functors}} (F),$%
$where $F$ is viewed as a $B$-functor.  Then $\AUT(F)$ is an $A$-functor.

\begin{corollary}
For any \MC\ $A$-functor $F$, $\AUT(F)$ is coherent.
\end{corollary}

\begin{proof}
Take the kernel of the map
\dmapx[[ \END(F) \times \END(F) || \END(F) \times \END(F) ]]%
given by
$(\alpha,\beta) \mapsto (\alpha \circ \beta - \id, \beta \circ \alpha - \id)$.
\qed
\end{proof}

By an {\it algebra-valued\/} $A$-functor $F$, we shall mean a \MV\ $A$-functor
$F$ which has the additional structure of a $B$-algebra on $F(B)$, for each
$B$.  We do not assume that these algebras $F(B)$ are commutative.

If $F$, $G$, and $H$ are \MV\ $A$-functors, then $\BIL(F \times G, H)$ will
denote the functor of bilinear maps from $F \times G$ to $H$, which sends $B$
to the $B$-module consisting of all morphisms of $B$-functors from
$F \times G$ to $H$ which are bilinear.  Then $\BIL(F \times G, H)$ is a
\MV\ $A$-functor.

\begin{corollary}
Let $F$, $G$, and $H$ be \MC\ $A$-functors.  Then
\snap\hbox{$\BIL(F \times G, H)$} is \MC.
\end{corollary}

\begin{proof}
We may identify $\BIL(F \times G, H)$ with $\HOM(F,\HOM(G,H))$.  \qed
\end{proof}

If $F$ and $G$ are algebra-valued $A$-functors, we let $\HOMalg(F,G)$ denote
the $A$-functor given by
$$B \mapsto \Hom_{\kern2pt{\op algebra-valued} \ B-{\op functors}} (F, G),$%
$where $F$ and $G$ are viewed as algebra-valued $B$-functors.  Similarly, we
may define $\AUTalg(F)$.

\begin{corollary}
Let $R$ and $S$ be \MC, algebra-valued $A$-functors.  Then
$\HOMalg(R,S)$ is coherent.
\end{corollary}

\begin{proof}
Let \mp[[ \mu_R || R \times R || R ]] and \mp[[ \mu_S || S \times S || S ]]
denote the multiplication maps.  Then $\HOMalg(R,S)$ is the kernel of the map:
\dmapx[[ \HOM(R,S) || \BIL(R \times R, S) ]]%
given by $f \mapsto [f \circ \mu_R] - [\mu_S \circ (f \times f)]$.  \qed
\end{proof}

\begin{example}
Let $R$ and $S$ be module-finite $A$-algebras (not necessarily commutative).
Then the $A$-functor defined by
$$B \mapsto \Hom_{B-{\op alg}}(R \o*_A B, S \o*_A B)$%
$is coherent.
\end{example}

\begin{corollary}\label{autalg-is-coherent}
For any \MC, algebra-valued $A$-functor $R$, \snap $\AUTalg(R)$ is coherent.
\end{corollary}

\begin{example}
Let $R$ be a module-finite $A$-algebra (not necessarily commutative).  Then the
$A$-functor given by $B \mapsto \Aut_{B-{\op alg}}(R \o*_A B)$
is coherent.
\end{example}
\block{The global case}\label{global-section}

We have defined the notion of $X$-functor in the introduction; these form a
category \cat{$X$-functors}.  If $X = \Spec(A)$, then there is a canonical
equivalence of categories between \cat{$A$-functors which are Zariski sheaves}
and \cat{$X$-functors which are Zariski sheaves}; we can pass back and forth
freely between these two categories.  Similarly, for an arbitrary noetherian
scheme $X$, we may (instead of looking at $X$-functors which are Zariski
sheaves) look at functors which are Zariski sheaves and whose source is
$$\opcat{$X$-schemes which are quasi-compact}.$$

Let \mp[[ \phi || X || Y ]] be a morphism of noetherian schemes.  We
consider pull-back and push-forward of functors:
\begin{itemize}
\item Let $F$ be a $Y$-functor.  Then there is an $X$-functor $\phi^*F$, given
      by $(\phi^*F)(T) = F(T)$ for all $X$-schemes $T$.  Sometimes we will
      write $\phi|_X$ instead of $\phi^*F$, and refer to the
      {\it restriction\/} of $F$ to $X$.
\item Let $G$ be an $X$-functor.  Then there is a $Y$-functor $\phi_*G$, given
      by $(\phi_*G)(S) = F(X \times_Y S)$, for all $Y$-schemes $S$.
\item Let $F$ be a $Y$-functor.  Then the $Y$-functor $\phi_*\phi^*(F)$ is
      given by $S \mapsto F(X \times_Y S)$, for all $Y$-schemes $S$.  Instead
      of writing $\phi_*\phi^*(F)$, we may refer to ``$F|_X$, viewed as a
      $Y$-functor.''
\end{itemize}

\par\noindent Note that both $\phi^*$ and $\phi_*$ are {\it exact\/} functors.

If $Y$ is an $X$-scheme, \mp[[ \pi || Y || X ]] is a morphism of schemes, and
$\shM$ is a \QC\ $\O_X$-module, we let $\shM_Y$ denote $\pi^*\shM$.
For any such $\shM$, there is an $X$-functor $\ushM$ given by
$\ushM(Y) = \Gamma(Y,\shM_Y)$.

An $X$-functor is {\it strictly coherent\/} if it is isomorphic to $\ushM$
for some coherent $\O_X$-module $\shM$.  An $X$-functor is {\it coherent\/}
if it is an iterated finite limit of strictly coherent $X$-functors, where the
limits are taken in \cat{$X$-functors}.  Similarly, one may define
{\it \QC\/} $X$-functors, by allowing $\shM$ to be \QC.

We may define the {\it level\/} of a coherent $X$-functor, exactly as we have
done for coherent $A$-functors.  We may also define the {\it level\/} of a
\QC\ $X$-functor, and it is distantly conceivable that there exists a coherent
$X$-functor whose level is lower when viewed as a \QC\ $X$-functor.

It is very important to note that if $X = \Spec(A)$, then coherent $X$-functors
are essentially the same as coherent $A$-functors.  This follows from the
fact that coherent $X$-functors are sheaves for the Zariski topology.  Indeed
we have:

\begin{prop}\label{sheaf}
Let $F$ be a \QC\ $X$-functor.  Then $F$ is a sheaf for the fpqc topology.
\end{prop}

\begin{proof}
If $F = \ushM$, for some \QC\ $\O_X$-module $\shM$, then the statement is true.
The proposition follows because any limit of sheaves is a sheaf.  \qed
\end{proof}

\begin{definition}
A {\it \MV\ $X$-functor\/} is an (abelian group)-valued $X$-functor $F$,
together with the structure of a $\Gamma(Y,\O_Y)$-module on each set $F(Y)$,
with the property that for each map of $X$-schemes \mapx[[ Y_1 || Y_2 ]], the
induced map \mapx[[ F(Y_2) || F(Y_1) ]] is a homomorphism of
$\Gamma(Y_2, \O_{Y_2})$-modules.
\end{definition}

The \MV\ $X$-functors form an abelian category.  When we have $X = \Spec(A)$,
there is a canonical equivalence of categories between:
$$\cat{\MV\ $A$-functors which are Zariski sheaves}$%
$and
$$\cat{\MV\ $X$-functors which are Zariski sheaves}.$%
$
If $\shM$ is a \QC\ $\O_X$-module, we let $\ushM$ denote the \MV\ $X$-functor
given by $Y \mapsto \Gamma(Y,\shM_Y)$.  A \MV\ $X$-functor $F$ is
{\it strictly \MC\/} if there exists a coherent $\O_X$-module $\shM$ such that
$F \iso \ushM$.

\begin{definition}
A \MV\ $X$-functor is {\it \MC\/} if it may be obtained as an iterated finite
limit of strictly \MC\ $X$-functors.  These limits are all taken in
\cat{\MV\ $X$-functors}.
\end{definition}

In a similar way, one may define \MQC\ $X$-functors.
If $X = \Spec(A)$, then \MC\ $X$-functors are essentially the same as
\MC\ $A$-functors.  For arbitrary $X$, the theory of \MC\ $X$-functors runs
parallel to the theory of \MC\ $A$-functors, but there is one difference.  When
one takes the cokernel of a morphism of \MC\ $X$-functors, it is necessary to
take the associated sheaf (\WRT\ the Zariski topology), in order to obtain a
\MC\ $X$-functor.

The {\it level\/} of a \MC\ (or \MQC) $X$-functor is defined analogously to
the definition of level for a \MC\ $A$-functor.  As in that case, we will find
ultimately that the level is always $\leq 1$.

An $X$-functor $F$ is {\it locally coherent\/} if $F$ is a Zariski sheaf and if
there exists an open cover $\vec U1n$ of $X$ such that $F|_{U_i}$ is a coherent
$U_i$-functor for each $i$.   Similarly, a \MV\ $X$-functor $F$ is
{\it locally \MC\/} if $F$ is a Zariski sheaf and if there exists an
open cover $\vec U1n$ of $X$ such that $F|_{U_i}$ is a \MC\ $U_i$-functor for
each $i$.  We will show shortly that any locally \MC\ $X$-functor is \MC,
assuming that $X$ is separated.
(One can also define {\it locally \QC\/} and {\it locally \MQC\/}
$X$-functors.)  We do not know the answer to the analogous question for locally
coherent $X$-functors:

\begin{conjecture}\label{locally-coherent-conjecture}
Every locally coherent $X$-functor is coherent.
\end{conjecture}

If the conjecture were true, it would follow immediately that if an $X$-functor
$F$ represents an affine $X$-scheme of finite type, then $F$ is coherent.
(This is true if $X$ is affine: see example \pref{representable-example} from
\S\ref{examples-section}.)

More generally, one could ask:

\begin{problemx}
Let \mp[[ \phi || Y || X ]] be a faithfully flat morphism of noetherian
schemes.  Let $F$ be an $X$-functor, which is a sheaf for the fpqc
topology.  Assume that $\phi^*F$ is a coherent $Y$-functor.  Does it follow
that $F$ is a coherent $X$-functor?
\end{problemx}

We now consider push-forward and pull-back of coherent and quasi-coherent
$X$-functors.  These operations also make sense for \MV\ $X$-functors.

\begin{prop}\label{pull}
Let \mp[[ \phi || X || Y ]] be a morphism of noetherian schemes.  Let
$F$ be a coherent $Y$-functor.  Then $\phi^*F$ is coherent.  Similarly,
if $F$ is a \MC\ $Y$-functor, then $\phi^*F$ is \MC.
\end{prop}

\begin{proof}
Suppose that $F$ is a coherent $Y$-functor.  (The parallel case for
\MC\ $Y$-functors is left to the reader.)  Let $n$ be the level of $F$.
First suppose that $n = 0$, so $F \iso \ushM$ for some coherent
$\O_Y$-module $\shM$.  But then $\phi^*F \iso \ul{\phi^*\shM}$, so $\phi^*F$
is coherent.  Now suppose that $n \geq 1$.  By an unstated analog of
\pref{lesx-exists}, there is a left exact sequence:
\lesx{F}{G}{\ushM%
}of $X$-functors in which $G$ is coherent of level $n-1$ and $\shM$ is a
coherent $\O_X$-module.  By induction on $n$, \WMAT\ $\phi^*G$ is coherent.
Since $\phi^*$ is an exact functor, we have a left exact sequence:
\lesxdot{\phi^*F}{\phi^*G}{\phi^*\ushM%
}Hence $\phi^*F$ is coherent.  \qed
\end{proof}

\begin{prop}\label{push}
Let \mp[[ \phi || X || Y ]] be a morphism of noetherian schemes.
Let $F$ be an $X$-functor.
\begin{alphalist}
\item If $F$ is \QC\ and $\phi$ is affine, then $\phi_*F$ is \QC.
\item If $F$ is coherent and $\phi$ is finite, then $\phi_*F$ is coherent.
\item Parallel statements apply if $F$ is a \MV\ $X$-functor.
\end{alphalist}
\end{prop}

\begin{proof}
(a): Let $n$ be the level of $F$.  If $n = 0$, $F = \ushM$ for some
\QC\ $\O_X$-module $\shM$, so
$(\phi_*F)(T) = \Gamma(X \times_Y T, \shM_{X \times_Y T})$ for any $Y$-scheme
$T$.  By (\Lcitemark 20\Rcitemark \ I:9.1.1), it follows that
$(\phi_*F)(T) \iso \Gamma(T, (\phi_*\shM)_T)$.   Hence $\phi_*F$ is coherent.

Now suppose that $n \geq 1$.  By the (unstated) analog of
\pref{lesx-exists} for \QC\ $X$-functors, there is a left exact sequence:
\lesx{F}{G}{\ushM%
}in which $G$ is a \QC\ $X$-functor of level $n-1$ and $\shM$ is a
\QC\ $\O_X$-module.
Then we have a left exact sequence:
\lesx{\phi_*F}{\phi_*G}{\phi_*\ushM%
}By induction on $n$, \WMAT\ $\phi_*G$ is \QC.  By the $n=0$ case, we may
identify $\phi_*\ushM$ with $\ul{\phi_*\shM}$.  Hence $\phi_*F$ is \QC.

Parts (b) and (c) are left to the reader.  \qed
\end{proof}

The next result is key, since it permits us to reduce to the affine case, and
thereby obtain the analogs of the results for \MV\ $A$-functors.  In
particular, it will follow that many examples of $X$-functors are coherent.
However, the reader interested only in the Picard group results may ignore the
next result and its corollaries, since for purposes of the finiteness result
\pref{coherent-implies-linear}, it is sufficient to know that a given
$X$-functor is {\it locally coherent}.

\begin{theorem}\label{locally-mc}
Assume that $X$ is separated.  Let $F$ be a locally
\MQC\ [resp.\ locally \MC\/] $X$-functor.  Then $F$ is \MQC\ [resp.\ \MC].
\end{theorem}

\begin{proof}
In the course of the proof we refer to sheaves, which shall always mean
sheaves for the {\it Zariski topology}.  We work not with $X$-functors, but
with functors whose source is \opcat{$X$-schemes which are quasi-compact}, as
discussed briefly at the beginning of this section.

Let $\vec U1n$ be as in the definition of {\it locally \MQC\/}
(or {\it locally \MC}).
By \pref{pull}, \WMAT\ each $U_i$ is affine.  Since $X$ is separated, it
follows that the open subschemes $U_i \cap U_j$ are affine and that the
inclusions of $U_i$ in $X$ and of $U_i \cap U_j$ in $X$ are affine morphisms.

First we prove the \MQC\ case.  (This will be needed for the \MC\ case.)
Regard $F|_{U_i}$ and $F|_{U_i \cap U_j}$ as \MV\ $X$-functors.  It follows
from \pref{push} and \pref{pull} that these are \MQC.
Because $F$ is a sheaf, we have a left exact sequence:
\les{F}{\prod_{i=1}^n F|_{U_i}}{\prod_{1 \leq i < j \leq n} F|_{U_i \cap U_j}%
}of \MV\ $X$-functors.  Hence $F$ is \MQC.

Now we show that if $F$ is a \MQC\ $X$-functor, then there exists a morphism
\mp[[ \phi || \shM || \shN ]] of \QC\ $\O_X$-modules such that
$F \iso \uKer(\phi)$.  By an unstated analog of \pref{les-exists},
we may embed $F$ as a sub-\MV-functor of
$\ushM$ for some \QC\ $\O_X$-module $\shM$.  Let $G = (\ushM/F)^*$, where
the superscript $*$ denotes sheafification.
By \pref{quasi-coherence-of-cokernel}, it follows that
$G$ is locally \MQC, so (by the first part of the proof) $G$ is \MQC.  Embed
$G$ as a sub-\MV-functor of $\ushN$ for some \QC\ $\O_X$-module $\shN$.
Let \mp[[ \phi || \shM || \shN ]] be the induced map.  Then $F = \uKer(\phi)$,
as required.

Now we begin the proof of the \MC\ case.  By what we have already shown,
\WMAT\ there exists a morphism \mp[[ \phi || \shM || \shN ]] of
\QC\ $\O_X$-modules such that
$F = \uKer(\phi)$.  We will show that there exist coherent sub-$\O_X$-modules
$\shM_0 \IN \shM_1 \IN \shM$ and $\shN_0 \IN \shN$ such that
$\phi(\shM_0) \IN \shN_0$ and such that in the induced diagram:
\diagramx{\ul{\shM_1}\cr
          \mapN{g}\cr
          \ul{\shM_0}&\mapE{f}&\ul{\shN_0}\cr%
}we have $F \iso g[\ker(f)]^*$.

Let us verify that the construction of this data will complete the proof.  We
must show that $g[\ker(f)]^*$ is \MC.  Let $\shP$ denote the co-fiber product
of $\shM_1$ with $\shN_0$ over $\shM_0$, taken in the category of
\QC\ $\O_X$-modules.  Then in fact the induced diagram
\diagramx{\ul{\shM_1}&\mapE{h}&\ushP\cr
          \mapN{g}&&\mapN{}\cr
          \ul{\shM_0}&\mapE{f}&\ul{\shN_0}\cr%
}is cocartesian, if it is viewed as a diagram in
$$\cat{\MV\ $X$-functors which are sheaves}.$%
$Since we have $\ker(h) \iso g[\ker(f)]^*$, the theorem will follow.

It remains to construct the data.  Certainly, for any choice of data,
there is a canonical morphism
\dmap[[ \psi || g[\ker(f)]^* || F ]]%
of \MV\ $X$-functors.

We work on choosing $\shM_0$.  Let $\sets \shM\lambda\Lambda$ be the coherent
sub-$\O_X$-modules of $\shM$.  Let $H_\lambda$ be the sheafified image of the
map \mapx[[ \ul{\shM_\lambda} || \ushM ]].  Then the $H_\lambda$ form a
directed system of \MV\ subfunctors of $\ushM$ (which are sheaves), whose union
is $\ushM$.  (The validity of the last assertion depends on the simplifying
assumption made in the first paragraph of this proof, to the effect that we
work only with quasi-compact $X$-schemes.)
Let $F_\lambda = F \cap H_\lambda$.  Then the $F_\lambda$
form a directed system of \MV\ subfunctors of $F$ (which are sheaves),
whose union is $F$.  Then it follows from \pref{noetherian} that
$F = F_\lambda$ for some $\liL$.  Let $\shM_0 = \shM_\lambda$.  Then
$F$ is contained in the sheafified image of the map
\mapx[[ \ul{\shM_0} || \ushM ]].

Now we work on choosing $\shN_0$.  Let $\sets \shN\lambda\Lambda$ be the
coherent sub-$\O_X$-modules of $\shN$ which contain $\phi(\shM_0)$.  Let
$I_\lambda$ be the sheafified image of the map
\mapx[[ \ker[\ul{\shM_0}\ \mapE{}\ \ul{\shN_\lambda}] || \ushM ]].  Then the
$I_\lambda$ form a directed system of \MV\ subfunctors of $F$ (which are
sheaves).  We will show that the union of the $I_\lambda$ is $F$.  Let
$J_\lambda = \ker[\ul{\shM_0}\ \mapE{}\ \ul{\shN_\lambda}]$.  It suffices to
show that $\ker[\ul{\shM_0}\ \mapE{}\ \ushN]$ is the union of the $J_\lambda$.
Let $Y$ be a quasi-compact $X$-scheme.  Let $\alpha \in \Gamma[(\shM_0)_Y]$,
and assume that $\alpha \mapsto 0$ in $\Gamma(\shN_Y)$.
We must show that there exists some $\shN_\lambda$ such that
$\alpha \mapsto 0$ in $\Gamma[(\shN_\lambda)_Y]$.  By
(\Lcitemark 20\Rcitemark \ I.6.9.9), we know that $\shN$ is the direct limit of
the
$\shN_\lambda$.  It follows that $\shN_Y$ is the direct limit of the
$(\shN_\lambda)_Y$.  This implies the statement about $\alpha \mapsto 0$, and
hence that the union of the $I_\lambda$ is $F$.  Arguing as in the preceding
paragraph, we see that for some $\liL$, we have $I_\lambda = F$.  Let
$\shN_0 = \shN_\lambda$.  It follows now that now matter how we choose
$\shM_1$, the map $\psi$ will be an epimorphism, when viewed as a morphism
in the category of \MV\ $X$-functors which are sheaves.

Now we work on choosing $\shM_1$.  Let $G = \ker(f)$.  Then $G$ is \MC.  Let
$K = \ker[G\ \mapE{}\ \ushM ]$.  Then $K = \ker[G\ \mapE{}\ F ]$, so $K$ is
locally \MC.  Let $\sets \shM\lambda\Lambda$ be
the coherent sub-$\O_X$-modules of $\shM$ which contain $\shM_0$.  (These
$\shM_\lambda$ are not the same as those defined earlier.)  Let
$K_\lambda = \ker[ G\ \mapE{}\ \ul{\shM_\lambda} ]$.  Then the
$K_\lambda$ form a directed system of \MV\ subfunctors of $K$ (which are
sheaves).  Arguing as in the preceding paragraph, we see that the union of
the $K_\lambda$ is $K$.  But $K$ is locally \MC,
so it follows from \pref{noetherian} that $K_\lambda = K$ for some $\liL$.
Let $\shM_1 = \shM_\lambda$.
Then $\psi$ is a monomorphism.  Hence $\psi$ is an isomorphism.  \qed
\end{proof}

\begin{corollary}\label{mump}
Assume that $X$ is separated.
Let \mp[[ \phi || F || G ]] be a morphism of \MC\ $X$-funct\-ors.  Then
the Zariski sheaf associated to $\Coker(\phi)$ is \MC.
\end{corollary}

\begin{corollary}
Assume that $X$ is separated.  Then the category of \MC\ $X$-functors is
abelian.
\end{corollary}

\begin{corollary}
Assume that $X$ is separated.
Let $F$ be a \MC\ $X$-functor.  Then $F$ has level $\leq 1$.  That is,
there exists a left exact sequence:
\les{F}{\ushM}{\ushN%
}in which $\shM$ and $\shN$ are coherent $\O_X$-modules.
\end{corollary}

The constructions $\HOM$, $\AUT$ and so forth which we defined in
\S\ref{examples-section} make sense for $X$-functors.  For example,
if $F$ and $G$ are module-valued $X$-functors, we let
$\HOM(F,G)$ denote the $X$-functor given by
$$Y \mapsto \Hom_{\kern2pt{\op module-valued}\ Y-{\op functors}}
            (F|_Y, G|_Y).$%
$
\begin{prop}
Assume that $X$ is separated.
Let $F$ and $G$ be \MC\ $X$-functors.  Then \snap $\HOM(F,G)$ is \MC.
\end{prop}

\begin{proof}
By \pref{HOM-coherent} and \pref{locally-mc}, it suffices to show that
$\HOM(F,G)$ is a Zariski sheaf.  This can be directly checked from the
definition.  \qed
\end{proof}

Similarly, we have:

\begin{corollary}\label{lump}
Assume that $X$ is separated.  Then
for any \MC\ $X$-functor $F$, $\END(F)$ is \MC, and $\AUT(F)$ is coherent.
Let $R$ and $S$ be \MC, algebra-valued $X$-functors.  Then
$\HOMalg(R,S)$ and $\AUTalg(R)$ are coherent.
\end{corollary}

Now we consider the extent to which a \MQC\ $A$-functor can be viewed as a
direct limit of \MC\ $A$-functors.  These considerations will enter into an
analysis of extensions of \MC\ $X$-functors, which will be the last topic
discussed in this section.

Unfortunately, it is not the case that every \MQC\ $A$-functor $H$ is the
direct limit of its \MC\ subfunctors.  For an example, let $H = \Im(\phi)$,
where $\phi$ is as in remark \pref{not-coh-example}.  If there existed a
directed system $\sets H\lambda\Lambda$ of \MC\ subfunctors of $H$, with union
$H$, it would follow by \pref{noetherian}, applied with $F = \ul{\Z}$,
$F_\lambda = \phi^{-1}(H_\lambda)$, that $\phi$ factors through a
\MC\ subfunctor of $\ul{\Q}$, and hence that $\ker(\phi)$ is \MC, which is not
the case.

\begin{definition}
A \MV\ $A$-functor $C$ is {\it bar-module-coherent\/} if it is \MQC\ and if
there exists a \MC\ $A$-functor $F$, together with an epimorphism
\mapx[[ F || C ]].
\end{definition}

A bar-\MC\ $A$-functor need not be \MC.  For an example, let $H = \Im(\phi)$,
where $\phi$ is as in remark \pref{not-coh-example}.  Then $H$ is bar-\MC, but
not \MC, since otherwise $\ker(\phi)$ would be \MC.

\begin{lemma}\label{directed-system-of-images}
Let $F$ be a \MQC\ $A$-functor.  Then there exists a directed
system $\sets F\lambda\Lambda$ of bar-\MC\ subfunctors of $F$, with union $F$.
\end{lemma}

\begin{proof}
Choose $A$-modules $M$, $N$ and a left exact sequence:
\lesmaps{F}{}{\uM}{h}{\uN%
}of module-valued $A$-functors, in which $h$ is induced by a homomorphism
\mp[[ \phi || M || N ]] of $A$-modules.  Let $\shS$ denote the
collection $\setof{(M_\lambda,N_\lambda)}_{\liL}$ consisting of all pairs
\snap\hbox{$(M_\lambda,N_\lambda)$} in which $M_\lambda$ is a \FG\ submodule of
$M$, $N_\lambda$ is a \FG\ submodule of $N$, and
$\phi(M_\lambda) \IN N_\lambda$. For each $\liL$,
let \mp[[ h_\lambda || \ul{M_\lambda} || \ul{N_{\lambda}} ]] be the
induced morphism of $A$-functors.  Let $K_\lambda = \ker(h_\lambda)$.  Then
$K_\lambda$ is \MC.  There is a canonical map
\mp[[ f_\lambda || K_\lambda || F ]].  Let $F_\lambda = \Im(f_\lambda)$.  Then
$F_\lambda$ is bar-\MC, and the $F_\lambda$ form a directed system of
subfunctors of $F$.

Let $B$ be a commutative $A$-algebra, and let $c \in F(B)$.  Then $c \in M_B$.
Choose a \FG\ submodule $M_\lambda \IN M$ and an element
$c_\lambda \in (M_\lambda)_B$ such that $c_\lambda \mapsto c$.  There exists a
\FG\ submodule $N_\lambda$ of $N$ such that $\phi(M_\lambda) \IN N_\lambda$
and such that $c_\lambda \mapsto 0$ in $(N_\lambda)_B$.  It follows that
$\cup_{\liL}F_\lambda = F$.  \qed
\end{proof}

We close this section with some questions and a result about extensions:

\begin{problemx}
Let
\Rowfive{1}{F'}{F}{F''}{1%
}be a short exact sequence of group-valued $X$-functors.  Assume that $F'$ and
$F''$ are coherent.  Is $F$ coherent?
\end{problemx}

\begin{problemx}
Let
\ses{F'}{F}{F''%
}be a short exact sequence of \MV\ $X$-functors.  Assume that $F'$ and $F''$
are
\MC.  Is $F$ \MC?
\end{problemx}

We can prove this if we assume that $F$ is \MQC\ and that $X$ is separated:

\begin{prop}\label{extension-is-mc}
Assume that $X$ is separated.  Let
\ses{F'}{F}{F''%
}be a short exact sequence of \MV\ $X$-functors.  Assume that $F'$ and $F''$
are \MC.  Assume that $F$ is \MQC.  Then $F$ is \MC.
\end{prop}

\begin{proof%
}By \pref{locally-mc}, we may reduce to working with $A$-functors.  By
\pref{directed-system-of-images}, $F$ is the direct limit of its
bar-\MC\ subfunctors.  By \pref{noetherian}, it follows that there exists a
bar-\MC\ subfunctor $B$ of $F$ such that $B$ maps onto $F''$.  Since $F'$ is
\MC, we see that $F$ is itself bar-\MC.  Choose a \MC\ $A$-functor $C$ and
an epimorphism \mp[[ \beta || C || F ]].

Let $P$ be the fiber product of $F'$ and $C$ over $F$.  Then we have a
commutative diagram with exact rows:
\diagramx{\rowfive{0}{P}{C}{F''}{0}\cr
          &&\mapS{\alpha}&&\mapS{\beta}&&\mapS{=}\cr
          \rowfive{0}{F'}{F}{F''}{0}\nulldot%
}Hence $\ker(\alpha) \iso \ker(\beta)$.  Since $C$ and $F''$ are \MC, so is
$P$.  Since $F'$ and $P$ are \MC, so is $\ker(\alpha)$.  Hence $\ker(\beta)$
is \MC.  Since $\ker(\beta)$ and $C$ are \MC, it follows by
\pref{coherence-of-cokernel} that $F$ is also \MC.  \qed
\end{proof}

This fact will be used in the proof of \pref{mulch-material}.
\block{Continuity}\label{continuity-section}

We consider the extent to which \QC\ $A$-functors preserve limits, and briefly,
the extent to which they preserve colimits.  We prove that a \MQC\ $A$-functor
which preserves products is \MC.  Although these topics do not play
much of a role in the subsequent parts of this paper, they are very natural.
Some important examples of limit and colimit preservation which have
arisen previously are Grothendieck's theorem on formal functions
(see e.g.{\ }\Lcitemark 22\Rcitemark \ III:11), and Grothendieck's
notion of functors which are locally of finite presentation
(\Lcitemark 1\Rcitemark \ 1.5), which enters into Artin's
criterion for representability (\Lcitemark 3\Rcitemark \ 3.4).

We begin by recalling some definitions about continuity of functors.
Let $\shC$ and $\shA$ be complete (meaning small-complete) categories, and let
\mp[[ F || \shC || \shA ]] be any functor.  Then $F$ is {\it continuous\/} if
it {\it preserves limits}, i.e.\ if for every small category $\shD$,
and every functor \mp[[ H || \shD || \shC ]], the canonical map
\dmapx[[ F \left( \inverselim\ H \right) || \inverselim(F \circ H) ]]%
is an isomorphism.  (See e.g.{\ }\Lcitemark 27\Rcitemark \ V.4.)  It is also
of interest to know if $F$ preserves more restricted sorts of limits,
e.g.\ does it preserve arbitrary products, or does it preserve equalizers.
These conditions correspond to placing appropriate restrictions on $\shD$.

For particular sorts of limits, one can usually rephrase the continuity
condition in a simpler way.  For example, $F$ preserves products \IFF\ for
every set \sets XiI of objects in $\shC$, the canonical map
\dmapx[[ F \left( \prod_{i \in I} X_i \right) || \prod_{i \in I} F(X_i) ]]%
is an isomorphism.

Just as one can check a category for completeness by checking if it has
products and equalizers, so one can check a functor for continuity by
checking if it preserves products and equalizers.

We will study the continuity properties of $A$-functors.  There are two general
observations to be made.  The first observation is that an
(abelian group)-valued $A$-functor or a \MV\ $A$-functor is continuous (or
preserves a particular kind of limit) \IFF\ the same statement holds for the
underlying functor from \cat{commutative $A$-algebras} to \cat{sets}.  The
second observation is the following lemma, whose proof is left to the reader:

\begin{lemma}\label{kernel-preserves}
Let $M$ be an $A$-module and let
\lesx{F}{G}{\uM%
}be a left exact sequence of $A$-functors.  If $\uM$ and $G$ preserve a
particular type of limit, then so does $F$.
\end{lemma}

The phrase ``particular type of limit'' is to be construed as referring to
a class of limits which is constrained by some restriction on the categories
$\shD$ and/or the functors \mp[[ H || \shD || \cat{commutative $A$-algebras} ]]
which are allowed.

We proceed to investigate the extent to which \QC\ $A$-functors preserve
various types of limits.  First we consider finite limits, that is limits in
which the category $\shD$ has only finitely many objects and morphisms.
The most important examples are finite products (including terminal objects),
and equalizers.  Also, if a functor preserves finite products and equalizers,
then it preserves all finite limits.  As for finite products, one sees
easily (using \ref{kernel-preserves}) that:

\begin{prop}\label{number-1}
Let $F$ be a \QC\ $A$-functor.  Then $F$ preserves finite products.
\end{prop}

Now we consider equalizers.  An $A$-module $M$ is flat \IFF\ the functor
\fun[[ \o*_A M || $A$-modules || $A$-modules ]] preserves equalizers,
so the analogous fact for $A$-functors is hardly surprising:

\begin{prop}\label{commutes-eq-iff-flat}
Let $M$ be an $A$-module.  Then the $A$-functor $\uM$ preserves
equalizers \IFF\ $M$ is flat.
\end{prop}

\begin{proof}
The case where $M$ is flat is left to the reader.  So suppose that $M$ is not
flat.  Then (see e.g.{\ }\Lcitemark 28\Rcitemark \ 3.53) there exists an
ideal $I \IN A$ such that the induced map \mp[[ \phi || M \o*_A I || M ]] is
not injective.  Let $y \in \ker(\phi) - \setof{0}$.  Let $\vec a1n$ be
generators for $I$.  Choose $\vec m1n \in M$ such that
$y = \sum_{i=1}^n m_i \o* a_i$.  Let
$B = A[\vec x1n]/(\setof{x_i x_j}_{1 \leq i,j \leq n})$.  Let
\mp[[ f,g || B || A ]] be the $A$-algebra maps given by $f(x_i) = a_i$
and $g(x_i) = 0$ for each $i$.  Let $p \in M \o*_A B$ be
$\sum_{i=1}^n m_i \o* x_i$.  Then $(f \o*_A M)(p) = (g \o*_A M)(p)$.  Let
$\Eq(f,g)$ be the equalizer of $f$ and $g$.  Since $y \not= 0$, it follows
(after a little work) that $p$ does not lie in the image of the canonical map
\mp[[ \lambda || \Eq(f,g) \o*_A M || \Eq(f \o*_A M, g \o*_A M) ]], and hence
that $\lambda$ is not an isomorphism.  Hence $\uM$ does not preserve the
equalizer of $f$ and $g$.  \qed
\end{proof}

It follows that $\uM$ preserves finite limits \IFF\ $M$ is flat.  In
addition, \pref{kernel-preserves} yields the following corollary:

\begin{corollary}\label{number-2}
Any \QC\ $A$-functor which is built up from flat modules will preserve
finite limits.
\end{corollary}

It is worth noting that if $F$ is an $A$-functor which preserves finite limits,
then $F$ is a sheaf \WRT\ the fppf topology; indeed the sheaf axioms may be
viewed as a statement about continuity.  While \QC\ $A$-functors do not always
preserve finite limits, we do know that they are sheaves \WRT\ the
fppf topology.

\begin{remark}\label{number-4}
It is not difficult to see that any coherent $A$-functor preserves
products.  Any coherent $A$-functor which is built up from \FG\ projective
$A$-modules will also preserve equalizers, and hence all limits.  Of
course, the $A$-functors which are built up in this way are exactly the
$A$-functors which are representable by an $A$-algebra of finite type.
\end{remark}

In general, coherent $A$-functors do not preserve inverse limits.  For example,
if $A = \Z$, then the $A$-functor $\ul{\Z/2\Z}$ does not preserve the limit of
\diagramx{\Z[x]/(x^2) & \mapW{x\kern2pt\mapsto\kern2pt 3x} & \Z[x]/(x^2)
          & \mapW{x\kern2pt\mapsto\kern2pt 3x} & \cdots.%
}The transition maps here are not surjective.  One might hope that the limit
would be preserved if the transition maps were surjective, or at least if the
system satisfied the Mittag-Leffler condition (\Lcitemark 17\Rcitemark \
0:13.1).
Unfortunately this is not the case:

\begin{example}
Let $A = \Z$.  Let
$$B_n = \Z[x_1, x_2,\ldots, y, z_1,\ldots,z_{n-1}]
        / (2x_1,2x_2,\ldots,Q),$%
$where $Q$ denotes the set of homogeneous quadratic polynomials in all of the
given variables.  Form an inverse system of commutative $A$-algebras
\diagramx{B_1 & \mapW{} & B_2 & \mapW{} & B_3 & \mapW{} & \cdots%
}in which (for each $n > 1$) the transition map \mapx[[ B_n || B_{n-1} ]] is
given by $x_k \mapsto x_{k+1}$, $y \mapsto x_1 + y$, $z_1 \mapsto x_1$,
$z_k \mapsto z_{k-1}$ for $k \geq 2$.  The transition maps of this system are
surjective.

The elements $2y, 2y, \ldots$ form a coherent sequence.  Each element in
this sequence is divisible by $2$, but there is no coherent sequence which when
multiplied by $2$ yields $2y, 2y, \ldots$.  Hence the $A$-functor $\ul{\Z/2\Z}$
does not preserve the limit of this system.
\end{example}

{}From (\Lcitemark 4\Rcitemark \ 10.13) and \pref{kernel-preserves} it follows
that:

\begin{prop}\label{number-5}
Let $F$ be a coherent $A$-functor.  Let $B$ be a commutative noetherian
$A$-algebra.  Let $I \IN B$ be an ideal.  Then $F$ preserves the limit of
\diagramx{B/I & \mapW{} & B/I^2 & \mapW{} & B/I^3 & \cdots.}
\end{prop}

There may well be interesting situations in which coherent functors preserve
inverse limits, other than those given in \pref{number-5} and \pref{number-4}.

% pstricks used here

%The following diagram summarizes what we know about continuity properties
%of certain kinds of \QC\ $A$-functors.  For precise statements, see
%\pref{number-1}, \pref{number-2}, \pref{number-4}, and
%\pref{number-5}.

%\vspace{0.2in}

%\begin{center}
%\par\noindent\kern0.15in\begin{tabular}{c@{\hskip 0.25cm}c@{\hskip 0.25cm}c}
%& \ovalnode[linecolor=gray,linewidth=2pt]{a}{affine representable}
%\makenull{\ \it preserves all limits}\\[2cm]
%\ovalnode[linecolor=gray,linewidth=2pt]{b}{coherent}
%$\vcenter{\makenull{\ \it preserves products,}\kern4pt%
%\makenull{\ \it noetherian completions}}$
%&& \ovalnode[linecolor=gray,linewidth=2pt]{c}{built from flats}
%\makenull{\kern-0.9in\lower0.3in\hbox{\it preserves finite limits}}\\[2cm]
%& \ovalnode[linecolor=gray,linewidth=2pt]{d}{quasi-coherent}
%\makenull{\ \it preserves finite products}
%\end{tabular}

%\ncline[linecolor=gray,linewidth=2pt]{-}{a}{b}
%\ncline[linecolor=gray,linewidth=2pt]{-}{a}{c}
%\ncline[linecolor=gray,linewidth=2pt]{-}{b}{d}
%\ncline[linecolor=gray,linewidth=2pt]{-}{c}{d}
%\end{center}

%\vspace{0.15in}

Our next objective is to show that a \MQC\ $A$-functor which preserves
products is \MC.  This as well as \pref{commutes-eq-iff-flat} allow one to use
knowledge about continuity to deduce some information about how a
\QC\ $A$-functor is built up.

\begin{lemma}\label{fg-sufficient}
Let $L_1$ and $L_2$ be \MQC\ subfunctors of a \MQC\ $A$-functor $H$.  If
$L_1(B) \IN L_2(B)$ for every \FG\ commutative $A$-algebra $B$, then
$L_1 \IN L_2$.
\end{lemma}

\begin{proof}
Let $G = L_1 / (L_1 \cap L_2)$.  Since $G$ is \MQC, we have $G \iso \uKer(f)$
for some homomorphism \mp[[ f || M || N ]] of $A$-modules.  Then
$\Ker(f \o*_A B) = 0$ for every \FG\ commutative $A$-algebra $B$, from which it
follows that $\Ker(f \o*_A B) = 0$ for {\it every\/} commutative $A$-algebra
$B$.  Hence $G = 0$.  \qed
\end{proof}

\begin{lemma}\label{factor-fg}
Let \mp[[ \phi || H || N ]] be a homomorphism of $A$-modules.  Assume that $H$
is \FG.  Let $B$ be a commutative $A$-algebra.  Let $h \in \Ker(\phi \o*_A B)$.
Then there exists a \FG\ submodule $N_0$ of $N$ such that $\phi$ factors
through $N_0$ and such that $h \mapsto 0$ in $N_0 \o*_A B$.
\end{lemma}

\begin{proof}
Certainly $N$ is the direct limit of its \FG\ submodules which contain
$\phi(H)$.  The statement follows from the fact that tensor products commute
with direct limits.  \qed
\end{proof}

\begin{prop}\label{product-preserving}
Let $F$ be a \MQC\ $A$-functor which preserves products.  Then $F$ is \MC.
\end{prop}

\begin{proof}
All tensor products in this proof are over $A$.
We may assume that $F = \uKer(f)$ for some homomorphism \mp[[ f || M || N ]] of
$A$-modules.  We will show that there exists a \FG\ submodule $L \IN M$ such
that if \mp[[ i || L || M ]] is the inclusion, then $\uKer(f) \IN \uIm(i)$.
Choose a complete set of isomorphism class representatives
$\sets B\lambda\Lambda$ for the \FG\ commutative $A$-algebras.  Form the
disjoint union $T = \coprod_{\liL} \Ker(f \o* B_\lambda)$.  For any
$t \in T$, let $\lambda(t)$ denote the corresponding element of $\Lambda$.  Let
$B = \prod_{t \in T} B_{\lambda(t)}$.  The elements $t \in T$ define an element
$x \in \prod_{t \in T} \Ker(f \o* B_{\lambda(t)})$.  Since $F$ preserves
products, $x \in \Ker(f \o* B)$.  Write $x = \sum_{j=1}^r m_j \o* b_j$, where
$\vec m1r \in M$, $\vec b1r \in B$.  Let $L$ be the submodule of $M$ generated
by $\vec m1r$.  The expansion of $x$ defines an element $\ltX \in L \o* B$
with the property that $\ltX \mapsto x$.  From this it follows that
$\Ker(f \o* B_\lambda) \IN \Im(i \o* B_\lambda)$ for every $\liL$.  By
\pref{fg-sufficient}, we have $\uKer(f) \IN \uIm(i)$.

Since $\ltX \mapsto 0$ in $N \o* B$, it follows from \pref{factor-fg} that
there exists a \FG\ submodule $N_0$ of $N$ such that \mapx[[ L || N ]] factors
through $N_0$ and such that $\ltX \mapsto 0$ in $N_0 \o* B$.

For any $t \in T$, let $\lambda = \lambda(t)$, and let
$\ltT \in L \o* B_\lambda$ be the image of $\ltX$ under the canonical map
\mp[[ \pi_t || L \o* B || L \o* B_\lambda ]] which projects onto the \th{t}
factor.  Then $(i \o* B_\lambda)(\ltT) = t$, and $\ltT \mapsto 0$ in
$N_0 \o* B_\lambda$.  Let $\uK = \uKer(L\ \mapE{}\ N_0)$.  Then
$\uK(B_\lambda)$ maps {\it onto\/} $F(B_\lambda)$.  By \pref{fg-sufficient},
the map \mp[[ \psi || \uK || F ]] is an epimorphism.  In particular, $F$ is the
image of a \MC\ $A$-functor.

Let $Q = \Ker(\psi)$.  Since $\uK$ is \MC, it preserves products.  Since $F$
also preserves products, it follows by \pref{kernel-preserves} that $Q$
preserves products.  Replaying the first part of the proof, with $F$ replaced
by $Q$, we see that $Q$ is also the image of a \MC\ $A$-functor.  Hence $F$ is
the cokernel of a map of \MC\ $A$-functors, so $F$ is \MC.  \qed
\end{proof}

The last objective of this section is to consider (briefly) the extent to which
coherent functors preserve colimits.

Whether or not coherent functors preserve finite colimits and coproducts does
not seem to be an interesting question.  One reason for this is that the
forgetful functor from \cat{groups} to \cat{sets} does not preserve coproducts
or finite colimits.  (One can substitute various other categories for
\cat{groups} with the same outcome.)  Therefore a functor from
\cat{commutative $A$-algebras} to \cat{groups} might preserve such colimits,
but the induced functor from \cat{commutative $A$-algebras} to \cat{sets}
might not.  With either interpretation, preservation of coproducts or finite
colimits seems like a bizarre requirement.

On the other hand, the forgetful functor \funx[[ groups || sets ]] does
preserve direct limits\footnote{The reader is reminded that direct limits are
a type of {\it colimit}.}, and the same statement is valid with various other
categories substituted for \cat{groups}.  Therefore, the situation for
direct limits is just like the situation which holds for all limits: an
(abelian group)-valued $A$-functor or a \MV\ $A$-functor preserves direct
limits \IFF\ the same statement holds for the underlying functor from
\cat{commutative $A$-algebras} to \cat{sets}.

The analog of \pref{kernel-preserves} for direct limits is valid, and since
tensor products commute with direct limits, it follows that any
\QC\ $A$-functor preserves direct limits.  Artin remarks that nearly all
$A$-functors which occur in practice do this; in Artin's terminology an
$A$-functor which preserves direct limits is said to be
{\it locally of finite presentation} (\Lcitemark 1\Rcitemark \ 1.5).
This condition enters into his criterion for representability
(\Lcitemark 3\Rcitemark \ 3.4).
\block{Coherence of higher direct images as functors}\label{higher-section}

In this section we consider a question which was posed (in an equivalent
form) by Artin\Lspace \Lcitemark 2\Rcitemark \Rspace{}:

\begin{problem}
Let $X$ be a proper $A$-scheme, let $\shF$ be a coherent sheaf on $X$, and
fix $n \geq 0$.  Is the $A$-functor $H = H^n_\shF$ given by
$B \mapsto H^n(X_B, \shF_B)$ \MC?
\end{problem}

Taking the \u Cech resolution of $\shF$ relative to some affine open cover of
$X$ yields a complex $K$ of $A$-modules, and by (\Lcitemark 17\Rcitemark \
1.4.1) we have
$H^n(X_B, \shF_B) \iso H^n(K \o* B)$ for all commutative $A$-algebras $B$.  It
follows that at least $H$ is {\it \MQC}.  The issue of whether $H$ is \MC\ is
quite subtle.  We will show that if $\shF$ is $A$-flat, or $A$ is a Dedekind
domain, then $H$ is \MC.  By \pref{product-preserving}, we know also that $H$
is \MC\ \IFF\ $H$ preserves products, but it is not clear how to use this
statement.  We will give an example which shows that in general $H$ is not \MC.

Let us say that a \MV\ $A$-functor $F$ is {\it upper semicontinuous\/} if
for every commutative $A$-algebra $B$, and every $\lfP \in \Spec(B)$, there
is a neighborhood $U$ of $\lfP$ such that
$\dim_{k(\lfQ)} F(k(\lfQ)) \leq \dim_{k(\lfP)} F(k(\lfP))$ for all $\lfQ \in
U$.
Similarly, one defines {\it lower semicontinuous\/} by reversing the
inequality.  If $M$ is a \FG\ $A$-module, then it follows by Nakayama's lemma
that $\uM$ is upper semicontinuous.  If \mp[[ f || P_1 || P_2 ]] is a map of
\FG\ projective $A$-modules, then $\uKer(f)$ is upper semicontinuous, whereas
$\uIm(f)$ is lower semicontinuous.  For $I \IN A$ an ideal,
$\uKer(A\ \mapE{}\ A/I)$ is in general not upper semicontinuous.  (But it is
lower semicontinuous.)  If the numerator of
a quotient is upper semicontinuous and the denominator is lower semicontinuous,
then the quotient is itself upper semicontinuous.  It follows that if $K_0$ is
a complex of \FG\ free $A$-modules, and $n \in \Z$, then the \MV\ $A$-functor
given by $B \mapsto H^n(K_0 \o* B)$ is upper semicontinuous.  It is also \MC.

\begin{prop}\label{flat-implies-cohomologically-coherent}
If $X$ is proper over $A$ and $\shF$ is a coherent sheaf on $X$ which is
$A$-flat, then the \MV\ $A$-functor $H^n_\shF$ given by
$B \mapsto H^n(X_B, \shF_B)$ is \MC\ and upper semicontinuous.
\end{prop}

\begin{proofnodot}
(pointed out to me by Deligne and Ogus; cf.{\ }\Lcitemark 23\Rcitemark .)
Let $K$ be the \u Cech complex discussed above.  It is bounded and flat.
Since $X$ is proper over $A$, the complex $K$ has \FG\ cohomology modules.
It follows that there exists a complex $K_0$ of \FG\ free $A$-modules which is
bounded above and a quasi-isomorphism \mp[[ \phi || K_0 || K ]].  Since
$K$ is flat and bounded above, there is a spectral sequence
$$E_{p,q}^2\ =\ \Tor_p(H^{-q}(K), B)\ \Longrightarrow\ H^{-p-q}(K \o* B).$%
$(See e.g.{\ }\Lcitemark 28\Rcitemark \ 11.34.)
Similarly, one has such a spectral sequence for $K_0$.  Moreover, $\phi$
induces a morphism from the spectral sequence for $K_0$ to the spectral
sequence for $K$.  Since $\phi$ is a quasi-isomorphism, the induced maps
\mapx[[ H^{-q}(K_0) || H^{-q}(K) ]] are isomorphisms, and so the induced maps
\mapx[[ \Tor_p(H^{-q}(K_0), B) || \Tor_p(H^{-q}(K), B) ]] are isomorphisms.
Hence $H^n(\phi)$ is an isomorphism.  Hence the \MV\ $A$-functor
$B \mapsto H^n(X_B, \shF_B)$ is isomorphic to the \MV\ $A$-functor
$B \mapsto H^n(K_0 \o* B)$.  \qed
\end{proofnodot}

The {\it upper semicontinuity\/} part of the proposition is of course the
usual theorem on upper semicontinuity of cohomology
(see e.g.{\ }\Lcitemark 22\Rcitemark \ III\ 12.8).

\begin{theorem}\label{mulch-material}
If $X$ is proper over a Dedekind domain $A$ and $\shF$ is a coherent sheaf on
$X$, then the \MV\ $A$-functor $H^n_\shF$ given by $B \mapsto H^n(X_B, \shF_B)$
is \MC.
\end{theorem}

Before proving this, there are some preliminaries.  By a {\it truncated \DVR},
we shall mean a ring $A$ of the form $R/I$ where $R$ is a \DVR\ and $I$ is a
proper nonzero ideal.  By a {\it uniformizing parameter\/} for $A$, we shall
mean the image in $A$ of a uniformizing parameter for $R$.

\begin{lemma}\label{mulch-material-lemma}
Let $A$ be a truncated \DVR\ with uniformizing parameter $t$.  Let $C$ be a
bounded complex of $A$-modules.  Assume that $H^n(C \o* A/(t^l))$ is \FG\ for
all $n$ and all $l$.  Then for each $n$, the $A$-functor $F$ given by
$B \mapsto H^n(C \o* B)$ is \MC.
\end{lemma}

\begin{proof}
Any module $M$ over $A$ is a direct sum of cyclic modules.  (See
e.g.{\ }\Lcitemark 11\Rcitemark \ Ch.\ VII\ \S2\ exercise 12(b).)  It
follows that if $M_0$ is a \FG\ submodule of $M$, then there exists a
\FG\ direct summand $M_1$ of $M$ with $M_0 \IN M_1$.

For each $l$, let $A_l$ denote $A/(t^l)$.
By working from low indices to high indices, one can construct a subcomplex
$C_0$ of $C$ with the properties that for each $n$:
\begin{alphalist}
\item for each $l$,
$\Ker[ d^n \o* A_l ] \IN C_0^n \o* A_l + \Im[ d^{n-1} \o* A_l ]$ and
\item $C_0^n$ is \FG\ and is a direct summand of $C^n$.
\end{alphalist}
Property (b) comes from the first paragraph.

Let $\oC = C/C_0$.  It follows from (b) that the sequence
\ses{C_0}{C}{\oC%
}is universally exact.  Since any $A$-module is a direct sum of cyclic modules,
it follows from (a) that the induced map
\mapx[[ H^n(C_0 \o* B) || H^n(C \o* B) ]] is surjective for every $n$ and
every $B$, and hence that we have short exact sequences:
\sesdot{H^{n-1}(\oC \o* B)}{H^n(C_0 \o* B)}{H^n(C \o* B)%
}Now construct a subcomplex $\oC_0 \IN \oC$ in the same way that
we constructed $C_0 \IN C$.  It follows that $F$ is expressible as the
cokernel of a map of \MC\ functors, and hence that $F$ is itself \MC.  \qed
\end{proof}

Now we prove \pref{mulch-material}, using arguments provided by Deligne.

\begin{proof}
There is an exact sequence of coherent sheaves on $X$:
\ses{\shF'}{\shF}{\shF''%
}in which $\shF'$ is $A$-torsion and $\shF''$ is $A$-torsion-free.  Since $A$
is a Dedekind domain, $\shF''$ is $A$-flat.  Hence this
sequence remains exact after tensoring over $A$ by anything, so by
the long exact sequence of cohomology and \pref{extension-is-mc}, it suffices
to prove the theorem when $\shF$ is either $A$-torsion or $A$-flat.  The
second case is taken care of by \pref{flat-implies-cohomologically-coherent}.
Therefore \WMAT\ $\shF$ is $A$-torsion.  Then in fact $a\shF = 0$ for some
nonzero $a \in A$.  We may then reduce to the case where $A$ is a truncated
\DVR.  Apply \pref{mulch-material-lemma} to the \u Cech complex of $\shF$.
\qed
\end{proof}

Finally we give a counterexample which shows that in general, $H^n_\shF$ is not
\MC.  The counterexample to be given here is based on
examples constructed by H.\ Cohen\Lspace \Lcitemark 13\Rcitemark \Rspace{}.
The same
examples appear in his thesis\Lspace \Lcitemark 14\Rcitemark \Rspace{}, where
Cohen remarks briefly
that the techniques therein yield a counterexample to Artin's problem.
However, he says nothing more about the matter.  It seems likely that Cohen
and/or Verdier did construct a counterexample, but it has now (apparently) been
lost.

Let $k$ be a field, let $A = k[[s,t]]$,
and let $X = \P3_A$, with coordinates $x,y,z,w$.  Let $\shF$ be the cokernel of
the map \mapx[[ \O_X || \O_X(1) ]] given by multiplication by $sx-ty$%
.\footnote{In terms of Cohen's construction, we have $r = 3$, $c_0 = s$,
$c_1 = -t$, and $m = 0$.}
Let $n = 1$.  We will show that $H$ is not \MC, using \pref{ding-dong}.

In his proposition 2, Cohen shows (in effect)
that if $B$ is a commutative $A$-algebra, and
$R = B[x^{\pm 1}, y^{\pm 1}, z^{\pm 1}, w^{\pm 1}]$,
then $H(B)$ is isomorphic to the degree $0$ part of the quotient of
$$\left\{ {f \over x^a y^b z^c w^d} \in R: f \in B[x,y,z,w],\ a,b,c,d \in \N,
  \hbox{\ and\ } (sx-ty)f = 0 \right\}$%
$by the sub-$B[x,y,z,w]$-module generated by
$$\left\{ {f \over yzw}, {f \over xzw}, {f \over xyw}, {f \over xyz} \in R:
  f \in B[x,y,z,w] \hbox{\ and\ } (sx-ty)f = 0 \right\}.$%
$
Let $B = A/(s^k, t^k)$, for some $k \in \N$.  We proceed to compute $H(B)$.

\begin{lemma}
The $B[x,y,z,w]$-module $\Ann_{B[x,y,z,w]}(sx-ty)$ is generated by
$$(st)^{k-1},\ (st)^{k-2} \sum_{i=0}^1 (sx)^i (ty)^{1-i},\ \ldots,
\ (st)^0 \sum_{i=0}^{k-1} (sx)^i (ty)^{k-1-i}.$$
\end{lemma}

\begin{sketch}
Let $C = k[s,t,x,y]/(s^k,t^k)$.  Since $B[x,y,z,w]$ is a flat $C$-algebra, it
suffices to show that the $C$-module $\Ann_C(sx-ty)$ admits the given
generators.  Let $f \in \Ann_C(sx-ty)$.  Write
$$f = \sum_{0 \leq i,j \leq k-1} f_{ij} s^i t^j,$%
$where $f_{ij} \in k[x,y]$.  The following assertions are easily checked:
\begin{itemize}
\item $f_{i,0} = 0$ if $0 \leq i \leq k-2$;
\item $f_{0,j} = 0$ if $0 \leq j \leq k-2$;
\item $f_{i,j-1} = (x/y) f_{i-1,j}$ if $1 \leq i,j \leq k-1$.
\end{itemize}
Hence $f$ is completely determined by $f_{0,k-1}, \ldots, f_{k-1,k-1}$.
The lemma follows.  \qed
\end{sketch}

{}From this lemma it follows that
$H(B)$ is isomorphic to the sub-$B$-module of $R$ generated by
$$\bigcup_{j=5}^k \left\{ {(st)^{k-j} \sum_{i=0}^{j-1} (sx)^i (ty)^{j-1-i}
  \over x^a y^b z^c w^d}: a,b,c,d \in \N \hbox{\ and\ }
  a+b+c+d = j-1 \right\}.$%
$Since this is a minimal generating set, we have:
$$\mu[H(B)]\ =\ \sum_{j=5}^k {j-2 \choose 3}\ \geq\ O(k^3).$%
$Since $(s,t)^{2k-1} \IN (s^k, t^k)$, it follows from \pref{ding-dong} that $H$
is not \MC.
\block{Global sections%
}\label{torsion-section}

Let us say that a group is {\it linear\/} if it may be embedded as a
subgroup of $\GL_n(k_1) \manytimes \GL_n(k_r)$ for some $n$ and some fields
$\vec k1r$.   We shall want to have some control over the fields: a group is
{\it $X$-linear\/} if there exist points $\vec x1r \in X$ (not necessarily
closed) and \FG\ field extensions $k_i$ of $k(x_i)$ (for each $i$) such that
the group may be embedded as a subgroup of
$\GL_n(k_1) \manytimes \GL_n(k_r)$ for some $n$.

The purpose of this section is two-fold.  The first purpose is to develop a
tool (the next theorem) for proving that groups are $X$-linear. The second
purpose is to study the torsion in $X$-linear groups.  Our result is
\pref{linear-implies-finite-torsion}, or in a slightly different form
\pref{linear-bound}.

\begin{theorem}\label{coherent-implies-linear}
Let $G$ be a group-valued locally coherent $X$-functor.  Then $G(X)$ is
$X$-linear.
\end{theorem}

Before proving this theorem, we will study the torsion in $X$-linear groups.
For any abelian group $H$, one can try to determine for which $n \in \N$
one has $\abs{{}_nH} < \infty$.  If $n = p_1^{k_1} \cdots p_r^{k_r}$, where
$\vec p1r$ are prime numbers, then $\abs{{}_nH} < \infty$ \IFF\
$\abs{{}_{p_i} H} < \infty$ for each $i$.  Therefore we may as well restrict to
the problem of determining when $\abs{{}_pH} < \infty$, where $p$ is prime.

If $C$ is a commutative ring, let us say that a morphism
\mp[[ \pi || X || \Spec(C) ]] is {\it essentially of finite type\/} if there
exists a commutative ring $D$, a homomorphism \mp[[ \phi || C || D ]]
which is essentially of finite type, and a morphism of finite type
\mp[[ \pi_0 || X || \Spec(D) ]], such that $\pi = \Spec(\phi) \circ \pi_0$.
Note that if $X$ is essentially of finite type over $\Z$, and $x \in X$, then
$k(x)$ is a \FG\ field extension of its prime subfield.

\begin{prop}\label{linear-implies-finite-torsion}
Let $H$ be an $X$-linear abelian group.
\begin{alphalist}
\item There are only finitely many prime numbers $p$ such that $H$ has infinite
      $p$-torsion.  Moreover, such a $p$ cannot be invertible in
      $\Gamma(X,\O_X)$.
\item If $X$ is essentially of finite type over $\Z$ or over $\Z_p$ (for some
      prime number $p$), then there exist prime numbers $\vec p1n$, none of
      which are invertible in $\Gamma(X,\O_X)$, such that the subgroup of $H$
      consisting of torsion prime to $p_1 \manycdot p_n$ is finite.
\end{alphalist}
\end{prop}

Clearly, if in the above proposition, $X$ is essentially of finite type over
$\Z_p$, then the list of primes $\vec p1n$ may be taken to be the single prime
$p$.  Also we have:

\begin{corollary}
Let $H$ be an $X$-linear abelian group.  If $X$ is essentially of finite
type over $\Q$ or over $\Q_p$ (for some prime number $p$), then the torsion
subgroup of $H$ is finite.
\end{corollary}

We now state a generalization of the proposition to the non-abelian case:

\begin{prop}\label{linear-bound}
Let $H$ be an $X$-linear group.
\begin{alphalist}
\item Let $n \in \N$ be invertible in $\Gamma(X,\O_X)$.  Then there exists some
      $N \in \N$, such that whenever $K$ is an $n$-torsion abelian subgroup of
      $H$, we have $\abs{K} \leq N$.
\item If $X$ is essentially of finite type over $\Z$ or over $\Z_p$ (for some
      prime number $p$), then there exist prime numbers $\vec p1n$, none of
      which are invertible in $\Gamma(X,\O_X)$, and some $N \in \N$, such that
      if $K$ is an abelian subgroup of $H$, and every element of $K$ is torsion
      prime to $p_1 \manycdot p_n$, then $\abs{K} \leq N$.
\end{alphalist}
\end{prop}

\begin{lemma}\label{injective-on-roots}
Let $(A,\lfM,k)$ be a local ring.  Let $n \in \N$, and assume that $n$ is
invertible in $A$.  Then the canonical map:
\dmapx[[ \setofh{\th{n} roots of unity in $A$}
      || \setofh{\th{n} roots of unity in $k$} ]]%
is injective.
\end{lemma}

\begin{proof}
Let $n \in \N$, $x \in A$, $a \in \lfM$, and suppose that $x^n = (x+a)^n = 1$.
Then $0 = [(x+a)^n - x^n] = a[nx^{n-1} + c]$, for some $c \in \lfM$.  But
$nx^{n-1}$ is a unit, so $nx^{n-1}+c$ is a unit, so $a = 0$.  \qed
\end{proof}

It is known (\Lcitemark 11\Rcitemark \ \S14, \#7, Cor.\ 2 to Prop.\ 17)
that a field \FG\ over its prime subfield (as a field extension) contains
only finitely many roots of unity.  This also holds for a field \FG\ over
$\Q_p$.  We will need a modest generalization of these statements:

\begin{lemma}\label{root-bound}
Let $K$ be a \FG\ field extension of $F$, where $F$ is $\Q$, or $\F_p$,
or $\Q_p$, for some prime number $p$.  Then there exists a constant $c$ such
that for every $m \in \N$, and every finite field extension $L$ of $K$ with
$[L:K] \leq m$, the number of roots of unity in $L$ is $\leq c m$
(if $F = \Q$), and is $\leq c^m$ if $F$ is $\F_p$ or $\Q_p$.
\end{lemma}

\begin{proof}
Let $\vec x1r$ be a transcendence basis for $K$ over $F$.  Let
$s = [K:F(\vec x1r)]$.  If $F = \Q$, let $c = 2s$.  If $F = \F_p$, let
$c = p^s$.  If $F = \Q_p$, let $c = 2sp^s$.

Let $L_0$ be the subfield of $L$ consisting of
elements algebraic over $F$.  Then $[L:F(\vec x1r)] = s[L:K] \leq sm$, so
$[L_0:F] \leq sm$.  If $F = \Q$ or $F = \F_p$, it is clear that
the given $c$ works.

Suppose that $F = \Q_p$.
Extend the standard absolute value on $\Q_p$ to $L_0$.  Let $(A,\lfM,k)$ be the
valuation ring of $L_0$.  If $x \in L_0$ is a root of unity, $\abs{x} = 1$, so
$x \in A$.  Since $k$ is an extension of $\F_p$ of degree $\leq sm$, the
number of elements in $k$ is bounded by $p^{sm}$.  By
\pref{injective-on-roots}, it follows that for any $r \in \N$ which is prime to
$p$ (and hence invertible in $A$), the number of \th{r} roots of unity in $A$
is $\leq p^{sm}$.

For $p \not= 2$, $\Q_p$ has no \th{p} roots of unity other than $1$ (see
\Lcitemark 24\Rcitemark \Rspace{}\ p.\ 20\ exercise 14).  For $p = 2$, $\Q_p$
contains no square root of $-1$.  For any field $M$, let $M'$ denote its
subfield generated by
\setof{x \in M: x^{(p^n)} = 1 \hbox{\ for some\ } n \in \N}.  Then
$\Q_p' = \Q$.  Hence $[L_0':\Q] \leq sm$.  Hence
$$\abs{\setof{x \in L_0: x^{(p^n)} = 1 \hbox{\ for some\ } n \in \N}}
       \ \leq\ 2sm.$%
$Hence the number of roots of unity in $L_0$ is $\leq (2sm)p^{sm} \leq c^m$.
\qed
\end{proof}

\begin{proofnodot}
(of \ref{linear-bound})
For part (a), \WMAT\ $H \IN \GL_r(k)$ for some field $k$, where $n$ is
invertible in $k$.  Let $\vec g1l \in \GL_r(k)$ be distinct commuting elements
with $g_i^n = 1$ for each $i$.  We need to prove that there is some $N \in \N$
(independent of $\vec g1l$) such that $l \leq N$.  Let $C$ be the subalgebra of
$\Mat_{r \times r}(k)$ generated by $\vec g1l$.  Then $C$ is a commutative,
artinian $k$-algebra, and $\Spec(C)$ has at most $r$ components, since
$\Mat_{r \times r}(k)$ has at most $r$ distinct nonzero orthogonal idempotents.
It follows from \pref{injective-on-roots} that the equation $x^n = 1$ has at
most $n^r$ solutions in $C$.  Let $N = n^r$.

For part (b), the field $k$ will be a \FG\ field extension of $\Q$, $\F_p$,
or $\Q_p$, for some prime number $p$.
Construct $C$ as in the preceding paragraph.  We have to bound
$l$ in terms of $r$ alone, and not in terms of $n$.  We may assume that $C$ is
local.  The residue field $L$ of $C$ is a finite extension of $k$, and
$[L:k] \leq r^2$.  Apply \pref{root-bound} and \pref{injective-on-roots}.  \qed
\end{proofnodot}

We now work towards a proof of \pref{coherent-implies-linear}.

\begin{lemma}\label{artinian-exists}
Let $M$ be a \FG\ $A$-module.  Then there exists a commutative artinian
$A$-algebra $B$, such that $B$ is essentially of finite type over $A$, and such
that the canonical map \mapx[[ M || M \o*_A B ]] is injective.
\end{lemma}

\begin{proof}
Let us say that an ideal $I \IN A$ is {\it good\/} if
there exists a commutative artinian $(A/I)$-algebra $B_{[I]}$ which is
essentially of finite type over $A/I$ such that the canonical map
\dmapx[[ M \o*_A (A/I) || M \o*_A B_{[I]} ]]%
is injective.

Let $I \IN A$ be an ideal, and suppose that every ideal properly containing
$I$ is good.  To prove the lemma, it suffices (by a sort of noetherian
induction) to show that $I$ is good.  Replacing $A$ by $A/I$, \WMAT\ $I = 0$.

Choose a primary decomposition $0 = Q_1 \manycap Q_r$ of $0$ in
$M$.  Then the map \mapx[[ M || M/Q_1 \manytimes M/Q_r ]] is injective and each
module $M/Q_i$ has a unique associated prime $\lfP_i$.

Since each module $M/Q_i$ admits a filtration with quotients isomorphic to
$A/\lfP_i$, it follows that there exists an integer $N$ such that the maps
\mapx[[ M/Q_i || (M/Q_i) \o*_A A/\lfP_i^N ]] are injective.  Note also that
the maps \mapx[[ M/Q_i || (M/Q_i) \o*_A A_{\lfP_i} ]] are injective.

We may assume that $\vec \lfP1r$ are arranged so that $\vec \lfP1k$ are
minimal primes of $A$ and $\vec \lfP{k+1}r$ are not.  Then
$\lfP_{k+1}^N,\ldots,\lfP_r^N$ are nonzero.  Let
\formulaqed{B = A_{\lfP_1} \times A_{\lfP_k} \times B_{[\lfP_{k+1}^N]}
                \times B_{[\lfP_r^N]}.}
\end{proof}

We use Witt rings in this and the next section.  The reader may find treatments
of the subject in
(\Lcitemark 26\Rcitemark \ VIII\ exercises\ 42--44),
(\Lcitemark 29\Rcitemark \ II\ \S5,\ \S6), and\Lspace \Lcitemark 9\Rcitemark
\Rspace{}.  Some
remarks about notation are in order.  There is a version of the Witt ring which
does {\it not\/} depend on the choice of a prime number $p$, as discussed for
example in (\Lcitemark 26\Rcitemark \ VIII\ exercise\ 42).  It
seems reasonable to denote this version of the Witt ring by $W(A)$.  There is a
second version of the Witt ring which does depend on the choice of a prime
number $p$, as discussed for example in
(\Lcitemark 26\Rcitemark \ VIII\ exercise\ 43).  To avoid
confusion, we will denote this version of the Witt ring by $W^p(A)$, as is done
in (\Lcitemark 9\Rcitemark \ p.\ 179).  However, the ring we denote by $W^p(A)$
is the same as the ring denoted $W(A)$ in\Lspace \Lcitemark 29\Rcitemark
\Rspace{}.  We will be
using the ring $W^p(A)$, as well as the truncated version $W^p_n(A)$.

\begin{lemma}\label{witt-coherent}
Let $k$ be a perfect field of positive characteristic $p$.  Fix $n \in \N$, and
let $A = W^p_n(k)$.  Let $F$ be a coherent $A$-functor.  Define a $k$-functor
$\tF$ by $\tF(B) = F(W^p_n(B))$, for all (commutative) $k$-algebras $B$.  Then
$\tF$ is coherent.
\end{lemma}

\begin{proof}
We let\ $\tilde{}$\ denote the operation which is in effect defined in the
statement.  Let $r$ be the level of $F$.

Suppose that $r = 0$, so \WMAT\ $F = \uM$ for some \FG\ $A$-module $M$.  Since
$A = W^p(k)/(p^n)$, and $W^p(k)$ is a complete \DVR\ with uniformizing
parameter $p$, it follows that $M$ may be expressed as a direct sum of modules
of the form $A/(p^i)$, where $i \in \setof{1,\ldots,n}$.  We may assume that in
fact $M = A/(p^i)$.  Then:
$$\tF(B)\ =\ A/(p^i) \o*_A W_n^p(B)\ =\ W_i^p(B),$%
$which may be identified (as a set) with $B^i$.  Hence $\tF$ is coherent.

Now suppose that $r \geq 1$.  By \pref{lesx-exists}, we may choose a
left exact sequence:
\lesx{F}{G}{\uN%
}of $A$-functors where $G$ is a coherent $A$-functor of level $r-1$
and $N$ is a \FG\ $A$-module.  We obtain a left exact sequence:
\lesx{\tF}{\tG}{\tilde{\uN}%
}of $k$-functors.  By induction on $r$, \WMAT\ $\tG$ is coherent,
and by the case $r = 0$, $\tilde{\uN}$ is coherent.  Hence $\tF$ is coherent.
\qed
\end{proof}

\begin{proofnodot}
(of \ref{coherent-implies-linear})
Since $G$ is locally coherent, it is a sheaf \WRT\ the Zariski topology, so we
may immediately reduce to the case where $X$ is affine, say $X = \Spec(A)$.  We
work with $A$-functors.

By \pref{lesx-exists}, there exists a \FG\ $A$-module $M$ and an
embedding $G\ \includeE{}\ \uM$ of $A$-functors.  (Note that this map need
not preserve the group structure; otherwise the proof would be much easier!)
By \pref{artinian-exists}, we can choose an artinian $A$-algebra $B$ which
is essentially of finite type over $A$ such that
the map \mapx[[ \uM(A) || \uM(B) ]] is injective.  Let $G|_B$ denote the
$B$-functor given by $G|_B(C) = G(C)$.  Then $G|_B$ is coherent by \pref{pull}.
Therefore it suffices to show that $G|_B$ has the desired property.  Replacing
$A$ by $B$, we may reduce to the case where $A$ is artinian.  Since $A$ is a
product of Artin local rings, we may in fact reduce to the case where $A$ is an
Artin local ring.

First suppose that $A$ is a field.  Then $A$ is a coherent
$A$-functor, so $G$ is representable by an affine group scheme of finite type
over $A$.  Hence (\Lcitemark 10\Rcitemark \ 11.11) $G(A)$ embeds in $\GL_n(A)$
for
some $n$.

Now suppose that $A$ contains a field.  Then from the Cohen structure theorem
for complete local rings, we know that $A$ contains a coefficient field $k$.
Since $A$ is artinian, it follows that $A$ is module-finite over $k$.  By
\pref{push}, we may reduce to the case $A = k$.

Finally, suppose that $A$ is an Artin local ring which does not contain a
field.  Then $A$ has mixed characteristic.  Let $\lfM$ be its maximal ideal,
and let $k$ be its residue field.  By (\Lcitemark 20\Rcitemark \ 0.6.8.3),
there exists a
(commutative) faithfully flat noetherian local $A$-algebra $(\tA,\tlfM, \ltK)$
with $\ltK$ being an algebraic closure of $k$, such that $\lfM\tA = \tlfM$.
{}From the latter fact, it follows that $\tA$ is artinian.  Since $G$ is
coherent, by \pref{sheaf} it is a sheaf for the fpqc topology, so
\mapx[[ F(A) || F(\tA) ]] is injective.  Hence \WMAT\ $A = \tA$
and so $k$ is algebraically closed.  By the Cohen structure theorem,
$A \iso W^p_n(k)[[\vec x1r]]/I$ for some $r,n \in \N$ and some ideal $I$.

Since $W^p_n(k)$ maps onto the residue field of $A$, and since
$A$ is artinian, it follows that $A$ is module-finite over $W^p_n(k)$.
By \pref{push}, we may reduce to the case $A = W^p_n(k)$.  Apply
\pref{witt-coherent} to reduce to the case $A = k$.  \qed
\end{proofnodot}
\block{Global sections -- arithmetic case}

In this section we refine the results of the last section, in the special case
where $X$ is of finite type over $\Z$.  In particular,
(\ref{linear-implies-finite-torsion}b) is supplanted by
\pref{arith-linear-structure}.

Let us say that a group is
{\it arithmetically linear\/} if it may be embedded as a subgroup of
$\GL_n(C)$, for some $n$ and some \FG\ commutative $\Z$-algebra $C$.
We will prove:

\begin{theorem}\label{coherent-implies-arithmetically-linear}
Assume that $X$ is of finite type over $\Z$.  Let $G$ be a group-valued locally
coherent $X$-functor.  Then $G(X)$ is arithmetically linear.  Moreover, if
$n$ is invertible in $\Gamma(X,\O_X)$, then the $n$-torsion in $G(X)$ is
finite.
\end{theorem}

First we analyze the structure of arithmetically linear abelian groups.
Recall that an abelian group $H$ is {\it bounded\/} if $nH = 0$ for some
$n \in \N$.  It is known [see\Lspace \Lcitemark 16\Rcitemark \Rspace{}\ 11.2 or
\Lcitemark 11\Rcitemark \Rspace{}\ Ch.\ VII\ \S2\ exercise 12(b)] that any
bounded
abelian group is a direct sum of cyclic groups.  Thus one may characterize the
bounded abelian groups as those which can be expressed as direct sums of cyclic
groups, in which the orders of the summands are bounded.

For purposes of this paper, let us say that an abelian group $H$ is
{\it cobounded\/} if it may be embedded as a subgroup of a direct sum of
(possibly infinitely many) copies of $\Z[1/n]$ for some $n \in \N$.  This is
equivalent to saying that $H$ is torsion-free and that $H \o*_\Z \Z[1/n]$ is a
free $\Z[1/n]$-module for some $n$.

\begin{remark}
We do not know of a structure theorem for abelian groups which are countable
and cobounded.  Certainly such groups can be rather complicated.  For
example, not every subgroup $H$ of $\o+_{k=1}^\infty \Z[1/2]$ can be expressed
as a direct sum of copies of $\Z$ and copies of $\Z[1/2]$; consider:
$$H\ = \ \Bigl\{ a \in \o+_{k=1}^\infty \Z[1/2]:
         \sum_{k=1}^\infty {a_k \over k} \in \Z \Bigr\}.$%
$Let $L$ be the maximal $2$-divisible subgroup of $H$.  Then $H/L \iso \Q$.
\end{remark}

Let us say that an abelian group $H$ is {\it\bxc\/} if $H \iso B \times C$ for
some bounded group $B$ and some cobounded group $C$.  We shall see
\pref{arith-linear-structure} that arithmetically linear abelian groups are
the same as countable (\bxc) abelian groups.

\begin{lemma}\label{bxc}
\
\begin{alphalist}
\item Let $M$ be a \bxc\ abelian group, and let $H$ be a subgroup of $M$.  Then
$H$ is \bxc.
\item Let
\ses{M'}{M}{M''%
}be a short exact sequence of abelian groups, in which $M'$ and $M''$ are \bxc.
Then $M$ is \bxc.
\end{alphalist}
\end{lemma}

\begin{proof}
Part (a) follows from (\Lcitemark 16\Rcitemark \ 50.3): if the torsion
subgroup of an abelian group is bounded, then it is a direct summand.

For part (b),
write $M' = B' \times C'$ and $M'' = B'' \times C''$ where $B'$, $B''$ are
bounded and $C'$, $C''$ are cobounded.  Let $M_{\op tor}$ be the torsion
subgroup of $M$.  We have a left exact sequence:
\les{B'}{M_{\op tor}}{B''%
}from which it follows that $M_{\op tor}$ is bounded, and hence that
$M_{\op tor}$ is a direct summand of $M$.  Therefore it suffices to show
that $M/M_{\op tor}$ is cobounded.  We have an exact sequence:
\ses{C'}{M/M_{\op tor}}{\ol{B''} \times C''%
}in which $\ol{B''}$ is a quotient of $B''$ and hence is bounded.  Choose
$n \in \N$ such that $n\ol{B''} = 0$, and such that $n$ satisfies the property
of $n$ in the definition of cobounded, for both $C'$ and $C''$.  Tensoring by
$\Z[1/n]$ yields an exact sequence:
\sesdot{C'[1/n]}{(M/M_{\op tor})[1/n]}{C''[1/n]%
}By construction, $C'[1/n]$ and $C''[1/n]$ are submodules of free
$\Z[1/n]$-modules and hence are themselves free.  Hence $(M/M_{\op tor})[1/n]$
is free so $M/M_{\op tor}$ is cobounded.  \qed
\end{proof}

\begin{lemma}\label{pork}
Let $A$ be a commutative $\Z$-algebra of finite type.  Let $M$ be a
\FG\ $A$-module.  Then the abelian group $M$ is \bxc.
\end{lemma}

\begin{proof}
As an abelian group, we may embed $M$ as a subgroup of the additive group of
the symmetric algebra of $M$, which is a \FG\ $A$-algebra.  Therefore we may
reduce to the case $M = A$.

Take a primary decomposition $0 = \lfQ_1 \manycap \lfQ_m$ of $0$ in $A$.
Then the canonical map \mapx[[ A || \prod_{i=1}^n A/\lfQ_i ]] is injective, so
we may reduce to the case where $A$ has a unique associated prime.

Let $q$ be the characteristic of $A$.  First suppose that $q > 0$.  Then
$qA = 0$ so $A$ is bounded.

Now suppose that $q = 0$.  By the Noether normalization lemma, there exist
elements $\vec x1r \in A \o*_\Z \Q$ which are algebraically independent over
$\Q$ and such that $A \o*_\Z \Q$ is module-finite over $\Q[\vec x1r]$.  We may
assume that $\vec x1r \in A$.  Let $S = \Spec(\Z)$, $X = \Spec(A)$,
$Y = \Spec(\Z[\vec x1r])$, so we have a morphism \mp[[ \phi || X || Y ]] of
$S$-schemes.  If $\eta \in S$ is the generic point, then $\phi_\eta$ is finite,
so it follows from (\Lcitemark 19\Rcitemark \ 8.1.2(a), 8.10.5(xii), 8.11.1,
9.6.1(vii)) that
for some $n \in \N$, $\phi \o*_S \Spec(\Z[1/n])$ is finite, i.e.\ that
$A_n = A \o*_\Z \Z[1/n]$ is a module-finite $\Z[1/n, \vec x1r]$-algebra.  Since
$A$ has characteristic zero and has a unique associated prime, it follows that
$A_n$ is a torsion-free $\Z[1/n, \vec x1r]$-module and that the canonical map
\mapx[[ A || A_n ]] is injective.  Hence $A$ embeds (as an abelian group) in
$(\Z[1/n, \vec x1r])^k$ for some $k$, so $A$ is cobounded.  \qed
\end{proof}

\begin{prop}\label{arith-linear-structure}
Let $H$ be an abelian group.  Then $H$ is arithmetically linear \IFF\ it is
countable and \bxc.
\end{prop}

\begin{proof}
First suppose that $H$ is countable and \bxc.  For some $n,k \in \N$ and some
prime numbers $\vec p1k$, we may embed $H$ as a subgroup of a countable direct
sum $K$ of copies of
$$\Z[1/n] \o+ \left( \o+_{i=1}^k \Z/p_i^n\Z \right).$%
$Then $K$ is the additive group of the ring
$$A = \Z[1/n,t] \times (\Z/p_1^n\Z)[t] \manytimes (\Z/p_k^n\Z)[t],$%
$so $K$ may be embedded as a subgroup of $\GL_2(A)$.  Hence $H$ is
arithmetically linear.

Now suppose that $H$ is arithmetically linear.  The countability of $H$ is
clear.  Embed $H$ as a subgroup of $\GL_r(A)$, for some \FG\ commutative
$\Z$-algebra $A$.  Let $R$ be the sub-$A$-algebra of $\Mat_{r \times r}(A)$
generated by $H$.  Then $R$ is a finite $A$-algebra, and $H$ is a subgroup of
$R^*$.  By (\ref{bxc}a), it suffices to show that $R^*$ is \bxc.  We may as
well view $R$ as an arbitrary \FG\ commutative $\Z$-algebra.

Let $J$ be the nilradical of $R$.  For each $n \in \N$, there is a short
exact sequence of abelian groups:
\diagramx{0&\mapE{}&J^n/J^{n+1}&\mapE{}&(R/J^{n+1})^*&\mapE{}&(R/J^n)^*&
          \mapE{}&1.%
}The group $(R/J)^*$ is \FG\ (see \ref{units-over-Z}), and hence is \bxc.
Hence
by (\ref{bxc}b), it suffices to show that $J^n/J^{n+1}$ is \bxc.  Apply
\pref{pork}.  \qed
\end{proof}

We now work towards a proof of
\pref{coherent-implies-arithmetically-linear}.
\def\oSinv{\overline{S}^{\kern2pt\lower3pt\hbox{$\scriptstyle -1$}}}

\begin{lemma}\label{localize-now}
Fix $n \in \N$ and a prime number $p$.  Let $S \IN W_n^p(A)$ be a
multiplicatively closed set.  Let \mp[[ \mu || W_n^p(A) || A ]] be the
canonical map.  Let $\oS = \mu(S)$.  Then the canonical map
\mp[[ i || W_n^p(A) || W_n^p(\oSinv A) ]] factors through $S^{-1}W_n^p(A)$.
\end{lemma}

\begin{proof}
Let \mp[[ \nu || W_n^p(\oSinv A) || \oSinv A ]] be the canonical map.  Let
$f \in S$.  Then the image of $\mu(f)$ in $\oSinv A$ is invertible, so
$\nu(i(f))$ is invertible.  But $\nu$ is surjective and has nilpotent kernel,
so $i(f)$ is invertible.  \qed
\end{proof}

Note that if $A$ is any $\F_p$-algebra, there is a canonical map
\mapx[[ W_n^p(\F_p) || W_n^p(A) ]], and since $W_n^p(\F_p) = \Z/p^n\Z$, we see
that $p^n = 0$ in $W_n^p(A)$.

\begin{prop}\label{witt-ring-of-poly-ring}
For any $n, k \in \N$, and any prime number $p$, the ring
$$W_n^p(\F_p[\vec t1k])$%
$is isomorphic to the subring of
$(\Z/p^n\Z)[t_1^{1/p^{n-1}},\ldots,t_k^{1/p^{n-1}}]$ generated by the elements
$p^r t_i^{j/p^r}$ for $0 \leq r \leq n-1$, $1 \leq i \leq k$, and
$1 \leq j \leq p-1$.
\end{prop}

\begin{proof}
Let $R = \F_p[t_1^{p^{-\infty}},\ldots,t_k^{p^{-\infty}}]$, and let
$V = W_n^p(R)$.  Let
$$W = (\Z/p^n\Z)[t_1^{p^{-\infty}},\ldots,t_k^{p^{-\infty}}].$%
$It is easily seen that there exists a unique ring homomorphism
\mp[[ \eta || W || V ]] with the property that
$\eta(t_i^{p^{-m}}) = (t_i^{p^{-m}},0,\ldots,0)$ for each $i$ and each
$m \geq 0$.

We show that $\eta$ is surjective.  Let \mp[[ \mu || V || R ]] be the
canonical map.  Since $R$ is perfect, $\Ker(\mu) = (p)$.  Because of this,
because $p^n = 0$ in $V$, and because $t_i^{1/p^r}$ (for various $i$, $r$)
generate $R$ as a $\Z$-algebra, it follows that the $\eta(t_i^{1/p^r})$
generate $V$ as a $\Z$-algebra.  Hence $\eta$ is surjective.

We show that $\eta$ is injective.  Suppose otherwise.  Let $x \in \Ker(\eta)$,
$x \not= 0$.  We may assume that $px = 0$.  Then $x = p^{n-1}y$ for some
$y \in W$; \WMAT\ the coefficients which appear in $y$ lie in the set
\setof{1,\ldots,p-1}.  It follows that the $\th{0}$ component of $\eta(y)$ is
nonzero.  Hence $p^{n-1} \eta(y) \not= 0$, so $\eta(x) \not= 0$: contradiction.
Hence $\eta$ is injective.

We return to the proof of the proposition.  Let
$A_n = W_n^p(\F_p[\vec t1k])$.  From what we have just done, it follows that
$A_n$
may be identified with a subring of $W$.  Since
$\eta(p^r t_i^{j/p^r}) = (0,\ldots,0,t_i^j,0,\ldots,0)$, where $t_i^j$ appears
in the \th{r} spot, it follows that $p^r t_i^{j/p^r} \in A_n$, for each $i$,
$r$, and $j$.  To complete the proof, we must show that these elements generate
$A_n$.  Let $A_n'$ be the subring of $A_n$ generated by the elements
$p^r t_i^{j/p^r}$.  Consider the canonical map \mp[[ \tau || A_n || A_{n-1} ]].
By induction on $n$ \WMAT\ the elements $p^r t_i^{j/p^r}$ generate $A_{n-1}$,
and so that $\tau(A_n') = A_{n-1}$.  Let $f \in \F_p[\vec t1k]$.  Write
$f = \sum_I a_I t^I$, where $I$ is a multi-index.  Since
$a_I \in \setof{0,\ldots,p-1}$, we may view $a_I$ as an element of
$\Z/p^n\Z$.  Let $\ltF = p^{n-1} \sum_I a_I (t^I)^{1/p^{n-1}} \in W$.  Then
$\ltF$ corresponds to the element $(0,\ldots,0,f)$ of $A_n$.  Hence
$\Ker(\tau) \IN A_n'$.  Since $\tau(A_n') = A_{n-1}$, it follows that
$A_n' = A_n$.  \qed
\end{proof}

\begin{corollary}\label{witt-finite}
Fix $n \in \N$ and a prime number $p$.  Let $A$ be an $\F_p$-algebra of finite
type.  Then $W_n^p(A)$ is a $\Z$-algebra of finite type.
\end{corollary}

\begin{proof}
Choose a surjection \mp[[ \pi || \F_p[\vec x1k] || A ]].  Then $W_n^p(\pi)$ is
a surjection, and since \pref{witt-ring-of-poly-ring} $W_n^p(\F_p[\vec x1k])$
is of finite type over $\Z$, so is $W_n^p(A)$.  \qed
\end{proof}

The following result is a variant of \pref{witt-coherent}.

\begin{lemma}\label{sneezewort}
Fix $n \in \N$ and a prime number $p$.  Assume that $A$ is an $\F_p$-algebra of
finite type.  Let $C = W_n^p(A)$.  (By \pref{witt-finite} $C$ is noetherian.)
Let $F$ be a coherent $C$-functor, built up from
$\setof{C/(p^k)}_{1 \leq k \leq n}$.  Let $G$ be the $A$-functor given by
$G(B) = F(W_n^p(B))$.  Let $\tG$ be the sheaf associated to $G$ for the ffqc
topology.  Then $\tG$ is representable.
\end{lemma}

\begin{remarks}
The functor $G$ is not in general coherent, as it is not in general an ffqc
sheaf.  For an arbitrary coherent $C$-functor $F$, with no restrictions on how
it is built up, it may be that the corresponding functor $\tG$ is always
representable.
\end{remarks}

\begin{proofnodot}
(of \ref{sneezewort}.)
Since any limit of finitely many representable functors is representable,
\WMAT\ $F = \ul{C/(p^k)}$ for some $k$.  Then we may describe $G$ by:
$$G(B)\ =\ {\setof{(\vec b0{n-1}): b_i \in B \hbox{\ for all\ } i} \over
            \setof{(0,\ldots,0,d_k^{p^k},\ldots,d_{n-1}^{p^k}):
            d_i \in B \hbox{\ for all\ } i}},$%
$where it is to be understood that this quotient of abelian groups takes place
\WRT\ the abelian group structure on $W_n^p(B)$.  For for any
$b \in B$, there exists a faithfully flat ring extension
\mp[[ \phi || B || B' ]] such that $\phi(b)$ is a \th{(p^k)} power.  Therefore
to complete the proof, it suffices to show that the $A$-functor $H$ given by
$$H(B)\ =\ {\setof{(\vec b0{n-1}): b_i \in B \hbox{\ for all\ } i} \over
            \setof{(0,\ldots,0,e_k,\ldots,e_{n-1}):
            e_i \in B \hbox{\ for all\ } i}}$%
$is representable.  (Then we will have $H = \tG$.)  But
\begin{eqnarray*}
H(B) & \iso & W_n^p(B) / \Ker[W_n^p(B)\ \mapE{}\ W_k^p(B)]\\
     & \iso & W_k^p(B)\ \iso\ B^k,
\end{eqnarray*}
so $H$ is representable.  \qed
\end{proofnodot}

Let $A = \F_p[\vec t1k]$.  Let $C = (\Z/p^n\Z)[\vec t1k]$.  There is a
canonical map \mp[[ \phi || C || W_n^p(A) ]] given by setting
$\phi(t_i) = (t_i,0,\ldots,0)$ for each $i$.

\begin{prop}\label{bindweed}
Let $A = \F_p[\vec t1k]$.  Let $C = (\Z/p^n\Z)[\vec t1k]$.
Let $S$ be a multiplicatively closed subset of $C$.  Let \mp[[ \pi || C || A ]]
be the canonical map.  Let $F$ be a coherent $S^{-1}C$-functor.  Assume that
$F$ is built up from $\setof{S^{-1}C/(p^r)}_{1 \leq r \leq n}$.  Let
$\oS = \pi(S)$.  Let $G$ be the ffqc sheaf associated to the
$\oS^{-1}A$-functor given by $B\mapsto F(W_n^p(B))$.  (This makes sense by
\ref{localize-now}.) Then:
\begin{alphalist}
\item $G$ is representable;
\item the canonical map \mp[[ i || F(S^{-1}C) || G(\oS^{-1}A) ]] is injective.
\end{alphalist}
\end{prop}

\begin{proof}
{\bf (a):} Let $D = W_n^p(\oS^{-1}A)$.
Let \mp[[ \phi || \Spec(D) || \Spec(S^{-1}C) ]] be
the canonical map.  Then $\phi^*F$ is a $D$-functor, which is coherent (see
\ref{pull}), and is in fact built up from $\setof{D/(p^r)}_{1 \leq r \leq n}$.
Apply \pref{sneezewort}.

{\bf (b):} The construction is functorial in $F$, so we may reduce
to the case where $F = \ul{S^{-1}C/(p^r)}$.  In that case, one sees that $i$ is
isomorphic to the canonical map
\mapx[[ S^{-1}(\Z/p^r\Z)[\vec t1k] || W_r^p(\oS^{-1}\F_p[\vec t1k]) ]].  We may
reduce to showing that the canonical map
\mp[[ j || (\Z/p^r\Z)[\vec t1k] || W_r^p(\F_p[\vec t1k]) ]] is injective.
This follows from the proof of \pref{witt-ring-of-poly-ring}.  \qed
\end{proof}

\begin{lemma}\label{generic-representability}
Assume that $A$ is a domain.  Let $F$ be a coherent $A$-functor.  Then
for some $f \in A - \setof{0}$, the pullback [along
\mapx[[ \Spec(A_f) || \Spec(A) ]]{}] of $F$ to $A_f$ is representable
by an $A$-algebra of finite type.
\end{lemma}

\begin{sketch}
The functor $F$ is built up from finitely many $A$-modules, each \FG.  Pick
some $f \in A - \setof{0}$ such that the localization of each such module at
$f$ is free.  \qed
\end{sketch}

\begin{lemma}\label{generic-linearity}
Assume that $X$ is integral.  Let $G$ be an affine group scheme of finite type
over $X$.  Then there exists some $n \in \N$, a nonempty open subscheme
$U \IN X$, and a closed immersion \mapx[[ G_U || \GL_n(U) ]] of $U$-schemes
which is also a homomorphism.
\end{lemma}

\begin{sketch}
Let $\eta$ be the generic point of $X$.  For some $n \in \N$, we have a closed
immersion (and a homomorphism) \mp[[ h || G_\eta || \GL_n(X)_\eta ]] of
$\Spec k(\eta)$-schemes.  There is a nonempty open subscheme $V \IN X$ and a
closed immersion \mapx[[ G_V || \GL_n(V) ]] of $V$-schemes which induces $h$.
Replacing $V$ by a sufficiently small nonempty open subscheme $U$, we obtain
(by restriction) a homomorphism (and a closed immersion)
\mapx[[ G_U || \GL_n(U) ]] of $U$-schemes.  \qed
\end{sketch}

\begin{lemma}\label{soup-bone}
Let $A = (\Z/p^n\Z)[\vec t1k]$.  Let $M$ be a \FG\ $A$-module.  Then there
exists a non-zero-divisor $f \in A$ and positive integers $\vec l1r$ such
that $M_f \iso A_f/(p^{l_1}) \manyo+ A_f/(p^{l_r})$ as $A_f$-modules.
\end{lemma}

\begin{sketch}
Let $S$ be the set of non-zero-divisors of $A$.  In $S^{-1}A$, the only
elements (up to associates) are $1,p,\ldots,p^{n-1},0$.  Hence every ideal of
$S^{-1}A$ is principal, so every \FG\ $S^{-1}A$-module is a direct sum of
cyclic modules, necessarily of the form $S^{-1}A/(p^j)$ for various $j$.  The
lemma follows.  \qed
\end{sketch}

\begin{proofnodot}
(of \ref{coherent-implies-arithmetically-linear}.)
The comment about what happens when $n$ is invertible in $\Gamma(X,\O_X)$
follows from (\ref{linear-implies-finite-torsion}b), so a direct proof is
omitted.

Some of the steps here follow the proof of \pref{coherent-implies-linear}.
We may assume that $X$ is affine, $X = \Spec(A)$.  Choose $M$ and $B$ as in the
proof of \pref{coherent-implies-linear}.  Write $B = S^{-1}C$ for some
\FG\ $A$-algebra $C$ and some multiplicatively closed subset $S \IN C$.
Certainly we may replace $A$ by $C$, so $B = S^{-1}A$.  Replacing $A$ by $A_g$
for some suitably chosen $g \in S$, \WMAT\ the connected components of
$\Spec(A)$ are irreducible and that they correspond bijectively with the points
of $\Spec(S^{-1}A)$.  Write $A = A_1 \manytimes A_m$, where $\Spec(A_i)$ is
irreducible for each $i$.  Localizing further if necessary, \WMAT\ $A_i$ has a
unique associated prime for each $i$.

We may reduce to the following situation: $A$ has a unique associated prime,
$S^{-1}A$ is an Artin local ring.  It suffices to show that for some $f \in S$,
$G(A_f)$ is arithmetically linear.  Let $r$ be the characteristic of $A$.

Since $A$ has a unique associated prime and $S^{-1}A$ is an Artin local ring,
every non-nilpotent element of $A$ lies in $S$.

First suppose that $r = 0$.  By Noether normalization, we may find
algebraically
independent elements $\vec x1s \in A$ such that $A \o*_\Z \Q$ is module-finite
over $\Q[\vec x1s]$.  As in the proof of \pref{pork}, there is some $n \in \N$
such that $A \o*_\Z \Z[1/n]$ is a module-finite $\Z[1/n, \vec x1s]$-algebra.
Every nonzero element of $\Z[\vec x1s]$ lies in $S$.  We may replace $A$ by
$A \o*_\Z \Z[1/n]$.  By (\ref{push}b), the pushforward of $G$ to
$\Z[1/n, \vec x1s]$ is coherent.  By \pref{generic-representability}, there
is some $f \in \Z[1/n,\vec x1s] - \setof{0}$ such that the pullback of
$G$ to $\Z[1/n,\vec x1s,f^{-1}]$ is representable by an affine
$\Z[1/n,\vec x1s,f^{-1}]$-algebra of finite type.  Apply
\pref{generic-linearity}.

Now suppose that $r > 0$.  Then $r = p^m$ for some prime $p$ and some $m$.  By
the Noether normalization theorem, we may find algebraically independent
elements $\vec x1s \in A$ such that $A \o*_\Z \F_p$ is module-finite over
$\F_p[\vec x1s]$.  It follows that $A$ is module-finite over
$E = (\Z/p^m\Z)[\vec x1s]$.  By (\ref{push}b), the pushforward $H$ of $G$ to
$E$ is coherent, so we may reduce to the case where $A = (\Z/p^m\Z)[\vec x1s]$.
Let $\vec M1d$ be $A$-modules from which $G$ can be built up.  By
\pref{soup-bone}, there is some $f \in S$ such that each $(M_i)_f$ is a direct
sum of modules of the form $A_f/(p^r)$, for various $r$.  Let $G_f$ be the
pullback of $G$ to $A_f$.  Then $G_f$ is built up from
$\setof{A_f/(p^r)}_{1 \leq r \leq n}$.  The theorem follows now by applying
\pref{bindweed} and \pref{generic-linearity}.  \qed
\end{proofnodot}

\begin{problemx}
Which arithmetically linear groups arise as $G(X)$, for some scheme $X$ of
finite type over $\Z$ and some group-valued coherent $X$-functor $G$?
\end{problemx}

Certain groups can be shown to be quotients of arithmetically linear abelian
groups by \FG\ subgroups.  For example, in the next section we shall see that
this is the case for $\Pic(X)$, where $X$ is a reduced scheme of finite type
over $\Z$.  Although it may in fact be the case that $\Pic(X)$ is itself
arithmetically linear, we have not been able to show this, so we are lead to
the following lemma:

\begin{lemma}\label{arith-linear-over-finite}
Let $G$ be an arithmetically linear abelian group.  Let $H$ be a
\FG\ subgroup of $G$.  Then the torsion subgroup of $G/H$ is supported at a
finite set of primes.
\end{lemma}

\begin{proof}
Write $G = B \times C$, where $B$ is bounded and $C$ is cobounded.  For any
abelian group $M$, let $M[1/n]$ denote $M \o*_\Z \Z[1/n]$.  Choose $n \in \N$
such that $nB = 0$ and such that $C[1/n]$ is a free $\Z[1/n]$-module.  Let
$Q = G/H$.  It suffices to show that the torsion subgroup of $Q[1/n]$ is
supported at a finite set of primes.  We have an exact sequence:
\ses{H[1/n]}{C[1/n]}{Q[1/n]%
}of $\Z[1/n]$-modules.  Since $H[1/n]$ is contained in a \FG\ direct summand of
$C[1/n]$, it suffices to show that for any \FG\ $\Z[1/n]$-module $M$, the
torsion subgroup of $M$ is supported at a finite set of primes.  This is easily
checked.  \qed
\end{proof}
\block{Application to the Picard group}

Let $M$ be a \FG\ $A$-module.  Let $\uPGL(M)$ denote the Zariski sheaf
associated to the $A$-functor given by
$$B \mapsto {\Aut_B(M \o*_A B) \over B^*}.$%
$More generally, let $\shM$ be a coherent $\O_X$-module.  Let $\uPGL(\shM)$
denote the Zariski sheaf associated to the $X$-functor given by
$$Y\ \mapsto\ {\Aut_Y(\shM_Y) \over \Gamma(\O_Y^*)}.$%
$Since $\uAut(\shM)$ acts by conjugation on $\uEnd(\shM)$, we obtain
a canonical morphism of group-valued $X$-functors:
\dmap[[ \psi || \uPGL(\shM) || \AUTalg(\uEnd(\shM)).\footnote{For the
actual arguments which we use, one could substitute the simpler functor
$\AUT(\uEnd(\shM))$, but we use $\AUTalg(\uEnd(\shM))$ instead for asthetic
reasons because it makes $\psi$ closer to an isomorphism.} ]]%
The $X$-functor $\AUTalg(\uEnd(\shM))$ is locally coherent by
\pref{autalg-is-coherent} and example \pref{end-example} from
\S\ref{examples-section}; it is coherent, at least assuming that $X$ is
separated, by \pref{lump}.

We would like to show that $\uPGL(\shM)$ is coherent.  We could do this by
showing that $\psi$ is an isomorphism.  However, other than the case where
$\shM$ is locally free, we do not know if $\psi$ is an isomorphism.  Therefore,
we settle for showing that (under certain special circumstances) $\psi(X)$ is
injective.  This is a weak substitute for showing that $\uPGL(\shM)$ is
coherent.\footnote{In trying to show that $\psi$ is a monomorphism, one comes
to the following question: Let $M$ be a \FG\ $A$-module.  Let
$\sigma$ be an automorphism of $M$ as an $A$-module.  Assume that for every
commutative $A$-algebra $B$, $\sigma \o*_A B$ lies in the center of
$\End_B(M \o*_A B)$.  Does it follow that $\sigma$ is a homothety, i.e.\ that
$\sigma$ is given by multiplication by an element of $A$?}

First we prove a lemma, then a corollary which
says something directly about $\psi$.

\begin{lemma}\label{homothety}
Assume that $A$ is reduced.  Let $A'$ be a ring, with
$A \IN A' \IN \nor{A}$.  Assume that the ideal $[A:A']$ of $A$ is prime.  Let
$\sigma$ be an endomorphism of $A'$ as an $A$-module.  Assume that for every
$\lfP \in \Spec(A)$, $\sigma \o*_A k(\lfP)$ is a homothety of
$A' \o*_A k(\lfP)$ as a $k(\lfP)$-module.  Then $\sigma$ is a homothety of $A'$
as an $A$-module.
\end{lemma}

\begin{proof}
First we show that $\sigma$ is given by multiplication by $\sigma(1)$.
By subtracting the endomorphism of $A'$ given by multiplication by $\sigma(1)$,
we may reduce to showing that if $\sigma(1) = 0$, then $\sigma = 0$.  Let
$\vec\lfP1r$ be the minimal primes of $A$.  Using the fact that
$\sigma \o*_A k(\lfP_i)$ is a homothety, we conclude that
$\sigma \o*_A k(\lfP_i) = 0$.  Let $x \in A'$.  Then $\sigma(x) \mapsto 0$
in $A' \o*_A k(\lfP_i)$.  Since the map
\dmapx[[ A' || \oplus_{i=1}^r A' \o*_A k(\lfP_i) ]]%
is injective, it follows that $\sigma(x) = 0$, and hence that $\sigma = 0$.
Hence (reverting to the original problem), we see that $\sigma$ is given
by multiplication by $\sigma(1)$.

Let $Q = A'/A$, which is an $A$-module.  Let $\lfP = [A:A']$.
Let $\ol{\sigma(1)}$ denote the image of $\sigma(1)$ in $Q$.  It follows
that $\ol{\sigma(1)} \mapsto 0$ in $Q \o*_A k(\lfP)$.  But by the construction
of $\lfP$, the canonical map \mapx[[ Q || Q \o*_A k(\lfP) ]] is injective.
Hence $\ol{\sigma(1)} = 0$.  Hence $\sigma(1) \in A$.  \qed
\end{proof}

\begin{corollary}\label{key-embedding}
Assume that $A$ is reduced.  Let $A'$ be a ring, with
$A \IN A' \IN \nor{A}$.  Assume that the ideal $[A:A']$ of $A$ is prime.
Assume that $A'$ is a \FG\ $A$-module.  Let
\dmap[[ \psi || \uPGL(A') || \AUTalg(\uEnd(A')) ]]%
be the canonical (conjugation) morphism of group-valued $A$-functors.  Then
$\psi(A)$ is injective.
\end{corollary}

\begin{proof}
Consider the morphism \mp[[ \psi_0 || \uAut(A') || \AUTalg(\uEnd(A')) ]]
which induces $\psi$.  Then $\ker(\psi_0(A))$ consists of those automorphisms
$\sigma$ of $A'$ as an $A$-module with the property that
$\sigma \o*_A B \in Z[\End_B(A' \o*_A B)]$ for every commutative $A$-algebra
$B$.  If $B$ is a field, it follows that for such a $\sigma$,
$\sigma \o*_A B$ is given by multiplication by an element of $B$.
By \pref{homothety}, such a $\sigma$ is itself a homothety.  Hence
$\ker(\psi_0(A)) = A^*$.

Let $f \in A - \setof{0}$.  Then the ideal $[A_f:A'_f]$ of $A_f$ is prime
and it follows that $\ker(\psi_0(A_f)) = A_f^*$.  Hence the map
\dmapx[[ {\Aut_{A_f}(A' \o*_A A_f) \over A_f^*} || \AUTalg(\uEnd(A'))(A_f) ]]%
is injective.  Considering the sheafification which occurs in the definition
of $\uPGL(A')$, we see that $\psi(A)$ is injective.  \qed
\end{proof}

\begin{lemma}\label{filtration-of-normalization}
Let $C$ be an overring of $A$, with $C$ \FG\ as an $A$-module.
Then there exists a chain:
$$A\ =\ A_0\ \IN\ A_1\ \IN\ \cdots\ \IN\ A_n\ =\ C$%
$of rings such that for each $k = 1, \ldots, n$, the ideal $[A_{k-1}:A_k]$
of $A_{k-1}$ is prime.
\end{lemma}

\begin{proof}
We may assume that $C \not= A$.  For each $y \in C - A$, let
$I_y = [A:A[y]]$.  (All conductors are to be computed as ideals in $A$.)
Choose $y \in C - A$ so that $I_y$ is maximal amongst all such ideals.
We will show that $I_y$ is prime.  This will complete the proof.

Pick $a \in A$ such that $[I_y:a]$ is prime.  Choose $n$ so that $y$ satisfies
a monic polynomial of degree $n$ with coefficients in $A$.  Then
\snap\hbox{$I_y = [A: y,y^2,\ldots,y^{n-1}]$}, and
$[I_y:a] = [A: ay,ay^2,\ldots,ay^{n-1}]$.

First we show that $ay \notin A$.  Suppose otherwise.  Since
$[I_y:a] \not= A$, $a \notin I_y$.  Hence $ay^k \notin A$, for some $k$
with $2 \leq k \leq n-1$.  Consider all pairs $(r,s) \in \N^2$,
$s < n$, such that $a^r y^s \notin A$.  Choose such a pair $(r_0,s_0)$ such
that the ratio $s_0/r_0$ is as small as possible.  (This is possible because
$ay \in A$, and hence $a^r y^s \in A$ for all $r \geq s$.)

Let $(r,s) \in \N^2$ be such that $a^r y^s \notin A$.  We will show that
$s_0/r_0 \leq s/r$.  Suppose otherwise: $s_0/r_0 > s/r$.  We may assume that
$s \geq n$.  We have:
$$a^r y^s = \sum_{i = 0}^{n-1} c_i a^r y^i$%
$for suitable $c_i \in A$.  For $i$ in the given range, $s/r > i/r$, so
$s_0/r_0 > i/r$.  Hence $a^r y^i \in A$.  Hence $a^r y^s \in A$: contradiction.
Hence $s_0/r_0 \leq s/r$.

Let $z = a^{r_0} y^{s_0}$.
We have $I_y \IN I_z$, so by the maximality of $I_y$, we have $I_y = I_z$.
For any $k \in \N$, $az^k = a^{r_0 k+1} y^{s_0 k}$, and
$${s_0 k \over r_0 k + 1} < {s_0 \over r_0}.$%
$Hence $az^k \in A$.  Hence $a \in I_z$.  Hence $a \in I_y$: contradiction.
Hence $ay \notin A$.

Note that $ay$ satisfies a monic
polynomial of degree $n$ with coefficients in $A$.  We have:
$$I_y\ \IN\ [I_y:a]\ \IN\ [A:ay,(ay)^2,\ldots,(ay)^{n-1}]\ =\ I_{ay},$%
$so by the maximality of $I_y$ we must have $I_y = I_{ay}$ and hence $I_y$ is
prime.  \qed
\end{proof}

Let us say that an abelian group $G$ is {\it pseudo-$X$-linear\/} if there
exists a filtration
$$0 = G_0 \IN G_1 \manyIN G_m = G,$%
$with the property that $G_i/G_{i-1}$ is $X$-linear for each $i$.  It is
conceivable that every such group $G$ is $X$-linear.

\begin{theorem}\label{filtration}
Assume that $X$ is reduced, and that the canonical map
\mp[[ \pi || \nor{X} || X ]] is finite.  Let $C$ be the quotient sheaf
$\pi_*\O_{\nor{X}}^*/\O_X^*$.  Then the group $C(X)$ is pseudo-$X$-linear.  If
$X$ is of finite type over $\Z$, then $C(X)$ is arithmetically linear.
\end{theorem}

\begin{proof}
The comment about what happens when $X$ is of finite type over $\Z$ is left
to the reader; the proof given below works with appropriate changes, provided
that one uses in addition (\ref{bxc}b) and \pref{arith-linear-structure}.  One
uses \pref{coherent-implies-arithmetically-linear} instead of
\pref{coherent-implies-linear}.

The following two facts are easily verified: any subgroup of a
pseudo-$X$-linear abelian group is pseudo-$X$-linear, and any product of
finitely many pseudo-$X$-linear abelian groups is pseudo-$X$-linear.  It
follows that we may reduce to the case where
$X$ is affine.  Then by \pref{filtration-of-normalization}, we may reduce to
the following situation: $A$ and $A'$ are reduced noetherian rings, with
$A \IN A' \IN \nor{A}$ (and $\nor{A}$ is module-finite over $A$), and the
ideal $[A:A']$ of $A$ is prime.  We must show that if $F$ is the Zariski sheaf
associated to the $A$-functor given by $B \mapsto (A' \o*_A B)^*/B^*$, then
$F(A)$ is $X$-linear.  But $F(A) \iso [\uPGL(A')](A)$, so by
\pref{key-embedding}, it suffices to show that if $G = \AUTalg(\uEnd(A'))$,
then $G(A)$ is $X$-linear.  But $G$ is coherent by \pref{autalg-is-coherent},
so the theorem follows from \pref{coherent-implies-linear}.  \qed
\end{proof}

One can check that \pref{linear-implies-finite-torsion} applies to a
pseudo-$X$-linear abelian group, so one obtains the following corollary:

\begin{corollary}\label{finiteness-of-C}
Assume that $X$ is reduced, and that the canonical map
\mp[[ \pi || \nor{X} || X ]] is finite.  Let $C$ be the quotient sheaf
$\pi_*\O_{\nor{X}}^*/\O_X^*$.  Then:
\begin{alphalist}
\item There are only finitely many prime numbers $p$ such that ${}_p C(X)$ is
      infinite.  Moreover, such a $p$ cannot be invertible in $\Gamma(X,\O_X)$.
\item If $X$ is essentially of finite type over $\Z$ or $\Z_p$ (for some prime
      number $p$), then there exist prime numbers $\vec p1n$, none of which are
      invertible in $\Gamma(X,\O_X)$, such that the subgroup of $C(X)$
      consisting of torsion prime to $p_1 \manycdot p_n$ is finite.
\end{alphalist}
\end{corollary}

\begin{corollary}\label{pic-main-theorem}
Assume that $X$ is reduced, and that the canonical map
\mp[[ \pi || \nor{X} || X ]] is finite.  Let
$Q = \Gamma(\O_{\nor{X}}^*)/\Gamma(\O_X^*)$.  Let $K$ be the kernel of the
canonical map \mapx[[ \Pic(X) || \Pic(\nor{X}) ]].
\begin{alphalist}
\item Fix $n \in \N$, and assume that $n$ is invertible in $\Gamma(X,\O_X)$.
      If $Q$ is \FG, or more generally if it admits an $n$-divisible subgroup
      with \FG\ cokernel, then ${}_n K$ is finite.
\item If $Q$ admits a divisible subgroup with \FG\ cokernel, then there are
      only finitely many prime numbers $p$ such that ${}_p K$ is infinite.
\end{alphalist}
\end{corollary}

\begin{proof}
Any \FG\ projective module of constant rank over a semilocal ring is free, so
any line bundle on $\nor{X}$ may be trivialized by an open cover pulled back
from $X$.  (This argument was shown to me by R.\ Wiegand.)  Hence
$R^1\pi_*(\O_{\nor{X}}^*) = 0$, and so by the Leray spectral sequence we see
that $H^1(X,\pi_* \O_{\nor{X}}^*) \iso \Pic(\nor{X})$.

Run the long exact sequence of cohomology on
\sesdot{\O_X^*}{\pi_*\O_{\nor{X}}^*}{C%
}We obtain a short exact sequence:
\sesdot{Q}{C(X)}{K%
}By \pref{finiteness-of-C}, it suffices to show that if
\ses{M'}{M}{M''%
}is a short exact sequence of abelian groups, $\abs{{}_n M} < \infty$, and
there exists an $n$-divisible subgroup $H \IN M'$ such that $M'/H$ is \FG,
then $\abs{{}_n M''} < \infty$.  This is easily proved -- see
\Lcitemark 21\Rcitemark \Rspace{}.  \qed
\end{proof}

\begin{example}
Let $X = \Spec \Q[x,y]/(x^2 - 2y^2)$.  Then
$\Gamma(\O_{\nor{X}}^*) / \Gamma(\O_X^*) \iso \Q[\sqrt{2}]^*/\Q^*$, which is
not \FG.
\end{example}

We consider what happens when $X$ is of finite type over $\Z$.  We need the
following well-known result, which is apparently due to Roquette.  A proof of
the key case ($X$ integral, affine) may be found in
(\Lcitemark 8\Rcitemark \ p.\ 39).

\begin{prop}\label{units-over-Z}
Let $X$ be a reduced scheme of finite type over $\Z$.  Then the group
$\Gamma(X,\O_X)^*$ is \FG.
\end{prop}

\begin{theorem}\label{rabbit-food}
Let $X$ be a reduced scheme of finite type over $\Z$.  Then the torsion
subgroup of $\Pic(X)$ is supported at a finite set of primes, and if
${}_p \Pic(X)$ is infinite, then the prime $p$ is not invertible in
$\Gamma(X,\O_X)$.
\end{theorem}

\begin{proof}
By (\Lcitemark 25\Rcitemark \ 2.7.6) we know that $\Pic(\nor{X})$ is \FG.
Therefore it suffices to show that if
$K = \Ker[\Pic(X)\ \mapE{}\ \Pic(\nor{X})]$, then $K$ is supported at a finite
set of primes, and if ${}_p K$ is infinite, then the prime $p$ is not
invertible in $\Gamma(X,\O_X)$.

By \pref{filtration}, \pref{units-over-Z} and the argument of
\pref{pic-main-theorem}, one sees that $K$ is isomorphic to the quotient of
an arithmetically linear group by a \FG\ subgroup.  Hence
\pref{arith-linear-over-finite} tells us that the torsion subgroup of $\Pic(X)$
is supported at a finite set of primes.  The last assertion of the theorem
follows from (\ref{pic-main-theorem}a).  \qed
\end{proof}

For $X$ a non-reduced scheme of finite type over $\Z$, we do not know if the
torsion subgroup of $\Pic(X)$ is supported at a finite set of primes.  However,
for any commutative noetherian ring $A$, the canonical map
\mapx[[ \Pic(A) || \Pic(\RED{A}) ]] is an isomorphism, so we have:

\begin{corollary}\label{yyyyy}
Let $A$ be a \FG\ commutative $\Z$-algebra.  Then the torsion subgroup of
$\Pic(A)$ is supported at a finite set of primes, and if ${}_p \Pic(A)$ is
infinite, then the prime $p$ is not invertible in $A$.
\end{corollary}

Let $A$ be a \FG\ commutative $\Z$-algebra.  If is natural to ask if there
exist prime numbers
$\vec p1n$, none of which are invertible in $A$, such that the subgroup of
$\Pic(A)$ consisting of torsion prime to $\vec p1n$ is finite.
Unfortunately, the answer is no.  For a counterexample, see
(\Lcitemark 21\Rcitemark \ 6.3).

Now we consider what happens when $X$ is of finite type over $\F_p$.

\begin{theorem}\label{F-sub-p}
Let $X$ be a scheme of finite type over $\F_p$.  Then there exists a
\FG\ abelian group $H$ and a finite $p$-group $F$ such that
$$\Pic(X)\ \iso\ H\ \o+ \ [\o+_{n=1}^\infty F].$$
\end{theorem}

\begin{sketch}
Consider the class $\shC$ of groups of the form ascribed to $\Pic(X)$ in the
theorem.  The groups in $\shC$ are all \bxc.  One can check without great
difficulty that $\shC$ is closed under formation of subgroups, quotient groups,
and extensions.

First suppose that $X$ is reduced.  Since $\Pic(\nor{X})$ is \FG, it suffices
to show that $K = \Ker[\Pic(X)\ \mapE{}\ \Pic(\nor{X})]$ is in $\shC$.  We may
assume that $X$ is affine.  In fact, it suffices to show that if
$X' = \Spec(A')$ is a partial normalization of $X = \Spec(A)$, and if
$[A':A]$ is prime, then $K' = \Ker[\Pic(X)\ \mapE{}\ \Pic(X')]$ is in $\shC$.
The proof of \pref{rabbit-food} shows that $K'$ is a quotient of an
arithmetically linear group.  In fact, the arguments used to arrive at this
result show that there exists an $\F_p$-algebra $A$ of finite type such that
$K'$ is a quotient of an abelian subgroup $H \IN \GL_n(A)$.  Therefore there
exists an $\F_p$-algebra $C$ of finite type such that $K'$ is a quotient of a
subgroup of $C^*$.  The usual methods show that $C^*$ may be built up via
extensions from a \FG\ abelian group and some $\F_p$-vector spaces.  It follows
that $C^* \in \shC$, and hence that any quotient of $C^*$ is in $\shC$.

Now suppose that $X$ is arbitrary, not necessarily reduced.  Let $\shJ$ be the
nilradical of $X$.  There is an exact sequence:
\Rowfive{0}{\shJ^n/\shJ^{n+1}}{(\O_X/\shJ^{n+1})^*}{(\O_X/\shJ^n)^*}{1%
}of sheaves of abelian groups on $X$.  Let $X_n$ be the closed subscheme of $X$
corresponding to the ideal $\shJ^n$.  By taking cohomology, one sees that
$$\Ker[\Pic(X_{n+1}) \ \mapE{}\ \Pic(X_n)]$%
$is $p$-torsion.  The theorem follows.  \qed
\end{sketch}

Finally, we compute the Picard group of a few simple examples.  The
calculations are an easy consequence of an exact sequence of Milnor
(\Lcitemark 6\Rcitemark \ IX\ 5.3) applied to the cartesian square
\squareSE{A}{\nor{A}}{A/\lfC}{\nor{A}/\lfC\nullcomma%
}where $\lfC$ is the conductor of $\nor{A}$ into $A$.  The exact sequence is:
\splitdiagram{A^*&\mapE{}&\nor{A}^* \times (A/\lfC)^*&\mapE{}%
&(\nor{A}/\lfC)^*}{\mapE{}&\Pic(A)&\mapE{}%
&\Pic(\nor{A}) \times \Pic(A/\lfC)&\mapE{}&\Pic(\nor{A}/\lfC).%
}Let $p$ be a prime number.  The first four rings given below may be
viewed as subrings of $\Z[t,x]$.

\vspace*{0.1in}

\begin{center}
\renewcommand{\arraystretch}{1.2}
\begin{tabular}{||c|c||} \hline
ring $A$ & structure of $\Pic(A)$ \\ \hline\hline
$\Z[t^2,t^3]$ & $\Z$ \\ \hline
$\Z[t^2, t^3, x]$ & free abelian of countably infinite rank \\ \hline
$\Z[pt, t^2, t^3]$ & $\F_p$ \\ \hline
$\Z[pt, t^2, t^3, x]$ & $\F_p$-vector space of countably infinite rank
     \\ \hline
$\Z[1/p, t^2, t^3]$ & $\Z[1/p]$ \\ \hline
$\F_p[t^2, t^3, x]$ & $\F_p$-vector space of countably infinite rank
     \\ \hline
\end{tabular}
\end{center}
\section*{References}
\addcontentsline{toc}{section}{References}
\ \par\noindent\vspace*{-0.25in}
\hfuzz 5pt

\bgroup\Resetstrings%
\def\Loccittest{}\def\Abbtest{ }\def\Capssmallcapstest{}\def\Edabbtest{
}\def\Edcapsmallcapstest{}\def\Underlinetest{ }%
\def\NoArev{1}\def\NoErev{0}\def\Acnt{1}\def\Ecnt{0}\def\acnt{0}\def\ecnt{0}%
\def\Ftest{ }\def\Fstr{1}%
\def\Atest{ }\def\Astr{Artin\Revcomma M\Initper }%
\def\Ttest{ }\def\Tstr{Algebraic approximation of structures over complete
local rings}%
\def\Jtest{ }\def\Jstr{Inst. Hautes \'Etudes Sci. Publ. Math.}%
\def\Vtest{ }\def\Vstr{36}%
\def\Dtest{ }\def\Dstr{1969}%
\def\Ptest{ }\def\Pcnt{ }\def\Pstr{23--58}%
\def\Qtest{ }\def\Qstr{access via "artin algebraic approximation"}%
\def\Xtest{ }\def\Xstr{In bound volume.}%
\Refformat\egroup%

\bgroup\Resetstrings%
\def\Loccittest{}\def\Abbtest{ }\def\Capssmallcapstest{}\def\Edabbtest{
}\def\Edcapsmallcapstest{}\def\Underlinetest{ }%
\def\NoArev{1}\def\NoErev{0}\def\Acnt{1}\def\Ecnt{0}\def\acnt{0}\def\ecnt{0}%
\def\Ftest{ }\def\Fstr{2}%
\def\Atest{ }\def\Astr{Artin\Revcomma M\Initper }%
\def\Ttest{ }\def\Tstr{Letter to Grothendieck}%
\def\Dtest{ }\def\Dstr{Nov.\ 5, 1968}%
\def\Qtest{ }\def\Qstr{access via "artin letter grothendieck"}%
\def\Astr{\Underlinemark}%
\Refformat\egroup%

\bgroup\Resetstrings%
\def\Loccittest{}\def\Abbtest{ }\def\Capssmallcapstest{}\def\Edabbtest{
}\def\Edcapsmallcapstest{}\def\Underlinetest{ }%
\def\NoArev{1}\def\NoErev{0}\def\Acnt{1}\def\Ecnt{0}\def\acnt{0}\def\ecnt{0}%
\def\Ftest{ }\def\Fstr{3}%
\def\Atest{ }\def\Astr{Artin\Revcomma M\Initper }%
\def\Ttest{ }\def\Tstr{Algebraization of formal moduli: I}%
\def\Btest{ }\def\Bstr{Global Analysis: Papers in Honor of K.\ Kodaira}%
\def\Itest{ }\def\Istr{Princeton Univ. Press}%
\def\Dtest{ }\def\Dstr{1969}%
\def\Ptest{ }\def\Pcnt{ }\def\Pstr{21--71}%
\def\Qtest{ }\def\Qstr{access via "artin formal moduli one"}%
\Refformat\egroup%

\bgroup\Resetstrings%
\def\Loccittest{}\def\Abbtest{ }\def\Capssmallcapstest{}\def\Edabbtest{
}\def\Edcapsmallcapstest{}\def\Underlinetest{ }%
\def\NoArev{1}\def\NoErev{0}\def\Acnt{2}\def\Ecnt{0}\def\acnt{0}\def\ecnt{0}%
\def\Ftest{ }\def\Fstr{4}%
\def\Atest{ }\def\Astr{Atiyah\Revcomma M\Initper \Initgap F\Initper %
  \Aand I\Initper \Initgap G\Initper  Macdonald}%
\def\Ttest{ }\def\Tstr{Introduction to Commutative Algebra}%
\def\Itest{ }\def\Istr{Addison-Wesley}%
\def\Ctest{ }\def\Cstr{Reading, Mass.}%
\def\Dtest{ }\def\Dstr{1969}%
\def\Qtest{ }\def\Qstr{access via "atiyah macdonald"}%
\Refformat\egroup%

\bgroup\Resetstrings%
\def\Loccittest{}\def\Abbtest{ }\def\Capssmallcapstest{}\def\Edabbtest{
}\def\Edcapsmallcapstest{}\def\Underlinetest{ }%
\def\NoArev{1}\def\NoErev{0}\def\Acnt{1}\def\Ecnt{4}\def\acnt{0}\def\ecnt{0}%
\def\Ftest{ }\def\Fstr{5}%
\def\Atest{ }\def\Astr{Auslander\Revcomma M\Initper }%
\def\Ttest{ }\def\Tstr{Coherent functors}%
\def\Btest{ }\def\Bstr{Proceedings of the Conference on Categorical Algebra (La
Jolla, 1965)}%
\def\Etest{ }\def\Estr{S\Initper  Eilenberg%
  \Ecomma D\Initper \Initgap K\Initper  Harrison%
  \Ecomma S\Initper  MacLane%
  \Eandd H\Initper  R\"ohrl}%
\def\Ptest{ }\def\Pcnt{ }\def\Pstr{189--231}%
\def\Itest{ }\def\Istr{Spring\-er-Ver\-lag}%
\def\Ctest{ }\def\Cstr{New York}%
\def\Qtest{ }\def\Qstr{access via "auslander coherent functors"}%
\Refformat\egroup%

\bgroup\Resetstrings%
\def\Loccittest{}\def\Abbtest{ }\def\Capssmallcapstest{}\def\Edabbtest{
}\def\Edcapsmallcapstest{}\def\Underlinetest{ }%
\def\NoArev{1}\def\NoErev{0}\def\Acnt{1}\def\Ecnt{0}\def\acnt{0}\def\ecnt{0}%
\def\Ftest{ }\def\Fstr{6}%
\def\Atest{ }\def\Astr{Bass\Revcomma H\Initper }%
\def\Ttest{ }\def\Tstr{Algebraic K-Theory}%
\def\Itest{ }\def\Istr{W.\ A.\ Benjamin}%
\def\Ctest{ }\def\Cstr{New York}%
\def\Dtest{ }\def\Dstr{1968}%
\Refformat\egroup%

\bgroup\Resetstrings%
\def\Loccittest{}\def\Abbtest{ }\def\Capssmallcapstest{}\def\Edabbtest{
}\def\Edcapsmallcapstest{}\def\Underlinetest{ }%
\def\NoArev{1}\def\NoErev{0}\def\Acnt{1}\def\Ecnt{0}\def\acnt{0}\def\ecnt{0}%
\def\Ftest{ }\def\Fstr{7}%
\def\Atest{ }\def\Astr{Bass\Revcomma H\Initper }%
\def\Ttest{ }\def\Tstr{Some problems in ``classical'' algebraic K-theory}%
\def\Stest{ }\def\Sstr{Lecture Notes in Mathematics}%
\def\Itest{ }\def\Istr{Spring\-er-Ver\-lag}%
\def\Ctest{ }\def\Cstr{New York}%
\def\Vtest{ }\def\Vstr{342}%
\def\Dtest{ }\def\Dstr{1973}%
\def\Ptest{ }\def\Pcnt{ }\def\Pstr{3--73}%
\def\Qtest{ }\def\Qstr{access via "bass classical problems"}%
\def\Xtest{ }\def\Xstr{I don't have this.}%
\def\Astr{\Underlinemark}%
\Refformat\egroup%

\bgroup\Resetstrings%
\def\Loccittest{}\def\Abbtest{ }\def\Capssmallcapstest{}\def\Edabbtest{
}\def\Edcapsmallcapstest{}\def\Underlinetest{ }%
\def\NoArev{1}\def\NoErev{0}\def\Acnt{1}\def\Ecnt{0}\def\acnt{0}\def\ecnt{0}%
\def\Ftest{ }\def\Fstr{8}%
\def\Atest{ }\def\Astr{Bass\Revcomma H\Initper }%
\def\Ttest{ }\def\Tstr{Introduction to some methods of algebraic K-theory}%
\def\Jtest{ }\def\Jstr{Conference Board of the Mathematical Sciences Regional
Conference Series in Mathematics}%
\def\Vtest{ }\def\Vstr{20}%
\def\Dtest{ }\def\Dstr{1974}%
\def\Qtest{ }\def\Qstr{access via "bass introduction methods"}%
\def\Xtest{ }\def\Xstr{MR 50:441.  (pamphlet on shelf) \par Let $R$ be either
$\Z$ or $F_q[t]$.  Let $A = R[\vec t2d]$, where $d \geq 1$.  Then $\GL_n(A)$ is
a finitely generated group for  $n \geq d+2$.  Also, let $A$ be any semilocal
ring or $\Z$ or $k[t]$, where $k$ is any field.  Then $\SL_n(A)$ is generated
by elementary matrices for any $n \geq 1$.  (See 2.6, 2.8.)}%
\def\Astr{\Underlinemark}%
\Refformat\egroup%

\bgroup\Resetstrings%
\def\Loccittest{}\def\Abbtest{ }\def\Capssmallcapstest{}\def\Edabbtest{
}\def\Edcapsmallcapstest{}\def\Underlinetest{ }%
\def\NoArev{1}\def\NoErev{0}\def\Acnt{1}\def\Ecnt{0}\def\acnt{0}\def\ecnt{0}%
\def\Ftest{ }\def\Fstr{9}%
\def\Atest{ }\def\Astr{Bergman\Revcomma G\Initper \Initgap M\Initper }%
\def\Ttest{ }\def\Tstr{Ring schemes; the Witt scheme}%
\def\Btest{ }\def\Bstr{Lectures on Curves on an Algebraic Surface {\rm by David
Mumford}}%
\def\Ptest{ }\def\Pcnt{ }\def\Pstr{171---187}%
\def\Itest{ }\def\Istr{Princeton Univ.\ Press}%
\def\Dtest{ }\def\Dstr{1966}%
\def\Qtest{ }\def\Qstr{access via "bergman witt scheme"}%
\Refformat\egroup%

\bgroup\Resetstrings%
\def\Loccittest{}\def\Abbtest{ }\def\Capssmallcapstest{}\def\Edabbtest{
}\def\Edcapsmallcapstest{}\def\Underlinetest{ }%
\def\NoArev{1}\def\NoErev{0}\def\Acnt{1}\def\Ecnt{0}\def\acnt{0}\def\ecnt{0}%
\def\Ftest{ }\def\Fstr{10}%
\def\Atest{ }\def\Astr{Bertin\Revcomma J\Initper \Initgap E\Initper }%
\def\Ttest{ }\def\Tstr{{\rm\tolerance=1000 Generalites sur les preschemas en
groupes,  \EXPOSE\ ${\rm VI}_{\rm B}$ in {\it\SGA} (SGA 3)}}%
\def\Stest{ }\def\Sstr{Lecture Notes in Mathematics}%
\def\Itest{ }\def\Istr{Spring\-er-Ver\-lag}%
\def\Ctest{ }\def\Cstr{New York}%
\def\Vtest{ }\def\Vstr{151}%
\def\Dtest{ }\def\Dstr{1970}%
\def\Ptest{ }\def\Pcnt{ }\def\Pstr{318--410}%
\def\Qtest{ }\def\Qstr{access via "bertin preschemas"}%
\def\Xtest{ }\def\Xstr{This includes the following result, announced by
Raynaud, but apparently not proved in this volume:  Theorem (11.11.1)  Let $S$
be a regular noetherian scheme of dimension $/leq 1$.  Let \map(\pi,G,S) be an
$S$-group scheme of finite-type.  Assume that $\pi$ is a flat, affine morphism.
 Then $G$ is isomorphic to a closed  subgroup of $\AUT(V)$, where $V$ is a
vector bundle on $S$.}%
\Refformat\egroup%

\bgroup\Resetstrings%
\def\Loccittest{}\def\Abbtest{ }\def\Capssmallcapstest{}\def\Edabbtest{
}\def\Edcapsmallcapstest{}\def\Underlinetest{ }%
\def\NoArev{1}\def\NoErev{0}\def\Acnt{1}\def\Ecnt{0}\def\acnt{0}\def\ecnt{0}%
\def\Ftest{ }\def\Fstr{11}%
\def\Atest{ }\def\Astr{Bourbaki\Revcomma N\Initper }%
\def\Ttest{ }\def\Tstr{Elements of Mathematics (Algebra II, Chapters 4--7)}%
\def\Itest{ }\def\Istr{Spring\-er-Ver\-lag}%
\def\Ctest{ }\def\Cstr{New York}%
\def\Dtest{ }\def\Dstr{1990}%
\def\Qtest{ }\def\Qstr{access via "bourbaki algebra part two"}%
\def\Xtest{ }\def\Xstr{I don't have this.}%
\Refformat\egroup%

\bgroup\Resetstrings%
\def\Loccittest{}\def\Abbtest{ }\def\Capssmallcapstest{}\def\Edabbtest{
}\def\Edcapsmallcapstest{}\def\Underlinetest{ }%
\def\NoArev{1}\def\NoErev{0}\def\Acnt{1}\def\Ecnt{0}\def\acnt{0}\def\ecnt{0}%
\def\Ftest{ }\def\Fstr{12}%
\def\Atest{ }\def\Astr{Claborn\Revcomma L\Initper }%
\def\Ttest{ }\def\Tstr{Every abelian group is a class group}%
\def\Jtest{ }\def\Jstr{Pacific J. Math.}%
\def\Vtest{ }\def\Vstr{18}%
\def\Dtest{ }\def\Dstr{1966}%
\def\Ptest{ }\def\Pcnt{ }\def\Pstr{219--222}%
\def\Qtest{ }\def\Qstr{access via "claborn"}%
\def\Xtest{ }\def\Xstr{I don't have this.}%
\Refformat\egroup%

\bgroup\Resetstrings%
\def\Loccittest{}\def\Abbtest{ }\def\Capssmallcapstest{}\def\Edabbtest{
}\def\Edcapsmallcapstest{}\def\Underlinetest{ }%
\def\NoArev{1}\def\NoErev{0}\def\Acnt{1}\def\Ecnt{0}\def\acnt{0}\def\ecnt{0}%
\def\Ftest{ }\def\Fstr{13}%
\def\Atest{ }\def\Astr{Cohen\Revcomma H\Initper }%
\def\Ttest{ }\def\Tstr{Un faisceau qui ne peut pas \^etre d\'etordu
universellement}%
\def\Jtest{ }\def\Jstr{C. R. Acad. Sci. Paris S\'er. I Math.}%
\def\Vtest{ }\def\Vstr{272}%
\def\Dtest{ }\def\Dstr{1971}%
\def\Ptest{ }\def\Pcnt{ }\def\Pstr{799--802}%
\def\Qtest{ }\def\Qstr{access via "cohen faisceau article"}%
\Refformat\egroup%

\bgroup\Resetstrings%
\def\Loccittest{}\def\Abbtest{ }\def\Capssmallcapstest{}\def\Edabbtest{
}\def\Edcapsmallcapstest{}\def\Underlinetest{ }%
\def\NoArev{1}\def\NoErev{0}\def\Acnt{1}\def\Ecnt{0}\def\acnt{0}\def\ecnt{0}%
\def\Ftest{ }\def\Fstr{14}%
\def\Atest{ }\def\Astr{Cohen\Revcomma H\Initper }%
\def\Ttest{ }\def\Tstr{Detorsion universelle de faisceaux coherents}%
\def\otest{ }\def\ostr{thesis (Docteur $3^\circ$ Cycle)}%
\def\Dtest{ }\def\Dstr{1972}%
\def\Itest{ }\def\Istr{Universit\'e de Paris}%
\def\Ctest{ }\def\Cstr{Orsay}%
\def\Qtest{ }\def\Qstr{access via "cohen thesis"}%
\def\Astr{\Underlinemark}%
\Refformat\egroup%

\bgroup\Resetstrings%
\def\Loccittest{}\def\Abbtest{ }\def\Capssmallcapstest{}\def\Edabbtest{
}\def\Edcapsmallcapstest{}\def\Underlinetest{ }%
\def\NoArev{1}\def\NoErev{0}\def\Acnt{1}\def\Ecnt{0}\def\acnt{0}\def\ecnt{0}%
\def\Ftest{ }\def\Fstr{15}%
\def\Atest{ }\def\Astr{Fossum\Revcomma R\Initper \Initgap M\Initper }%
\def\Ttest{ }\def\Tstr{The Divisor Class Group of a Krull Domain}%
\def\Itest{ }\def\Istr{Spring\-er-Ver\-lag}%
\def\Ctest{ }\def\Cstr{New York}%
\def\Dtest{ }\def\Dstr{1973}%
\def\Qtest{ }\def\Qstr{access via "fossum"}%
\def\Xtest{ }\def\Xstr{I don't have this.}%
\Refformat\egroup%

\bgroup\Resetstrings%
\def\Loccittest{}\def\Abbtest{ }\def\Capssmallcapstest{}\def\Edabbtest{
}\def\Edcapsmallcapstest{}\def\Underlinetest{ }%
\def\NoArev{1}\def\NoErev{0}\def\Acnt{1}\def\Ecnt{0}\def\acnt{0}\def\ecnt{0}%
\def\Ftest{ }\def\Fstr{16}%
\def\Atest{ }\def\Astr{Fuchs\Revcomma L\Initper }%
\def\Ttest{ }\def\Tstr{Abelian Groups}%
\def\Itest{ }\def\Istr{Hungarian Academy of Sciences}%
\def\Ctest{ }\def\Cstr{Budapest}%
\def\Dtest{ }\def\Dstr{1958}%
\def\Qtest{ }\def\Qstr{access via "fuchs abelian groups 1958"}%
\def\Xtest{ }\def\Xstr{I don't have this.}%
\Refformat\egroup%

\bgroup\Resetstrings%
\def\Loccittest{}\def\Abbtest{ }\def\Capssmallcapstest{}\def\Edabbtest{
}\def\Edcapsmallcapstest{}\def\Underlinetest{ }%
\def\NoArev{1}\def\NoErev{0}\def\Acnt{2}\def\Ecnt{0}\def\acnt{0}\def\ecnt{0}%
\def\Ftest{ }\def\Fstr{17}%
\def\Atest{ }\def\Astr{Grothendieck\Revcomma A\Initper %
  \Aand J\Initper \Initgap A\Initper  Dieudonn\'e}%
\def\Ttest{ }\def\Tstr{El\'ements de \geom\ alg\'e\-brique III (part one)}%
\def\Jtest{ }\def\Jstr{Inst. Hautes \'Etudes Sci. Publ. Math.}%
\def\Vtest{ }\def\Vstr{11}%
\def\Dtest{ }\def\Dstr{1961}%
\def\Qtest{ }\def\Qstr{access via "EGA3-1"}%
\Refformat\egroup%

\bgroup\Resetstrings%
\def\Loccittest{}\def\Abbtest{ }\def\Capssmallcapstest{}\def\Edabbtest{
}\def\Edcapsmallcapstest{}\def\Underlinetest{ }%
\def\NoArev{1}\def\NoErev{0}\def\Acnt{2}\def\Ecnt{0}\def\acnt{0}\def\ecnt{0}%
\def\Ftest{ }\def\Fstr{18}%
\def\Atest{ }\def\Astr{Grothendieck\Revcomma A\Initper %
  \Aand J\Initper \Initgap A\Initper  Dieudonn\'e}%
\def\Ttest{ }\def\Tstr{El\'ements de \geom\ alg\'e\-brique III (part two)}%
\def\Jtest{ }\def\Jstr{Inst. Hautes \'Etudes Sci. Publ. Math.}%
\def\Vtest{ }\def\Vstr{17}%
\def\Dtest{ }\def\Dstr{1963}%
\def\Qtest{ }\def\Qstr{access via "EGA3-2"}%
\def\Astr{\Underlinemark}%
\Refformat\egroup%

\bgroup\Resetstrings%
\def\Loccittest{}\def\Abbtest{ }\def\Capssmallcapstest{}\def\Edabbtest{
}\def\Edcapsmallcapstest{}\def\Underlinetest{ }%
\def\NoArev{1}\def\NoErev{0}\def\Acnt{2}\def\Ecnt{0}\def\acnt{0}\def\ecnt{0}%
\def\Ftest{ }\def\Fstr{19}%
\def\Atest{ }\def\Astr{Grothendieck\Revcomma A\Initper %
  \Aand J\Initper \Initgap A\Initper  Dieudonn\'e}%
\def\Ttest{ }\def\Tstr{El\'ements de \geom\ alg\'e\-brique IV (part three)}%
\def\Jtest{ }\def\Jstr{Inst. Hautes \'Etudes Sci. Publ. Math.}%
\def\Vtest{ }\def\Vstr{28}%
\def\Dtest{ }\def\Dstr{1966}%
\def\Qtest{ }\def\Qstr{access via "EGA4-3"}%
\def\Xtest{ }\def\Xstr{bertini: 9.7.7 A proper monomorphism of finite
presentation is a closed immersion (8.11.5).}%
\def\Astr{\Underlinemark}%
\Refformat\egroup%

\bgroup\Resetstrings%
\def\Loccittest{}\def\Abbtest{ }\def\Capssmallcapstest{}\def\Edabbtest{
}\def\Edcapsmallcapstest{}\def\Underlinetest{ }%
\def\NoArev{1}\def\NoErev{0}\def\Acnt{2}\def\Ecnt{0}\def\acnt{0}\def\ecnt{0}%
\def\Ftest{ }\def\Fstr{20}%
\def\Atest{ }\def\Astr{Grothendieck\Revcomma A\Initper %
  \Aand J\Initper \Initgap A\Initper  Dieudonn\'e}%
\def\Ttest{ }\def\Tstr{El\'ements de \geom\ alg\'e\-brique I}%
\def\Itest{ }\def\Istr{Springer-Verlag}%
\def\Ctest{ }\def\Cstr{New \snap York}%
\def\Dtest{ }\def\Dstr{1971}%
\def\Qtest{ }\def\Qstr{access via "EGA1"}%
\def\Astr{\Underlinemark}%
\Refformat\egroup%

\bgroup\Resetstrings%
\def\Loccittest{}\def\Abbtest{ }\def\Capssmallcapstest{}\def\Edabbtest{
}\def\Edcapsmallcapstest{}\def\Underlinetest{ }%
\def\NoArev{1}\def\NoErev{0}\def\Acnt{4}\def\Ecnt{0}\def\acnt{0}\def\ecnt{0}%
\def\Ftest{ }\def\Fstr{21}%
\def\Atest{ }\def\Astr{Guralnick\Revcomma R\Initper %
  \Acomma D\Initper \Initgap B\Initper  Jaffe%
  \Acomma W\Initper  Raskind%
  \Aandd R\Initper  Wiegand}%
\def\Ttest{ }\def\Tstr{The kernel of the map on Picard groups induced by a
faithfully flat homomorphism}%
\def\Rtest{ }\def\Rstr{preprint}%
\def\Qtest{ }\def\Qstr{access via "gang of four"}%
\Refformat\egroup%

\bgroup\Resetstrings%
\def\Loccittest{}\def\Abbtest{ }\def\Capssmallcapstest{}\def\Edabbtest{
}\def\Edcapsmallcapstest{}\def\Underlinetest{ }%
\def\NoArev{1}\def\NoErev{0}\def\Acnt{1}\def\Ecnt{0}\def\acnt{0}\def\ecnt{0}%
\def\Ftest{ }\def\Fstr{22}%
\def\Atest{ }\def\Astr{Hartshorne\Revcomma R\Initper }%
\def\Ttest{ }\def\Tstr{Algebraic Geometry}%
\def\Itest{ }\def\Istr{Spring\-er-Ver\-lag}%
\def\Ctest{ }\def\Cstr{New York}%
\def\Dtest{ }\def\Dstr{1977}%
\def\Qtest{ }\def\Qstr{access via "hartshorne algebraic geometry"}%
\Refformat\egroup%

\bgroup\Resetstrings%
\def\Loccittest{}\def\Abbtest{ }\def\Capssmallcapstest{}\def\Edabbtest{
}\def\Edcapsmallcapstest{}\def\Underlinetest{ }%
\def\NoArev{1}\def\NoErev{0}\def\Acnt{1}\def\Ecnt{0}\def\acnt{0}\def\ecnt{0}%
\def\Ftest{ }\def\Fstr{23}%
\def\Atest{ }\def\Astr{Illusie\Revcomma L\Initper }%
\def\Ttest{ }\def\Tstr{G\'en\'eralit\'es sur les conditions de finitude dans
les cat\'egories d\'eriv\'ees, \EXPOSE\ I in {\it\SGA} (SGA) 6}%
\def\Stest{ }\def\Sstr{Lecture Notes in Mathematics}%
\def\Itest{ }\def\Istr{Spring\-er-Ver\-lag}%
\def\Ctest{ }\def\Cstr{New York}%
\def\Vtest{ }\def\Vstr{225}%
\def\Dtest{ }\def\Dstr{1971}%
\def\Ptest{ }\def\Pcnt{ }\def\Pstr{78--159}%
\def\Qtest{ }\def\Qstr{access via "illusie finiteness first"}%
\Refformat\egroup%

\bgroup\Resetstrings%
\def\Loccittest{}\def\Abbtest{ }\def\Capssmallcapstest{}\def\Edabbtest{
}\def\Edcapsmallcapstest{}\def\Underlinetest{ }%
\def\NoArev{1}\def\NoErev{0}\def\Acnt{1}\def\Ecnt{0}\def\acnt{0}\def\ecnt{0}%
\def\Ftest{ }\def\Fstr{24}%
\def\Atest{ }\def\Astr{Koblitz\Revcomma N\Initper }%
\def\Ttest{ }\def\Tstr{$p$-adic Numbers, $p$-adic Analysis, and
Zeta-Functions}%
\def\Itest{ }\def\Istr{\snap\hbox{Springer-Verlag}}%
\def\Ctest{ }\def\Cstr{New York}%
\def\Dtest{ }\def\Dstr{1977}%
\def\Qtest{ }\def\Qstr{access via "koblitz p-adic analysis"}%
\Refformat\egroup%

\bgroup\Resetstrings%
\def\Loccittest{}\def\Abbtest{ }\def\Capssmallcapstest{}\def\Edabbtest{
}\def\Edcapsmallcapstest{}\def\Underlinetest{ }%
\def\NoArev{1}\def\NoErev{0}\def\Acnt{1}\def\Ecnt{0}\def\acnt{0}\def\ecnt{0}%
\def\Ftest{ }\def\Fstr{25}%
\def\Atest{ }\def\Astr{Lang\Revcomma S\Initper }%
\def\Ttest{ }\def\Tstr{Fundamentals\kern2.5pt\ of\kern2.5pt\
Diophantine\kern2.5pt\ Geometry}%
\def\Itest{ }\def\Istr{Spring\-er-Ver\-lag}%
\def\Ctest{ }\def\Cstr{New York}%
\def\Dtest{ }\def\Dstr{1983}%
\def\Qtest{ }\def\Qstr{access via "lang diophantine geometry"}%
\Refformat\egroup%

\bgroup\Resetstrings%
\def\Loccittest{}\def\Abbtest{ }\def\Capssmallcapstest{}\def\Edabbtest{
}\def\Edcapsmallcapstest{}\def\Underlinetest{ }%
\def\NoArev{1}\def\NoErev{0}\def\Acnt{1}\def\Ecnt{0}\def\acnt{0}\def\ecnt{0}%
\def\Ftest{ }\def\Fstr{26}%
\def\Atest{ }\def\Astr{Lang\Revcomma S\Initper }%
\def\Ttest{ }\def\Tstr{Algebra ($2^{\rm nd}$ edition)}%
\def\Itest{ }\def\Istr{Addison-Wesley}%
\def\Ctest{ }\def\Cstr{Reading, Mass.}%
\def\Dtest{ }\def\Dstr{1984}%
\def\Qtest{ }\def\Qstr{access via "lang algebra bible second edition"}%
\def\Xtest{ }\def\Xstr{I don't have this but Roger Wiegand does.}%
\def\Astr{\Underlinemark}%
\Refformat\egroup%

\bgroup\Resetstrings%
\def\Loccittest{}\def\Abbtest{ }\def\Capssmallcapstest{}\def\Edabbtest{
}\def\Edcapsmallcapstest{}\def\Underlinetest{ }%
\def\NoArev{1}\def\NoErev{0}\def\Acnt{1}\def\Ecnt{0}\def\acnt{0}\def\ecnt{0}%
\def\Ftest{ }\def\Fstr{27}%
\def\MacLane{Mac Lane}
{}%
\def\Atest{ }\def\Astr{\MacLane\Revcomma S\Initper }%
\def\Ttest{ }\def\Tstr{Categories for the Working Mathematician}%
\def\Itest{ }\def\Istr{Spring\-er-Ver\-lag}%
\def\Ctest{ }\def\Cstr{New York}%
\def\Dtest{ }\def\Dstr{1971}%
\def\Qtest{ }\def\Qstr{access via "working mathematician"}%
\Refformat\egroup%

\bgroup\Resetstrings%
\def\Loccittest{}\def\Abbtest{ }\def\Capssmallcapstest{}\def\Edabbtest{
}\def\Edcapsmallcapstest{}\def\Underlinetest{ }%
\def\NoArev{1}\def\NoErev{0}\def\Acnt{1}\def\Ecnt{0}\def\acnt{0}\def\ecnt{0}%
\def\Ftest{ }\def\Fstr{28}%
\def\Atest{ }\def\Astr{Rotman\Revcomma J\Initper \Initgap J\Initper }%
\def\Ttest{ }\def\Tstr{An Introduction to Homological Algebra}%
\def\Itest{ }\def\Istr{Academic Press}%
\def\Ctest{ }\def\Cstr{New York}%
\def\Dtest{ }\def\Dstr{1979}%
\def\Qtest{ }\def\Qstr{access via "rotman homological algebra"}%
\Refformat\egroup%

\bgroup\Resetstrings%
\def\Loccittest{}\def\Abbtest{ }\def\Capssmallcapstest{}\def\Edabbtest{
}\def\Edcapsmallcapstest{}\def\Underlinetest{ }%
\def\NoArev{1}\def\NoErev{0}\def\Acnt{1}\def\Ecnt{0}\def\acnt{0}\def\ecnt{0}%
\def\Ftest{ }\def\Fstr{29}%
\def\Atest{ }\def\Astr{Serre\Revcomma J-P\Initper }%
\def\Ttest{ }\def\Tstr{Local Fields}%
\def\Itest{ }\def\Istr{Spring\-er-Ver\-lag}%
\def\Ctest{ }\def\Cstr{New York}%
\def\Dtest{ }\def\Dstr{1979}%
\def\Qtest{ }\def\Qstr{access via "serre local fields"}%
\Refformat\egroup%

\end{document}